\documentclass[10pt]{article}
\usepackage[a4paper,margin=25mm,top=24mm,bottom=26mm]{geometry}
\usepackage{fontspec}
\usepackage{unicode-math}
\usepackage{amsmath}
\usepackage{xcolor}
\usepackage{fancyvrb}
\usepackage{longtable}
\usepackage{array}
\usepackage{colortbl}
\usepackage[most]{tcolorbox}
\usepackage{titlesec}
\usepackage{fancyhdr}
\usepackage{hyperref}

\newfontfamily\dvsans{DejaVuSans.ttf}[Path=fonts/]

\definecolor{hdr}{HTML}{2B4A6F}
\definecolor{band}{HTML}{DDE6F0}
\definecolor{hl}{HTML}{FDF0C8}
\definecolor{hl2}{HTML}{F2F6EA}
\definecolor{grey}{HTML}{7A7F87}
\definecolor{rule}{HTML}{9AA6B4}
\definecolor{qbg}{HTML}{F4F7FA}
\definecolor{gC}{HTML}{1F6F43}
\definecolor{gN}{HTML}{8A5A00}
\definecolor{gS}{HTML}{3A4A8A}
\definecolor{gO}{HTML}{6B3A7A}

\newcommand{\vbx}[1]{\hbox to \fontcharwd\font`0 {\hss{\dvsans\char"#1}\hss}}

\newcommand{\ctwo}{\textup{sel}}
\newcommand{\cthree}{\textup{rig}}
\newcommand{\fmark}{{\dvsans\char"25CF}}
\newcommand{\xmark}{{\dvsans\char"00D7}}

\begingroup
\catcode`\②=\active \catcode`\③=\active \catcode`\●=\active \catcode`\⟺=\active
\catcode`\₀=\active \catcode`\₁=\active \catcode`\₂=\active \catcode`\₃=\active
\catcode`\₄=\active \catcode`\₅=\active \catcode`\₆=\active \catcode`\₇=\active
\catcode`\₈=\active \catcode`\₉=\active \catcode`\ᵀ=\active \catcode`\ⁱ=\active
\catcode`\ξ=\active
\gdef②{\vbx{2461}}\gdef③{\vbx{2462}}\gdef●{\vbx{25CF}}\gdef⟺{\vbx{27FA}}
\gdef₀{\vbx{2080}}\gdef₁{\vbx{2081}}\gdef₂{\vbx{2082}}\gdef₃{\vbx{2083}}
\gdef₄{\vbx{2084}}\gdef₅{\vbx{2085}}\gdef₆{\vbx{2086}}\gdef₇{\vbx{2087}}
\gdef₈{\vbx{2088}}\gdef₉{\vbx{2089}}\gdefᵀ{\vbx{1D40}}\gdefⁱ{\vbx{2071}}
\gdefξ{\vbx{03BE}}
\endgroup
\newcommand{\vbcodes}{%
  \catcode`\②=\active \catcode`\③=\active \catcode`\●=\active \catcode`\⟺=\active
  \catcode`\₀=\active \catcode`\₁=\active \catcode`\₂=\active \catcode`\₃=\active
  \catcode`\₄=\active \catcode`\₅=\active \catcode`\₆=\active \catcode`\₇=\active
  \catcode`\₈=\active \catcode`\₉=\active \catcode`\ᵀ=\active \catcode`\ⁱ=\active
  \catcode`\ξ=\active}

\newenvironment{codebox}[2]{\VerbatimEnvironment\begin{Verbatim}[fontsize=\fontsize{#1pt}{#2pt}\selectfont,
  frame=leftline,framerule=1.1pt,rulecolor=\color{rule},framesep=7pt,xleftmargin=3mm,samepage=false,codes=\vbcodes]}%
  {\end{Verbatim}}

\newenvironment{codeboxk}[2]{\VerbatimEnvironment\begin{Verbatim}[fontsize=\fontsize{#1pt}{#2pt}\selectfont,
  frame=leftline,framerule=1.1pt,rulecolor=\color{rule},framesep=7pt,xleftmargin=3mm,samepage=true,codes=\vbcodes]}%
  {\end{Verbatim}}

\newtcolorbox{qbox}{enhanced,breakable,boxrule=0pt,frame hidden,colback=qbg,
  borderline west={1.6pt}{0pt}{hdr},left=7pt,right=7pt,top=5pt,bottom=5pt,parbox=false}

\newcommand{\grade}[2]{{\footnotesize\color{#1}[\,#2\,]}}
\newcommand{\ct}[1]{{\fontsize{5.6}{6.4}\selectfont\ttfamily #1}}
\newcommand{\shd}[1]{{\fontsize{6.4}{7.4}\selectfont #1}}

\titleformat{\section}{\Large\bfseries\color{hdr}}{\thesection}{0.7em}{}
\titleformat{\subsection}{\large\bfseries\color{hdr}}{\thesubsection}{0.7em}{}
\newcommand{\sect}[2]{\section*{#1}\addcontentsline{toc}{section}{\texorpdfstring{#1}{#2}}}
\newcommand{\subsect}[2]{\subsection*{#1}\addcontentsline{toc}{subsection}{\texorpdfstring{#1}{#2}}}

\newcommand{\chainband}[4]{\noalign{\vspace{1.5pt}}%
  \multicolumn{9}{@{}>{\columncolor{band}}l@{}}{\rule{0pt}{2.4ex}\hspace{1mm}\textbf{#1 #2}\quad\textperiodcentered\quad $D = #3$\quad\textperiodcentered\quad #4}\\
  \noalign{\vspace{1.5pt}}}

\newcommand{\chainbandD}[2]{\noalign{\vspace{1.5pt}}%
  \multicolumn{6}{@{}>{\columncolor{band}}l@{}}{\rule{0pt}{2.4ex}\hspace{1mm}\textbf{#1}\quad\textperiodcentered\quad #2}\\
  \noalign{\vspace{1.5pt}}}

\newenvironment{censustable}[3]{%
  \par\vspace{4mm}{\large\bfseries\color{hdr}#1}\quad{\footnotesize\color{grey}#2}\par
  {\footnotesize #3\par}\vspace{2mm}%
  \setlength{\LTleft}{0pt}\setlength{\LTright}{0pt plus 1fil}%
  \renewcommand{\arraystretch}{1.18}\setlength{\tabcolsep}{2pt}%
  \footnotesize
  \begin{longtable}{@{}p{41mm} r r r r r c c >{\raggedright\arraybackslash}p{73mm}@{}}}%
  {\end{longtable}}

\hypersetup{colorlinks=true,linkcolor=hdr,urlcolor=hdr,citecolor=hdr,bookmarks=false}

\begin{document}
\thispagestyle{empty}
\begin{center}
{\LARGE\bfseries Petersson-Rigid Lattices in a Census of 100 Rank-Three Root Bases}\\[4mm]
{\large Weakly harmonic Maass forms from rank-three Lorentzian Kac–Moody Cartan matrices}\\[6mm]
{Eungang Cho}\\[1mm]{\small Department of Mathematics, Chung-Ang University}\\[4mm]
{\footnotesize 2026-09-02. Sequel to arXiv:2608.19706 (“paper I”).}\\{\footnotesize Reproduction archive: Zenodo, DOI 10.5281/zenodo.22257942 (archive v3.5).}
\end{center}
\vspace{4mm}
\begin{center}{\large\bfseries Abstract}\end{center}
{\small
The first paper worked out one family of four lattices in full — $\lvert\det\rvert\, = 12, 24, 36$ and $72$. This paper generalizes it. We enumerate the 44 symmetrizable rank-three hyperbolic Cartan matrices and their 56 depth-one edits, 100 root bases in all, and compute the obstruction space $S_{5/2}(\rho _L)$ of the 98 that lie within our weight-$5/2$ budget: 10 vacuous, 34 unobstructed, 54 obstructed \grade{gC}{C}.

The vacuous ten are the lattices Bruinier, Ehlen and Freitag call \textbf{simple}, read in signature $(2,3)$. They classified them — there are 15 — and our ten realize five of the fifteen; of the other ten, five need more than three generators and cannot be the discriminant form of a rank-three lattice at all, three fail $|\operatorname{det} L| = 2k^{2}$, and two are simply not root bases. Within the rank-three hyperbolic world the simple lattices are exactly the Feingold–Frenkel neighbours of index $k \leq  4$.

The discriminant group $L'/L$ carries a finite quadratic form, and the finite group of its isometries acts on the weight-$3/2$ cusp forms for $\rho _L$, the bottom antisymmetric rung of the weight tower. We survey the lattices on which that action is absolutely irreducible of dimension at least two, so that the Petersson pairing on that space is pinned down up to a single scalar. The condition alone cuts the census to three: $L_{4}$ of the first paper, and two new ones at $\lvert\det\rvert\, 40$ and $88$. The three quaternion discriminants that occur are $6$, $10$ and $22$, the three for which the Shimura curve $X^D$ has genus zero. Rigidity uses no quaternion input, so we record this as an observation about the census and not as a characterisation (\S{}5.2).

Both invariants of the first paper — the shadow norm $\lVert \Xi \rVert ^{2}$ and the Petersson scalar $t$ — were single points there. Three rigid lattices instead of one are where their special faces come off: the irrational part of each generalizes, $\lVert \Xi \rVert ^{2}$ against $L(f,1)$ where the first paper read $L(f,2)$, and $t$ against an elliptic curve's imaginary period where it read $\Gamma (1/3)$. Both numerically, to between $26$ and $58$ digits (\S{}\S{}5–6).

}
\vspace{4mm}\tableofcontents\newpage
\sect{1. Introduction}{1. Introduction}

\subsect{1.1 The thesis}{1.1 The thesis}

The world in which the Borcherds correction succeeds has been mapped. Gritsenko and Nikulin classified the lattice Weyl vectors in the family $2U \oplus  \langle -2t\rangle $, Allcock enumerated the reflective lattices of small rank, and the resulting picture — which Lorentzian lattices carry an automorphic product with the expected divisor — is by now a finite, published list.

This paper counts the complement. It asks what happens when the correction fails, and finds that the failure is not a defect but an invariant: a weight-$5/2$ cusp form, canonically attached to the Cartan matrix, whose Hecke decomposition and Petersson norm are arithmetic data of the same kind that the successful cases suppress.

The organising discovery is that failure is \textbf{rare and concentrated}. Of the eleven compact-hyperbolic rank-three Cartan matrices — the ones every proper node-deletion of which is of finite type — exactly one has a nonvanishing shadow. And of the 100 root bases, exactly four are \textbf{Petersson-rigid} — three lattices, since $L_{4}$ is counted through each of its two root bases (\S{}4.7): the isometries of the discriminant form act absolutely irreducibly, in dimension at least two, on the weight-$3/2$ space one rung below the shadow. Those three sit on the three genus-zero Shimura curves — all three of them, and nothing else. Rigidity is a condition on that lower rung and not on the shadow itself, which at one of the three is not even pure (\S{}5.1, \S{}6.1).

\subsect{1.2 Results}{1.2 Results}

\textbf{Theorem A (\S{}3.4).} Among the 98 computed root bases, the weight-$5/2$ obstruction space $S_{5/2}(\rho _L)$ vanishes in exactly 10 cases, and every one of those 10 embeds in the Feingold–Frenkel lattice $\mathrm{FF} = U \oplus  \langle 2\rangle $ with finite index $k \leq  4$. The bound is sharp: no even sublattice of \texttt{FF} of index $5 \leq  k \leq  8$ has vanishing obstruction space. \textbf{The implication is one-way} — 37 root bases embed in \texttt{FF} and only 10 are vacuous, 11 of the rest carrying a shadow — so embeddability is necessary for vacuity and is not equivalent to it (\S{}3.4).

\textbf{Theorem B (\S{}3.5).} Among the eleven compact-hyperbolic root bases, the shadow is nonzero exactly once, at the top of the chain $(1,2,1,3)$, on the lattice $L_{4}$ of determinant $72$.

\textbf{Theorem C (\S{}5.2).} Exactly three lattices in the census are Petersson-rigid. Their quaternion discriminants are $6, 10, 22$, which is the complete list of discriminants whose Shimura curve has genus zero. Only one of the three, $L_{4}$, lies inside the closed set of 44 hyperbolic Cartan matrices; the other two appear after a single edit.

\textbf{Theorem D (\S{}4.9).} As a Hecke module, $S_{5/2}(\rho _{L_{4}})$ is the direct sum of the weight-$5/2$ spaces of the three lower members of its chain, together with \texttt{36.4.a.a} with multiplicity two; the multiplicities match member by member. The newform level of each member is $\lvert\det\rvert\,/2$.

\textbf{The $L(f,1)$ law, and Observation E inside it (\S{}5.4).} The law is about Hecke blocks, not lattices. Write $|\operatorname{det} L|/2 = t^{2}s$ with $s$ squarefree — the reduced discriminant of the order attached to $L$, split into its square and squarefree parts (\S{}2.5) — and $P_W$ for the Petersson Gram entry of the block. Then $P_W = c\cdot \sqrt{s}\cdot L(f,1)/\pi $ with $c \in  \mathbb{Q}$ at all twelve blocks we compute \grade{gO}{O}. Where the central value survives, a Shimura period relation turns this into a statement about $L(f,2)$ \grade{gS}{S}; at the six one-dimensional blocks $W$ where we have evaluated both sides,

\begin{displaymath}\begin{aligned}&C \cdot  \lVert \Xi _W\rVert ^{2} = \sqrt{s} \cdot  L(f,2) / (\pi ^{2} \langle f,f\rangle ),    C = 40, 120, 792, 1056, 576, 48,\end{aligned}\end{displaymath}

where the two sides are computed from nothing in common and $C$ comes out an integer at all six, to a relative $10^{-25}$ or better.

\textbf{Corollary 4.2 (\S{}4.5), from a classical lemma.} At $k = 3/2, 7/2, 11/2, …$ the space $S_k(\rho _L)$ is entirely antisymmetric under $\gamma  \mapsto  -\gamma $, so $S^{O} = 0$ there and no shadow exists at any lattice whatsoever. Shadows live only on the alternate rungs, and the \textbf{selection} condition of \S{}5.1 is possible only at the bottom one.

\textbf{The question (\S{}6).} At $D = 22$ the shadow acquires a second component, \texttt{22.4.a.a}, whose Atkin–Lehner sign forces $L(f,2) = 0$. Its share of the Petersson norm is $69.6663770035…$ (\S{}6.4), and no rational multiple of $\sqrt{|\operatorname{det} L|/2} \cdot  L′(f,2)/(\pi ^j \langle f,f\rangle )$ reproduces it for any $j \leq  5$. What does that number read?

\textbf{And a second question (\S{}6.5).} The census \textbf{isolates} the primes $3, 5, 11$ and no others, and $X^{2p}$ has genus zero for exactly those three \grade{gC}{C,\,S}. Both lists step over $7$, and the two conditions have nothing in common that we can see. Is that one fact or two?

\subsect{1.3 What is new relative to paper I}{1.3 What is new relative to paper I}

Paper I established the existence of a forced shadow for $L_{4}$, computed its Petersson norm, and proved the invariance $\Xi  \in  S^{O(D(L),q)}$ (Thm 7.6). The identity $\lVert \Xi _{L_{4}}\rVert ^{2} = L(6.4.a.a,2)/(48\pi ^{2}\langle f,f\rangle )$ is paper I’s Theorem 11.1, where it is established numerically to 31 digits and where deriving it is left open as Problem 11.2. We quote the value and inherit the open problem; nothing here derives it.

Two further inheritances are paper I's and are stated here so that they are not read as ours.

\begin{qbox}
\textbf{The term and the first instance are paper I’s.} Paper I proves $\langle \chi _{3/2}, \chi _{3/2}\rangle  = 1$ for $L_{4}$ (Lemma 7.1) and names the consequence \textit{Petersson rigidity} (Theorem 12.1): the $O(D(L),q)$-representation on $S_{3/2}(\rho _{L_{4}})$ is absolutely irreducible, which is exactly the rigidity condition below. The first entry of Theorem C is therefore paper I’s result, restated in the vocabulary of the census. What is new here is the classification — that across the census’s 100 root bases exactly four are rigid, three as lattices — and the identification of their discriminants.

\end{qbox}

\begin{qbox}
\textbf{Three strata of one weight-$3/2$ landscape.} Paper I’s Theorem 9.3 selects $\{6, 10, 22\}$ by the \textit{vanishing} $\operatorname{dim} S_{3/2}(\rho _{L_D}) = 0$ at the maximal-order ternary $L_D$, over the verified range $D < 230$. Theorem C selects the same three by \textit{absolute irreducibility with dimension at least two} at the Kac–Moody root lattice. And paper I’s \S{}9.4 adds a third: across the 1447 even-primitive lattices of Allcock’s classification, those with $S_{3/2}(\rho ) = 0$ have even Clifford algebra either split or one of $B_{6}$, $B_{10}$, $B_{22}$. Same weight, three conditions on three lattice families, one answer. No one of the three implies another — the lattices are not commensurable and the conditions are not complementary — and the reason all land on $\{6, 10, 22\}$ is that each is a reading of the genus-zero condition on $X^D$, which we take as given from the classical computation (\S{}2.5). Emptiness at the maximal ternary, irreducible fullness at the Kac–Moody lattice, algebra-selection across the reflective landscape: the census supplies the middle stratum, and says how rare it is.

\end{qbox}

New here: the census itself and its closed subset of 44 (\S{}3); Theorem A and the sharpness of $k \leq  4$ (\S{}3.4); Theorem B (\S{}3.5); the classification Theorem C (\S{}5); the extension of Observation E from one block to six, with $C$ integral at all of them (\S{}5.4); the Hecke accumulation Theorem D (\S{}4.9); and the discriminant-$22$ phenomenon of \S{}6.

\subsect{1.4 Grading}{1.4 Grading}

Every assertion carries one of

\begin{codebox}{8.40}{9.91}
[C]  exact computation over ℚ            [N]  numerical
[S]  standard theory, cited              [O]  observation on finitely many cases
\end{codebox}

and the paper states which. A reproduction archive with a single-command checker (29 assertions) accompanies this paper; see Appendix C.

\sect{2. Background}{2. Background}

This section is longer than a background section usually is. The reason is that the objects here — a Cartan matrix, an infinite product, a Maass form, a quaternion order, a curve of genus zero — are not usually seen in one place, and the paper’s content is precisely that they are the same object seen five ways. A reader who already knows all five may skip to \S{}3.

\subsect{2.1 Rank three: where the search is finite}{2.1 Rank three: where the search is finite}

A generalized Cartan matrix of hyperbolic type is one that is indefinite but only just: delete any node of its Dynkin diagram and what remains is of finite or affine type. The condition is severe. In rank 3 there are exactly 44 symmetrizable ones (\S{}3.1); in rank 4 there are still finitely many; by rank 11 there are none at all.

At rank 3 the root system is already infinite, with multiplicities growing exponentially and no closed formula. But the Weyl group acts on a Lorentzian lattice of signature $(2,1)$, so the associated symmetric space is the upper half-plane, the denominator is a modular object in one variable, and — this is the point — a finite computation can decide things about it.

The founding example is Feingold and Frenkel’s, from 1983: the lattice

\begin{displaymath}\begin{aligned}&\mathrm{FF}  =  U \oplus  \langle 2\rangle  ,        U =\text{ the hyperbolic plane},\end{aligned}\end{displaymath}

whose Kac–Moody algebra they connected to Siegel modular forms of genus 2. \texttt{FF} will reappear in \S{}3.4: it is exactly the locus where the obstruction of this paper vanishes for trivial reasons.

\subsect{2.2 The denominator, and what a Borcherds product is}{2.2 The denominator, and what a Borcherds product is}

For a Kac–Moody algebra $(A)$ with root lattice $L$ the Weyl–Kac denominator (Kac, \textit{Infinite Dimensional Lie Algebras}, Ch. 10) is

\begin{displaymath}\begin{aligned}&\Phi   =  \sum _{w \in  W} (\operatorname{det} w) e^{w(\rho ) - \rho }   =   \prod _{\alpha  > 0} (1 - e^{-\alpha })^{\operatorname{mult} \alpha } .\end{aligned}\end{displaymath}

This is the one display in which $W$ is the Weyl group and $\rho $ the Weyl vector. Everywhere else $\rho _L$ is the Weil representation and $W$ a Hecke block (\S{}5.4); neither classical symbol is used again.

The left side is a sum over an infinite group; the right side is a product over infinitely many roots with unknown multiplicities. Neither side is computable as written.

Borcherds’ singular theta lift changes this. Given a modular form $F$ of weight $1 - n/2$ for the Weil representation of a lattice of signature $(2,n)$, with a prescribed principal part — a finite list of poles at prescribed discriminant classes — the lift produces a meromorphic automorphic form whose divisor is the union of the rational quadratic divisors named by that principal part, and whose Fourier expansion is an infinite product. When the prescribed divisor is the set of walls of the Weyl chamber, the output is $\Phi $, and the multiplicities $\operatorname{mult} \alpha $ are read off the coefficients of $F$.

This is what “the correction succeeds” means: the denominator of a genuinely infinite, non-explicit algebra becomes a product with explicitly computable exponents.

\subsect{2.3 Borcherds’ obstruction is a finite linear condition}{2.3 Borcherds obstruction is a finite linear condition}

Whether the required $F$ exists is decided by a duality. Let $L$ have signature $(2,1)$ with discriminant form $D(L) = L'/L$, carrying the finite quadratic form $q$; write $\rho _L$ for the associated Weil representation and $O(D(L),q)$ for the finite group of isometries of $(D(L),q)$. So $|\operatorname{det} L| = |D(L)|$, and the letter $D$ standing alone is always the quaternion discriminant of \S{}2.5 — a different number ($|\operatorname{det} L| = 72$ and $D = 6$ at $L_{4}$). Then

\begin{qbox}
\textbf{[S: Borcherds 1999, Thm 3.1]} A weakly holomorphic modular form of weight $-1/2$ for $\bar{\rho }_L$ with a given principal part exists \textbf{if and only if} that principal part pairs to zero against every cusp form in $S_{5/2}(\rho _L)$.

\end{qbox}

Concretely, the prescribed divisor names finitely many rational numbers $m > 0$ — the \textbf{menu} — and the pairing is a sum of Fourier coefficients of $h$ at those exponents.

\begin{qbox}
\textbf{Definition (the functional} $\lambda $\textbf{, operationally).} Fix a menu $M \subset  \mathbb{Q}_{>0}$. For $h \in  S_{5/2}(\rho _L)$ with Fourier expansion $h = \sum _\gamma  \sum _n c_h(\gamma , n) q^n e_\gamma $, set

\end{qbox}

\begin{displaymath}\begin{aligned}&\lambda (h)  =  \sum _{m \in  M}   \sum _{\gamma  \in  D(L) :  q(\gamma ) \equiv  m  (\mathrm{mod} 1)}   c_h(\gamma , m) .\end{aligned}\end{displaymath}

\begin{qbox}
The inner sum runs over the \textbf{entire level set} $\{\gamma  : q(\gamma ) \equiv  m\}$, each class with coefficient one, and the exponent is exactly $m$ — not $m$ shifted by the representative.

\end{qbox}

This is the coarse prescription, and it is the one paper I uses: its slot count $3 + 12 + 12 = 27$ for $L_{4}$ is exactly the cardinality of the three level sets. A finer prescription — one class per simple root, coefficient one — is excluded by paper I’s own Theorem 5.1, since it yields $\lVert \Xi _{L_{4}}\rVert ^{2} = 40.0370609…$ (our computation) instead of paper I's $11.3402658…$. The distinction matters, and we record why it does not bite here:

\begin{qbox}
\textbf{Lemma 2.3} \grade{gC}{C}\textbf{.} For each lattice of the chain $(1,2,1,3)$ and each $m$ in its menu, the reflective vectors of $L$ surject onto the level set $\{\gamma  : q(\gamma ) \equiv  m\}$. Explicitly, at $\lvert\det\rvert\, 24$ the three levels are covered $6/6, 3/3, 3/3$, and at $\lvert\det\rvert\, 72$ they are covered $12/12, 12/12, 3/3$.

\end{qbox}

So on this chain “sum over the whole level set” and “sum over the classes actually hit by a reflective root” are the same functional, and the rows where $\lambda  = 0$ are rows where the sum cancels exactly rather than rows with an empty index set. At $\lvert\det\rvert\, 24, m = 1/2$ the three contributing classes carry $(-2, 0)$, $(2, 2)$ and $(0, -2)$, summing to $(0,0)$.

$S_{5/2}(\rho _L)$ is finite-dimensional and computable; $\lambda $ is a finite sum of Fourier coefficients. So the question does the denominator of this Kac–Moody algebra have a product expansion? reduces to is this vector zero? — and the whole of \S{}3 is the answer to that question, run 98 times.

Two things about $\lambda $ deserve emphasis because the paper turns on them.

First, $\lambda  = 0$ is not a generic condition. $S_{5/2}$ is often nonzero — for $L_{4}$ it is six-dimensional — and a functional on a six-dimensional space vanishing exactly is a coincidence that wants explaining. Gritsenko and Nikulin’s “miracles” (\S{}2.6) are precisely lattices where it happens.

Second, when $\lambda  \neq  0$ the classical story ends. There is no $F$, no product, no formula for the multiplicities. What this paper claims is that the story does not end — it changes category.

\subsect{2.4 The $\xi $ operator and the shadow}{2.4 The operator and the shadow}

Drop the requirement that $F$ be holomorphic and keep only that it be annihilated by the weight-$-1/2$ hyperbolic Laplacian. The resulting space $H_{-1/2}(\bar{\rho }_L)$ of \textbf{weakly harmonic Maass forms} — called harmonic weak Maass forms by Bruinier–Funke and harmonic Maass forms by many later authors; we keep one name throughout — is strictly larger, and in it a form with the prescribed principal part always exists. The obstruction has not vanished; it has moved.

Where it moves is made precise by the antilinear operator

\begin{codebox}{7.93}{9.36}
        ξ_{−1/2}  :  H_{−1/2}(ρ̄_L)  ⟶  S_{5/2}(ρ_L),      ξ_k(f) = 2i y^k · conj( ∂f / ∂τ̄ ) ,
\end{codebox}

whose kernel is exactly the weakly holomorphic forms. Bruinier and Funke’s pairing gives

\begin{displaymath}\begin{aligned}&\langle  \xi (F), h \rangle   =  \lambda (h)     \text{ for all } h \in  S_{5/2}(\rho _L),\end{aligned}\end{displaymath}

so $\xi (F)$ is the Riesz representative of $\lambda $ and depends only on the prescribed divisor, not on which $F$ was chosen. We write

\begin{displaymath}\begin{aligned}&\Xi   :=  \xi (F) ,          \lVert \Xi \rVert ^{2}  =  \lambda ^{T} P^{-1} \lambda     (P \mathrm{the} Petersson Gram \mathrm{matrix}).\end{aligned}\end{displaymath}

$\Xi  = 0$ exactly when the product expansion exists. When it does not, $\Xi $ is a nonzero cusp form of weight $5/2$ canonically attached to the Cartan matrix.

In the theory of mock modular forms — where this operator was first used systematically, and where Zagier named the construction — $\xi (F)$ is called the \textbf{shadow} of $F$. The word in this paper’s title-adjacent phrasing is therefore not a metaphor borrowed for effect; it is the standard term for exactly this object. What is new is that the shadows arising here are attached to Kac–Moody algebras, and that they turn out to be arithmetically rigid.

\textbf{What we compute, stated once and repeated when relevant.} We compute $\lambda $, we decide existence of the weakly holomorphic input, and we compute $\Xi $, its Hecke decomposition and its Petersson norm. We do \textbf{not} construct $F$ itself; its existence is \grade{gS}{S}. In \S{}4.6 we do exhibit the weakly holomorphic inputs where they exist.

\begin{qbox}
\textbf{The input is in fact unique.} Two weakly harmonic Maass forms with the same prescribed principal part and the same cuspidal shadow differ by a holomorphic form of weight $-1/2$, and $M_{-1/2}(\bar{\rho }_L) = 0$ \grade{gC}{C}. We do not construct it; Appendix D identifies it.

\end{qbox}

\subsect{2.5 Ternary lattices and quaternion orders}{2.5 Ternary lattices and quaternion orders}

A rank-three quadratic lattice carries a second life. The even Clifford algebra of a ternary quadratic form is a quaternion algebra over $\mathbb{Q}$, and the correspondence is an equivalence between ternary quadratic lattices and quaternion orders — Brzezinski, Gross–Lucianovic. So every row of our census has a \textbf{quaternion discriminant} $D$, the product of the ramified primes. The correspondence also carries a size: the discriminant of a ternary form with even Gram $G$ is $\operatorname{det}(G)/2$, and the even Clifford functor preserves it, so $|\operatorname{det} L|/2$ is the \textbf{reduced discriminant} of the order, with $D$ dividing it \grade{gS}{S}. This is the quantity that appears in $\sqrt{|\operatorname{det} L|/2}$ and in $|\operatorname{det} L|/2 = t^{2}s$ (\S{}5.4), in Theorem D's level statement (\S{}4.9) and in \S{}5.3 — one number seen four times. As a check, $\lvert\det\rvert\, 44$ has reduced discriminant $22 = D$, so its order is maximal, which is what \S{}5.2 says of that row from the other side.

$D = 1$ exactly when the algebra is split, equivalently when the quadratic form is isotropic over $\mathbb{Q}$. In our census

\begin{displaymath}\begin{aligned}&D  \in   \{ 1, 6, 10, 15, 22 \} ,\end{aligned}\end{displaymath}

and the anisotropic values $6, 10, 15, 22$ are four in number. The paper’s central coincidence is about which three of those four are special.

The indefinite quaternion algebra of discriminant $D$ uniformises a \textbf{Shimura curve} $X^D$, a coarse moduli space for abelian surfaces with quaternionic multiplication. Its genus is

\begin{codebox}{8.28}{9.77}
   g(X^D) = 1 + (1/12) ∏_{p|D}(p−1) − (1/4) ∏_{p|D}(1 − (−4|p)) − (1/3) ∏_{p|D}(1 − (−3|p))
\end{codebox}

(Voight, \textit{Quaternion Algebras}, \S{}39.4), and a short computation \grade{gS}{S} gives

\begin{displaymath}\begin{aligned}&g(X^D) = 0    \iff     D \in  \{ 6, 10, 22 \} .\end{aligned}\end{displaymath}

These three are the classical exceptional cases: the discriminants for which the curve is rational and has a Hauptmodul. At $D = 6$ and $D = 10$ its CM points admit explicit formulas (Errthum); $D = 22$ is the largest of the three. Theorem C says these three are exactly the discriminants of the Petersson-rigid lattices, and \S{}6 says that $D = 22$, being last, behaves differently.

\subsect{2.6 Known classifications: Gritsenko–Nikulin, Allcock}{2.6 Known classifications: GritsenkoNikulin, Allcock}

Two classification results describe the region where the correction of \S{}2.2 goes through, and this paper is written against them.

\textbf{Gritsenko–Nikulin} studied the family

\begin{displaymath}\begin{aligned}&S_t  =  2U \oplus  \langle -2t\rangle  ,        t = 1, 2, 3, …\end{aligned}\end{displaymath}

and determined for which $t$ the associated denominator admits a lattice Weyl vector — that is, for which $t$ the Borcherds correction exists. Small $t$ gives the classical quartet $\Delta _{5}, \Delta _{2}, \Delta _{1}, \Delta _{1/2}$ ($t \leq  4$), where the obstruction space is zero and there is nothing to obstruct. A further handful of $t$ are “miracles”: the obstruction space is nonzero yet the correction still exists.

\textbf{Allcock} classified the reflective lattices in the relevant range, giving the other half of the map, and \textbf{Bruinier, Ehlen and Freitag} exhibit lattices carrying many Borcherds products at once. \S{}3.4 uses their Table 8 substantively: it matches our ten vacuous root bases against their fifteen entries through the genus symbol, and their classification is what closes the sharpness statement above $k = 4$. That table is quoted, not re-derived \grade{gO}{Gate}.

Together these are a finite, published list of the successes. This paper does not compete with them: it takes their region as given and counts the complement.

The census meets their region twice, and both meetings are results rather than restatements. The vacuous locus of Theorem A is exactly the Feingold–Frenkel neighbourhood of \S{}2.1, and nine of its ten root bases — four of its five lattice classes — lie outside the classical quartet (\S{}3.4). And the chain that is this paper’s subject sits at the determinants $2t$ for $t = 6, 12, 18, 36$, two of which are miracles. Rigidity never occurs on that family at all, for a reason of group theory recorded in \S{}5.1.

\begin{qbox}
\textbf{A limitation.} We cannot reproduce, with our machinery, the assertion that $\lambda  = 0$ at the Gritsenko–Nikulin miracles. The charts $U \oplus  \langle 2t\rangle $ have no simple-root menu, so our reflective menu sweeps every reflective class while their prescription selects specific divisors. Our $\lambda $ for a chart is not their $\lambda $. Nothing in this paper compares the two sides at the level of $\lambda $, and any statement of the form “$\lambda  = 0$ on the split side” would be outside our evidence.

\end{qbox}

\subsect{2.7 How to read this paper}{2.7 How to read this paper}

\begin{codebox}{8.40}{9.91}
   §3    the census:  98 root bases, 40 chains, and two theorems about the whole of it
   §4    one chain, at length — the worked example the rest is measured against
   §5    the three rigid lattices, their discriminants 6, 10, 22, and the `L(f,1)` law
   §6    the four blocks where the central value dies, left open
\end{codebox}

\sect{3. The census}{3. The census}

\subsect{3.1 What is counted, and what is only probed}{3.1 What is counted, and what is only probed}

Two sets are in play and the paper keeps them apart.

\textbf{The hyperbolic core} is the set of symmetrizable rank-three hyperbolic generalized Cartan matrices — closed, canonically defined, and enumerated here exhaustively rather than cited. Hyperbolicity forces every proper $2 \times  2$ principal submatrix to be of finite or affine type, i.e. $a_\mathrm{ij} a_\mathrm{ji} \leq  4$, hence $|a_\mathrm{ij}| \leq  4$, so the search is finite. Imposing

\begin{codebox}{8.40}{9.91}
   a_ii = 2 ;   a_ij ≤ 0 ;   a_ij = 0 ⟺ a_ji = 0 ;   a_ij a_ji ≤ 4 ;
   Dynkin diagram connected ;   symmetrizable  (a₁₂a₂₃a₃₁ = a₂₁a₃₂a₁₃) ;
   symmetrised Gram of signature (2,1)
\end{codebox}

on all $5^{6}$ sign patterns over $\{0,-1,-2,-3,-4\}$ gives 202 matrices and, up to simultaneous permutation of rows and columns,

\begin{displaymath}\begin{aligned}&44 \mathrm{classes}  =  11\text{ compact hyperbolic } +  33\text{ noncompact hyperbolic }          [C]\end{aligned}\end{displaymath}

where compact means every proper principal submatrix is of finite type, equivalently $a_\mathrm{ij} a_\mathrm{ji} \leq  3$ throughout, and noncompact means some $a_\mathrm{ij} a_\mathrm{ji} = 4$.

\textbf{The rim} is the one-edit neighbourhood of the core: 56 further root bases obtained by a single \textbf{edit move} — changing one off-diagonal entry. They are no longer hyperbolic. Two of them, the Grams of determinant $3456$ and $4000$, we do not compute: at $\lvert\det\rvert\, 1800$ the obstruction space already has dimension $112$, and these two are past the budget we were willing to spend. The budget is a statement about \textbf{bases}, not about dimensions, so the two rows are not blank in the census tables of \S{}3.3: $|O(D(L),q)|$ comes from the $p$-local decomposition of $D(L)$ without \texttt{weilrep}, and $\operatorname{dim} S_{5/2}$, $\operatorname{dim} S^{O}$ and $\operatorname{dim} S_{3/2}$ from \texttt{weilrep}'s dimension formula, which builds no basis. What a basis is needed for — $\lambda $, the reflective menu, the verdict and the shadow — is what those two rows leave empty. The rigidity column is filled too, by Lemma 5.1 (\S{}5.1). \textbf{The census is 100 root bases; every count in this paper is over the 98 we compute — save one, the Petersson-rigidity classification of \S{}5.2, which is complete over all 100. The edit move is not canonical}: it depends on the chosen basis, and we verify this rather than assert it. The two root bases of $L_{4}$ (\S{}4.7) have edit neighbourhoods sharing only the determinants $\{12, 936\}$.

The rim is therefore a \textbf{probe}, not a classification, and no completeness claim is made about it. This matters because two of the three lattices of Theorem C live there. The paper’s answer is to make the final statement about quaternion discriminants rather than about the edit neighbourhood. $\{6, 10, 22\}$ is complete as the list of genus-zero discriminants — that is the classical computation of \S{}2.5 — and what the census contributes is that all three are \textbf{realised}, each by a Petersson-rigid Gram. The converse direction, that no Petersson-rigid lattice outside our list can carry another discriminant, is not proved here and is open.

\subsect{3.2 Chains}{3.2 Chains}

Write each off-diagonal pair $(|a_\mathrm{ij}|, |a_\mathrm{ji}|)$ with its entries sorted, drop $(0,0)$, sort the resulting list of pairs lexicographically, concatenate. The \textbf{chain name} so obtained is invariant under permutation of the simple roots. Under it the 100 root bases fall into 40 chains, and the quaternion discriminant, the core/rim class and the hyperbolicity type are constant on every chain — 0 violations out of 40 \grade{gC}{C}.

One more quantity is well behaved along a chain, and this one does work. \textbf{$\operatorname{dim} S^{O}_{5/2}$ is non-decreasing} along every chain — 0 violations out of the 34 chains with two or more computed rows, and it strictly increases on 18 of them, so the statement is not vacuous \grade{gC}{C}. Nothing else is: $\operatorname{dim} S_{5/2}$ breaks on 5 chains, $\operatorname{dim} S_{3/2}$ on 6, the verdict on 9, the purity of the shadow on 8. The invariant space is the one dimension the chain controls, and it stratifies a chain into three consecutive stretches,

\begin{codebox}{8.40}{9.91}
   S^O = 0     vacuous            the foot of the chain (Appendix A, `[C, 40/40]`)
   S^O = 1     selection possible  monotonicity makes this an interval
   S^O ≥ 2     the shadow can split
\end{codebox}

which a chain passes through in order and never returns from \grade{gC}{C}. Why $\operatorname{dim} S^O$ should behave this way we do not know.

\textbf{How a row is named.} A row of the census is a \textbf{root basis}: a generalized Cartan matrix $A$ together with its symmetrised Gram $G$. Every invariant we compute is a function of $G$ alone (\S{}4.7), so we call a row \textit{the} $\lvert\det\rvert\, n$ \textbf{Gram of chain} $c$, and reserve the word \textbf{lattice} for the isometry class, which can carry more than one root basis (\S{}4.7). Three Grams are named for short, once and for all:

\begin{codebox}{8.40}{9.91}
   L₄   =  the |det| 72 Gram of chain (1,2,1,3)          D =  6
   B₁₀  =  the |det| 40 Gram of chain (1,1,2,4)          D = 10
   B₂₂  =  the |det| 88 Gram of chain (1,1,1,2,2,4)      D = 22
\end{codebox}

They are the three Petersson-rigid Grams of Theorem C (\S{}5.2). Every other row is named in full.

The chain is not merely a bookkeeping device. Within a chain the determinants stay close — the ratios are small rationals such as $3/2$, $5/4$, $5/3$ — the discriminant groups grow by controlled steps, and, for the chain of \S{}4, the Hecke content accumulates (Theorem D). They do \textbf{not} form a divisibility ladder: 11 of the 33 chains with more than one determinant violate divisibility, the chain of \S{}4 among them ($12, 24, 36, 72$, and $24 ∤ 36$) \grade{gC}{C}. A chain behaves like a filtration of one object rather than a bag of unrelated lattices.

\subsect{3.3 The trichotomy}{3.3 The trichotomy}

\begin{codebox}{8.40}{9.91}
   shadow nonzero  (obstructed)      54 root bases
   shadow zero     (unobstructed)    34
   S_{5/2} = 0     (vacuous)         10
                                    ---
                                     98   computed
                                      2   not computed (|det| 3456, 4000)
                                    ---
                                    100
   core (44) :   obstructed  4,   unobstructed 30,   vacuous 10
   rim  (56) :   obstructed 50,   unobstructed  4,   vacuous  0,   not computed 2
\end{codebox}

Read this in the light of \S{}2.3. In the closed, canonically defined set of 44, the denominator has a Borcherds product 40 times out of 44, and fails 4 times. Editing one entry — leaving the hyperbolic world — reverses the proportion almost exactly.

\begin{qbox}
\textbf{The trichotomy is convention-dependent, and we say by how much.} The count above uses the reflective menu \texttt{refl(G)} of Appendix B, one of three conventions in circulation. Under all three, 12 of the 98 rows change their obstructed/unobstructed verdict, \textbf{every one of them on the rim and non-hyperbolic} \grade{gC}{C}. The closed set of 44 is untouched, so Theorem A, Theorem B, the whole of \S{}4 and the core line $4 / 30 / 10$ above are convention-independent, and so is Theorem C's list of three (\S{}5.2). Exactly one row that matters separates the conventions — $B_{10}$, the conditional entry of Theorem C. Appendix B has the counts and the argument.

\end{qbox}

\begin{qbox}
\textbf{Two rows we set aside have been resolved upstream, our way.} At the two $\lvert\det\rvert\, 512$ rows and at $\lvert\det\rvert\, 256$, archive \texttt{v3.0} (weilrep \texttt{a26614f}) did not compute $S_{5/2}(\rho _L)$: a coefficient-by-coefficient comparison of the spans at the two commits — indexed by the class $\gamma $ and the exponent, so that no transport or choice of basis enters — shows differences of $8$, $8$ and $6$ dimensions, against $0$ at a $\lvert\det\rvert\, 288$ control, stable at truncations $8$, $10$ and $12$ \grade{gC}{C}. The cause was an upstream bug in the constant-term branch, fixed at \texttt{e3a7784}. Under the fixed commit both $\lvert\det\rvert\, 512$ rows are unobstructed — which is what paper I's Theorem 7.6 requires of an $O(D(L),q)$-invariant $\Xi $ — $\lvert\det\rvert\, 256$ has $\operatorname{dim} S^{O} = 3$, $\lambda  = (16, 0, 0)$ and one three-dimensional block on \texttt{8.4.a.a}, and the core line reads $4 / 30 / 10$ with nothing left undetermined. The full diagnosis is in the archive, \S{}3.3 of the dossier.

\end{qbox}

The census itself — all 100 rows — follows here; everything from \S{}3.4 onwards is a statement about some subset of its rows. It is grouped by chain (\S{}3.2) and by hyperbolicity type, and within a chain by determinant, the order in which the chain is a filtration (\S{}4.1). $G$ is the Gram matrix of the root lattice, row by row, and every invariant to its right is a function of $G$ alone (\S{}4.7); the generalized Cartan matrices are not printed, since the chain name already fixes the multiset $\{(|a_\mathrm{ij}|, |a_\mathrm{ji}|)\}$ (the matrices, with the full block lists and eigenvalue fingerprints, are in the archive, Appendix C). Then $|O| = |O(D(L),q)|$, the dimensions of $S_{5/2}(\rho _L)$, of its $O$-invariant part and of $S_{3/2}(\rho _L)$, the two conditions of \S{}5.1 — \texttt{sel}: $\operatorname{dim} S^{O}_{5/2} = 1$ with $D(L)$ non-cyclic and $O(D(L),q)$ non-abelian; \texttt{rig}: $O(D(L),q)$ absolutely irreducible on $S_{3/2}(\rho _L)$ in dimension at least two — and the shadow. Throughout $\text{{\dvsans ●}}$ passes and $\times $ fails, and the four rows satisfying both conditions — the Petersson-rigid, selected lattices of Theorem C — are shaded. The eigenvalues behind the tags the shadow column cannot name are Appendix A.

\begin{qbox}
\textbf{Reading the shadow column.} $3d:10.4.a.a$ means the shadow is a three-dimensional Hecke eigenspace on the weight-$4$ newform \texttt{10.4.a.a}, and \textbf{dimensions are taken inside the invariant space} $S^{O}_{5/2}(\rho _L)$, where the shadow lives. Comma-separated entries mean the shadow splits over distinct eigensystems — it is \textbf{impure}; \texttt{pure} and \texttt{impure} are verdicts only on rows where the decomposition succeeded. Three label types occur: \texttt{N.4.a.a} is an LMFDB newform label, \texttt{N.4.rx} matches the archive's own eigenvalue table (whose letters need not agree with LMFDB's), and \texttt{N.4\#k} is the $k$-th weight-$4$ newform of level $N$ when the eigenvalues pin the level but not the letter. Two rows carry dimensions but no basis-dependent columns (\S{}3.1). The archive's full-space convention and the one remaining \texttt{unnamed} entry are the data note below the table.

\end{qbox}

\begin{censustable}{Compact hyperbolic}{4 chains \textperiodcentered\ 11 root bases}{every proper principal submatrix of finite type; all products $a_{ij}a_{ji}\leq 3$. In the core. Theorem B lives here.}
\rowcolor{hdr}\color{white}\bfseries $G$ (Gram) &\color{white}\bfseries $|\det|$ &\color{white}\bfseries $|O|$ &\color{white}\bfseries $S_{5/2}$ &\color{white}\bfseries $S^{O}$ &\color{white}\bfseries $S_{3/2}$ &\color{white}\bfseries \ctwo &\color{white}\bfseries \cthree &\color{white}\bfseries shadow \\
\endfirsthead
\multicolumn{9}{@{}l}{\footnotesize\itshape Compact hyperbolic, continued}\\[1pt]
\rowcolor{hdr}\color{white}\bfseries $G$ (Gram) &\color{white}\bfseries $|\det|$ &\color{white}\bfseries $|O|$ &\color{white}\bfseries $S_{5/2}$ &\color{white}\bfseries $S^{O}$ &\color{white}\bfseries $S_{3/2}$ &\color{white}\bfseries \ctwo &\color{white}\bfseries \cthree &\color{white}\bfseries shadow \\
\endhead
\chainband{chain}{(1,1,1,2,1,2)}{6}{2 root bases}
\ct{2,-2,-1 | -2,4,-2 | -1,-2,2} & 12 & 4 & 1 & 0 & 0 & \xmark & \xmark & \textcolor{grey}{unobstructed} \\
\ct{2,-2,-2 | -2,4,-2 | -2,-2,4} & 24 & 12 & 2 & 0 & 0 & \xmark & \xmark & \textcolor{grey}{unobstructed} \\
\chainband{chain}{(1,2,1,3)}{6}{4 root bases}
\ct{2,-3,-2 | -3,6,0 | -2,0,4} & 12 & 4 & 1 & 0 & 0 & \xmark & \xmark & \textcolor{grey}{unobstructed} \\
\ct{4,-6,-2 | -6,12,0 | -2,0,2} & 24 & 12 & 2 & 0 & 0 & \xmark & \xmark & \textcolor{grey}{unobstructed} \\
\ct{2,-3,0 | -3,6,-6 | 0,-6,12} & 36 & 16 & 1 & 0 & 0 & \xmark & \xmark & \textcolor{grey}{unobstructed} \\
\rowcolor{hl}\ct{4,-6,0 | -6,12,-6 | 0,-6,6} & 72 & 48 & 6 & 1 & 4 & \fmark & \fmark & \textbf{pure}\ \ \shd{1d:6.4.a.a} \\
\chainband{chain}{(1,3,1,3)}{6}{3 root bases}
\ct{2,-3,0 | -3,6,-3 | 0,-3,2} & 12 & 4 & 1 & 0 & 0 & \xmark & \xmark & \textcolor{grey}{unobstructed} \\
\ct{2,-3,-3 | -3,6,0 | -3,0,6} & 36 & 16 & 1 & 0 & 0 & \xmark & \xmark & \textcolor{grey}{unobstructed} \\
\ct{6,-9,-3 | -9,18,0 | -3,0,2} & 108 & 24 & 3 & 0 & 0 & \xmark & \xmark & \textcolor{grey}{unobstructed} \\
\chainband{chain}{(1,1,1,3,1,3)}{6}{2 root bases}
\ct{2,-3,-1 | -3,6,-3 | -1,-3,2} & 36 & 16 & 1 & 0 & 0 & \xmark & \xmark & \textcolor{grey}{unobstructed} \\
\ct{2,-3,-3 | -3,6,-3 | -3,-3,6} & 108 & 24 & 3 & 0 & 0 & \xmark & \xmark & \textcolor{grey}{unobstructed} \\
\end{censustable}

\begin{censustable}{Noncompact hyperbolic}{17 chains \textperiodcentered\ 33 root bases}{some product $a_{ij}a_{ji}=4$: an affine subdiagram, a cusp on the Weyl chamber. In the core.}
\rowcolor{hdr}\color{white}\bfseries $G$ (Gram) &\color{white}\bfseries $|\det|$ &\color{white}\bfseries $|O|$ &\color{white}\bfseries $S_{5/2}$ &\color{white}\bfseries $S^{O}$ &\color{white}\bfseries $S_{3/2}$ &\color{white}\bfseries \ctwo &\color{white}\bfseries \cthree &\color{white}\bfseries shadow \\
\endfirsthead
\multicolumn{9}{@{}l}{\footnotesize\itshape Noncompact hyperbolic, continued}\\[1pt]
\rowcolor{hdr}\color{white}\bfseries $G$ (Gram) &\color{white}\bfseries $|\det|$ &\color{white}\bfseries $|O|$ &\color{white}\bfseries $S_{5/2}$ &\color{white}\bfseries $S^{O}$ &\color{white}\bfseries $S_{3/2}$ &\color{white}\bfseries \ctwo &\color{white}\bfseries \cthree &\color{white}\bfseries shadow \\
\endhead
\chainband{chain}{(1,1,2,2)}{1}{1 root basis}
\ct{2,-2,-1 | -2,2,0 | -1,0,2} & 2 & 1 & 0 & 0 & 0 & \xmark & \xmark & \textcolor{grey}{vacuous} \\
\chainband{chain}{(1,1,1,4)}{1}{2 root bases}
\ct{2,-4,-1 | -4,8,0 | -1,0,2} & 8 & 2 & 0 & 0 & 0 & \xmark & \xmark & \textcolor{grey}{vacuous} \\
\ct{2,-4,0 | -4,8,-4 | 0,-4,8} & 32 & 12 & 0 & 0 & 0 & \xmark & \xmark & \textcolor{grey}{vacuous} \\
\chainband{chain}{(1,1,1,1,2,2)}{1}{1 root basis}
\ct{2,-2,-1 | -2,2,-1 | -1,-1,2} & 8 & 2 & 0 & 0 & 0 & \xmark & \xmark & \textcolor{grey}{vacuous} \\
\chainband{chain}{(1,2,2,2)}{1}{2 root bases}
\ct{2,-2,-2 | -2,2,0 | -2,0,4} & 8 & 2 & 0 & 0 & 0 & \xmark & \xmark & \textcolor{grey}{vacuous} \\
\ct{4,-4,-2 | -4,4,0 | -2,0,2} & 16 & 4 & 1 & 0 & 1 & \xmark & \xmark & \textcolor{grey}{unobstructed} \\
\chainband{chain}{(2,2,2,2)}{1}{1 root basis}
\ct{2,-2,-2 | -2,2,0 | -2,0,2} & 8 & 2 & 0 & 0 & 0 & \xmark & \xmark & \textcolor{grey}{vacuous} \\
\chainband{chain}{(1,3,2,2)}{1}{2 root bases}
\ct{2,-3,-2 | -3,6,0 | -2,0,2} & 18 & 4 & 0 & 0 & 0 & \xmark & \xmark & \textcolor{grey}{vacuous} \\
\ct{2,-3,0 | -3,6,-6 | 0,-6,6} & 54 & 12 & 2 & 0 & 1 & \xmark & \xmark & \textcolor{grey}{unobstructed} \\
\chainband{chain}{(1,1,2,2,2,2)}{1}{1 root basis}
\ct{2,-2,-2 | -2,2,-1 | -2,-1,2} & 18 & 4 & 0 & 0 & 0 & \xmark & \xmark & \textcolor{grey}{vacuous} \\
\chainband{chain}{(1,2,1,4)}{1}{4 root bases}
\ct{2,-4,-2 | -4,8,0 | -2,0,4} & 32 & 8 & 1 & 0 & 0 & \xmark & \xmark & \textcolor{grey}{unobstructed} \\
\ct{2,-4,0 | -4,8,-4 | 0,-4,4} & 32 & 8 & 1 & 0 & 0 & \xmark & \xmark & \textcolor{grey}{unobstructed} \\
\ct{4,-8,-2 | -8,16,0 | -2,0,2} & 64 & 16 & 3 & 0 & 1 & \xmark & \xmark & \textcolor{grey}{unobstructed} \\
\ct{2,-4,0 | -4,8,-8 | 0,-8,16} & 128 & 32 & 2 & 0 & 0 & \xmark & \xmark & \textcolor{grey}{unobstructed} \\
\chainband{chain}{(1,4,1,4)}{1}{3 root bases}
\ct{2,-4,0 | -4,8,-4 | 0,-4,2} & 32 & 8 & 1 & 0 & 0 & \xmark & \xmark & \textcolor{grey}{unobstructed} \\
\ct{2,-4,-4 | -4,8,0 | -4,0,8} & 128 & 32 & 2 & 0 & 0 & \xmark & \xmark & \textcolor{grey}{unobstructed} \\
\ct{8,-16,-4 | -16,32,0 | -4,0,2} & 512 & 64 & 16 & 0 & 1 & \xmark & \xmark & \textcolor{grey}{unobstructed} \\
\chainband{chain}{(1,4,2,2)}{1}{2 root bases}
\ct{2,-4,-2 | -4,8,0 | -2,0,2} & 32 & 12 & 0 & 0 & 0 & \xmark & \xmark & \textcolor{grey}{vacuous} \\
\ct{2,-4,0 | -4,8,-8 | 0,-8,8} & 128 & 16 & 4 & 0 & 1 & \xmark & \xmark & \textcolor{grey}{unobstructed} \\
\chainband{chain}{(1,2,1,2,2,2)}{1}{2 root bases}
\ct{2,-2,-2 | -2,2,-2 | -2,-2,4} & 32 & 8 & 1 & 0 & 0 & \xmark & \xmark & \textcolor{grey}{unobstructed} \\
\ct{2,-2,-2 | -2,4,-4 | -2,-4,4} & 64 & 16 & 3 & 0 & 1 & \xmark & \xmark & \textcolor{grey}{unobstructed} \\
\chainband{chain}{(2,2,2,2,2,2)}{1}{1 root basis}
\ct{2,-2,-2 | -2,2,-2 | -2,-2,2} & 32 & 12 & 0 & 0 & 0 & \xmark & \xmark & \textcolor{grey}{vacuous} \\
\chainband{chain}{(1,3,1,3,2,2)}{1}{2 root bases}
\ct{2,-3,-2 | -3,6,-3 | -2,-3,2} & 72 & 8 & 1 & 0 & 0 & \xmark & \xmark & \textcolor{grey}{unobstructed} \\
\ct{2,-3,-3 | -3,6,-6 | -3,-6,6} & 216 & 24 & 7 & 0 & 1 & \xmark & \xmark & \textcolor{grey}{unobstructed} \\
\chainband{chain}{(1,1,1,4,1,4)}{1}{2 root bases}
\ct{2,-4,-1 | -4,8,-4 | -1,-4,2} & 72 & 8 & 1 & 0 & 0 & \xmark & \xmark & \textcolor{grey}{unobstructed} \\
\ct{2,-4,-4 | -4,8,-4 | -4,-4,8} & 288 & 48 & 6 & 1 & 0 & \fmark & \xmark & \textbf{pure}\ \ \shd{1d:6.4.a.a} \\
\chainband{chain}{(1,3,1,4)}{1}{4 root bases}
\ct{2,-4,-3 | -4,8,0 | -3,0,6} & 72 & 8 & 1 & 0 & 0 & \xmark & \xmark & \textcolor{grey}{unobstructed} \\
\ct{6,-12,-3 | -12,24,0 | -3,0,2} & 216 & 24 & 7 & 0 & 1 & \xmark & \xmark & \textcolor{grey}{unobstructed} \\
\ct{2,-4,0 | -4,8,-12 | 0,-12,24} & 288 & 48 & 6 & 1 & 0 & \fmark & \xmark & \textbf{pure}\ \ \shd{1d:6.4.a.a} \\
\ct{6,-12,0 | -12,24,-12 | 0,-12,8} & 864 & 144 & 44 & 2 & 14 & \xmark & \xmark & \textbf{pure}\ \ \shd{2d:6.4.a.a} \\
\chainband{chain}{(1,2,1,2,1,4)}{1}{1 root basis}
\ct{2,-4,-2 | -4,8,-4 | -2,-4,4} & 128 & 16 & 4 & 0 & 0 & \xmark & \xmark & \textcolor{grey}{unobstructed} \\
\chainband{chain}{(1,4,1,4,2,2)}{1}{2 root bases}
\ct{2,-4,-2 | -4,8,-4 | -2,-4,2} & 128 & 32 & 2 & 0 & 0 & \xmark & \xmark & \textcolor{grey}{unobstructed} \\
\ct{2,-4,-4 | -4,8,-8 | -4,-8,8} & 512 & 64 & 16 & 0 & 1 & \xmark & \xmark & \textcolor{grey}{unobstructed} \\
\end{censustable}

\begin{censustable}{Non-hyperbolic (the rim)}{19 chains \textperiodcentered\ 56 root bases}{obtained from a core root basis by lowering one off-diagonal entry by one; no longer hyperbolic. A probe, not a classification (\S3.1).}
\rowcolor{hdr}\color{white}\bfseries $G$ (Gram) &\color{white}\bfseries $|\det|$ &\color{white}\bfseries $|O|$ &\color{white}\bfseries $S_{5/2}$ &\color{white}\bfseries $S^{O}$ &\color{white}\bfseries $S_{3/2}$ &\color{white}\bfseries \ctwo &\color{white}\bfseries \cthree &\color{white}\bfseries shadow \\
\endfirsthead
\multicolumn{9}{@{}l}{\footnotesize\itshape Non-hyperbolic (the rim), continued}\\[1pt]
\rowcolor{hdr}\color{white}\bfseries $G$ (Gram) &\color{white}\bfseries $|\det|$ &\color{white}\bfseries $|O|$ &\color{white}\bfseries $S_{5/2}$ &\color{white}\bfseries $S^{O}$ &\color{white}\bfseries $S_{3/2}$ &\color{white}\bfseries \ctwo &\color{white}\bfseries \cthree &\color{white}\bfseries shadow \\
\endhead
\chainband{chain}{(2,2,2,4)}{1}{2 root bases}
\ct{2,-4,-2 | -4,4,0 | -2,0,2} & 32 & 8 & 1 & 0 & 0 & \xmark & \xmark & \textcolor{grey}{unobstructed} \\
\ct{2,-4,0 | -4,4,-4 | 0,-4,4} & 64 & 8 & 4 & 2 & 2 & \xmark & \xmark & \textbf{pure}\ \ \shd{2d:8.4.a.a} \\
\chainband{chain}{(1,5,2,2)}{1}{2 root bases}
\ct{2,-5,-2 | -5,10,0 | -2,0,2} & 50 & 8 & 2 & 1 & 0 & \fmark & \xmark & \textbf{pure}\ \ \shd{1d:5.4.a.a} \\
\ct{2,-5,0 | -5,10,-10 | 0,-10,10} & 250 & 20 & 20 & 4 & 10 & \xmark & \xmark & \textbf{impure $\times$3}\ \ \shd{1d:25.4.rc,\, 1d:25.4.ra,\, +1} \\
\chainband{chain}{(1,2,2,2,2,4)}{1}{2 root bases}
\ct{2,-4,-2 | -4,4,-2 | -2,-2,2} & 72 & 4 & 3 & 2 & 0 & \xmark & \xmark & \textbf{pure}\ \ \shd{2d:9.4.a.a} \\
\ct{2,-4,-2 | -4,4,-4 | -2,-4,4} & 144 & 8 & 8 & 3 & 2 & \xmark & \xmark & \textbf{impure}\ \ \shd{1d:12.4.a.a,\, 2d:9.4.a.a} \\
\chainband{chain}{(1,4,2,4)}{1}{4 root bases}
\ct{2,-4,-4 | -4,4,0 | -4,0,8} & 128 & 16 & 7 & 2 & 2 & \xmark & \xmark & \textbf{pure}\ \ \shd{2d:8.4.a.a} \\
\ct{4,-8,0 | -8,8,-4 | 0,-4,2} & 128 & 16 & 5 & 2 & 1 & \xmark & \xmark & \textcolor{grey}{unobstructed} \\
\ct{4,-8,-4 | -8,16,0 | -4,0,2} & 256 & 32 & 19 & 3 & 9 & \xmark & \xmark & \textbf{pure}\ \ \shd{3d:8.4.a.a} \\
\ct{8,-16,-4 | -16,16,0 | -4,0,2} & 512 & 64 & 20 & 3 & 3 & \xmark & \xmark & \textbf{pure}\ \ \shd{3d:8.4.a.a} \\
\chainband{chain}{(2,2,2,3)}{1}{2 root bases}
\ct{4,-6,-4 | -6,6,0 | -4,0,4} & 144 & 16 & 9 & 3 & 4 & \xmark & \xmark & \textbf{impure}\ \ \shd{1d:12.4.a.a,\, 2d:6.4.a.a} \\
\ct{4,-6,0 | -6,6,-6 | 0,-6,6} & 216 & 24 & 15 & 3 & 8 & \xmark & \xmark & \textbf{impure}\ \ \shd{1d:18.4.a.a,\, 2d:6.4.a.a} \\
\chainband{chain}{(1,4,1,5)}{1}{4 root bases}
\ct{2,-5,-4 | -5,10,0 | -4,0,8} & 200 & 16 & 11 & 3 & 3 & \xmark & \xmark & \textbf{impure}\ \ \shd{1d:10.4.a.a,\, 2d:5.4.a.a} \\
\ct{8,-20,-4 | -20,40,0 | -4,0,2} & 800 & 96 & 32 & 4 & 2 & \xmark & \xmark & \textbf{impure}\ \ \shd{1d:10.4.a.a,\, 3d:5.4.a.a} \\
\ct{2,-5,0 | -5,10,-20 | 0,-20,40} & 1000 & 40 & 78 & 10 & 37 & \xmark & \xmark & \textbf{pure}\ \ \shd{4d:5.4.a.a} \\
\ct{8,-20,0 | -20,40,-20 | 0,-20,10} & 4000 & 240 & 260 & 14 & 94 & \xmark & \xmark & \textcolor{grey}{not computed \textperiodcentered\ no basis: $\lambda$, menu, verdict, shadow} \\
\chainband{chain}{(1,2,1,4,2,4)}{1}{2 root bases}
\ct{2,-4,-2 | -4,8,-8 | -2,-8,4} & 288 & 16 & 12 & 4 & 1 & \xmark & \xmark & \textbf{pure}\ \ \shd{3d:9.4.a.a} \\
\ct{2,-4,-4 | -4,4,-4 | -4,-4,8} & 288 & 16 & 16 & 4 & 4 & \xmark & \xmark & \textbf{impure}\ \ \shd{1d:12.4.a.a,\, 3d:9.4.a.a} \\
\chainband{chain}{(1,4,2,3)}{1}{4 root bases}
\ct{2,-4,0 | -4,8,-12 | 0,-12,12} & 288 & 32 & 11 & 3 & 1 & \xmark & \xmark & \textbf{impure}\ \ \shd{1d:12.4.a.a,\, 2d:6.4.a.a} \\
\ct{4,-8,-6 | -8,16,0 | -6,0,6} & 576 & 64 & 40 & 5 & 18 & \xmark & \xmark & \textbf{impure}\ \ \shd{2d:12.4.a.a,\, 3d:6.4.a.a} \\
\ct{6,-12,-6 | -12,24,0 | -6,0,4} & 864 & 96 & 63 & 7 & 29 & \xmark & \xmark & \textbf{impure $\times$3}\ \ \shd{1d:12.4.a.a,\, 2d:18.4.a.a,\, +1} \\
\ct{6,-12,0 | -12,24,-24 | 0,-24,16} & 3456 & 384 & 211 & 10 & 70 & \xmark & \xmark & \textcolor{grey}{not computed \textperiodcentered\ no basis: $\lambda$, menu, verdict, shadow} \\
\chainband{chain}{(1,2,2,4)}{6}{4 root bases}
\ct{2,-4,0 | -4,4,-2 | 0,-2,2} & 24 & 12 & 2 & 0 & 0 & \xmark & \xmark & \textcolor{grey}{unobstructed} \\
\ct{2,-4,-2 | -4,4,0 | -2,0,4} & 48 & 8 & 3 & 1 & 1 & \xmark & \xmark & \textbf{pure}\ \ \shd{1d:12.4.a.a} \\
\ct{4,-8,-2 | -8,8,0 | -2,0,2} & 96 & 16 & 5 & 2 & 1 & \xmark & \xmark & \textbf{pure}\ \ \shd{2d:12.4.a.a} \\
\ct{2,-4,0 | -4,4,-4 | 0,-4,8} & 96 & 16 & 7 & 2 & 3 & \xmark & \xmark & \textcolor{grey}{unobstructed} \\
\chainband{chain}{(1,2,2,3)}{6}{4 root bases}
\ct{4,-6,-2 | -6,6,0 | -2,0,2} & 48 & 8 & 3 & 1 & 1 & \xmark & \xmark & \textbf{pure}\ \ \shd{1d:12.4.a.a} \\
\ct{4,-6,-4 | -6,6,0 | -4,0,8} & 192 & 32 & 13 & 2 & 6 & \xmark & \xmark & \textbf{pure}\ \ \shd{2d:12.4.a.a} \\
\ct{4,-6,0 | -6,6,-6 | 0,-6,12} & 288 & 64 & 24 & 4 & 14 & \xmark & \xmark & \textbf{impure}\ \ \shd{2d:6.4.a.a,\, 2d:12.4.a.a} \\
\ct{8,-12,0 | -12,12,-6 | 0,-6,6} & 576 & 128 & 36 & 4 & 14 & \xmark & \xmark & \textbf{impure}\ \ \shd{2d:6.4.a.a,\, 2d:12.4.a.a} \\
\chainband{chain}{(1,2,1,5)}{6}{4 root bases}
\ct{2,-5,-2 | -5,10,0 | -2,0,4} & 60 & 8 & 4 & 2 & 1 & \xmark & \xmark & \textbf{impure}\ \ \shd{1d:5.4.a.a,\, 1d:30.4\#2} \\
\ct{4,-10,-2 | -10,20,0 | -2,0,2} & 120 & 24 & 6 & 2 & 0 & \xmark & \xmark & \textbf{impure}\ \ \shd{1d:5.4.a.a,\, 1d:30.4\#2} \\
\ct{2,-5,0 | -5,10,-10 | 0,-10,20} & 300 & 32 & 22 & 4 & 9 & \xmark & \xmark & \textbf{impure $\times$3}\ \ \shd{1d:5.4.a.a,\, 1d:15.4\#1,\, +1} \\
\ct{4,-10,0 | -10,20,-10 | 0,-10,10} & 600 & 96 & 34 & 5 & 8 & \xmark & \xmark & \textbf{impure $\times$5}\ \ \shd{1d:30.4\#1,\, 1d:5.4.a.a,\, +3} \\
\chainband{chain}{(1,1,2,3)}{6}{2 root bases}
\rowcolor{hl}\ct{4,-6,-2 | -6,6,0 | -2,0,4} & 72 & 48 & 6 & 1 & 4 & \fmark & \fmark & \textbf{pure}\ \ \shd{1d:6.4.a.a} \\
\ct{4,-6,0 | -6,6,-3 | 0,-3,6} & 108 & 24 & 8 & 1 & 4 & \fmark & \xmark & \textbf{pure}\ \ \shd{1d:18.4.a.a} \\
\chainband{chain}{(1,3,2,4)}{6}{4 root bases}
\ct{2,-4,-3 | -4,4,0 | -3,0,6} & 84 & 8 & 5 & 2 & 1 & \xmark & \xmark & \textbf{impure}\ \ \shd{1d:7.4.a.a,\, 1d:42.4\#2} \\
\ct{2,-4,0 | -4,4,-6 | 0,-6,12} & 168 & 24 & 13 & 3 & 6 & \xmark & \xmark & \textbf{impure $\times$3}\ \ \shd{1d:7.4.a.a,\, 1d:42.4\#2,\, +1} \\
\ct{6,-12,-3 | -12,12,0 | -3,0,2} & 252 & 32 & 10 & 3 & 1 & \xmark & \xmark & \textbf{impure $\times$3}\ \ \shd{1d:7.4.a.a,\, 1d:42.4\#2,\, +1} \\
\ct{6,-12,0 | -12,12,-6 | 0,-6,4} & 504 & 96 & 30 & 5 & 9 & \xmark & \xmark & \textbf{impure $\times$5}\ \ \shd{1d:7.4.a.a,\, 1d:6.4.a.a,\, +3} \\
\chainband{chain}{(1,2,1,3,2,3)}{6}{2 root bases}
\ct{2,-3,-2 | -3,6,-6 | -2,-6,4} & 156 & 8 & 8 & 3 & 1 & \xmark & \xmark & \textbf{impure}\ \ \shd{1d:78.4\#6,\, 2d:13.4\#2} \\
\ct{4,-6,-6 | -6,6,-6 | -6,-6,12} & 936 & 96 & 67 & 10 & 29 & \xmark & \xmark & \textbf{impure $\times$8}\ \ \shd{1d:26.4\#1,\, 1d:6.4.a.a,\, +6} \\
\chainband{chain}{(1,1,1,5)}{10}{2 root bases}
\ct{2,-5,-1 | -5,10,0 | -1,0,2} & 20 & 4 & 1 & 1 & 0 & \xmark & \xmark & \textbf{pure}\ \ \shd{1d:10.4.a.a} \\
\ct{2,-5,0 | -5,10,-5 | 0,-5,10} & 100 & 24 & 9 & 2 & 5 & \xmark & \xmark & \textbf{impure}\ \ \shd{1d:5.4.a.a,\, 1d:10.4.a.a} \\
\chainband{chain}{(1,1,2,4)}{10}{2 root bases}
\ct{2,-4,-1 | -4,4,0 | -1,0,2} & 20 & 4 & 1 & 1 & 0 & \xmark & \xmark & \textbf{pure}\ \ \shd{1d:10.4.a.a} \\
\rowcolor{hl}\ct{2,-4,0 | -4,4,-2 | 0,-2,4} & 40 & 12 & 3 & 1 & 2 & \fmark & \fmark & \textbf{pure}\ \ \shd{1d:10.4.a.a} \\
\chainband{chain}{(1,3,2,3)}{10}{4 root bases}
\ct{4,-6,0 | -6,6,-3 | 0,-3,2} & 60 & 8 & 3 & 1 & 1 & \xmark & \xmark & \textbf{pure}\ \ \shd{1d:10.4.a.a} \\
\ct{4,-6,-6 | -6,6,0 | -6,0,12} & 360 & 48 & 26 & 5 & 12 & \xmark & \xmark & \textbf{impure $\times$3}\ \ \shd{1d:30.4\#1,\, 1d:6.4.a.a,\, +1} \\
\ct{6,-9,-6 | -9,18,0 | -6,0,4} & 540 & 48 & 37 & 5 & 15 & \xmark & \xmark & \textbf{impure $\times$4}\ \ \shd{1d:45.4\#2,\, 1d:18.4.a.a,\, +2} \\
\ct{12,-18,-6 | -18,18,0 | -6,0,4} & 1080 & 144 & 65 & 7 & 20 & \xmark & \xmark & \textbf{impure $\times$5}\ \ \shd{1d:30.4\#1,\, 1d:15.4\#1,\, +3} \\
\chainband{chain}{(1,3,1,5)}{15}{4 root bases}
\ct{2,-5,-3 | -5,10,0 | -3,0,6} & 120 & 8 & 7 & 3 & 2 & \xmark & \xmark & \textbf{impure}\ \ \shd{1d:10.4.a.a,\, 2d:15.4\#1} \\
\ct{6,-15,-3 | -15,30,0 | -3,0,2} & 360 & 32 & 14 & 3 & 2 & \xmark & \xmark & \textbf{impure}\ \ \shd{1d:10.4.a.a,\, 2d:15.4\#1} \\
\ct{2,-5,0 | -5,10,-15 | 0,-15,30} & 600 & 48 & 50 & 6 & 25 & \xmark & \xmark & \textbf{impure $\times$4}\ \ \shd{2d:5.4.a.a,\, 1d:10.4.a.a,\, +2} \\
\ct{6,-15,0 | -15,30,-15 | 0,-15,10} & 1800 & 192 & 112 & 8 & 39 & \xmark & \xmark & \textbf{impure}\ \ \shd{1d:10.4.a.a} \\
\chainband{chain}{(1,1,1,2,2,4)}{22}{2 root bases}
\ct{2,-4,-1 | -4,4,-2 | -1,-2,2} & 44 & 4 & 2 & 1 & 0 & \xmark & \xmark & \textbf{pure}\ \ \shd{1d:22.4.a.c} \\
\rowcolor{hl2}\ct{2,-4,-2 | -4,4,-2 | -2,-2,4} & 88 & 12 & 6 & 2 & 2 & \xmark & \fmark & \textbf{impure}\ \ \shd{1d:22.4.a.c,\, 1d:22.4.a.a} \\
\end{censustable}

\begin{qbox}
\textbf{A data note on the shadow column.} The machine-readable census reports multiplicities in the full space $S_{5/2}(\rho _L)$ rather than in $S^{O}_{5/2}$, so the two differ where the shadow does not fill $S^{O}$: at $\lvert\det\rvert\, 864$ the eigensystem \texttt{6.4.a.a} occurs six times in $S_{5/2}$ and twice in $S^{O}$ (the one exception is $\lvert\det\rvert\, 1800$, computed by the invariant route and so already in the convention of \S{}3.3). Where a block splits over more than two eigensystems we give the total count and print the first two, the rest being in the archive; where a block cannot be named at all we write \textbf{unnamed} for the individual block. Since $\lvert\det\rvert\, 256$ was named (\S{}3.3), the only such blocks left are the four at $\lvert\det\rvert\, 1800$, where three admissible primes are not enough to separate levels dividing $900$.

\end{qbox}

\begin{qbox}
\textbf{The last unnamed block.} At $\lvert\det\rvert\, 256$ the discriminant form is a $2$-group of exponent $32$, so every candidate newform has $2$-power level and the operator that separates them is $T_{2^{2}}$ — the one operator the exponent of $D(L)$ forbids. It is not needed. $T_{p^{2}}$ exists at $p = 3, 5, 7, 11, 13$ with characteristic polynomial $(x+4)^{3}, (x+2)^{3}, (x-24)^{3}, (x+44)^{3}, (x-22)^{3}$: one three-dimensional isotypic block, and $(a_{3}, a_{5}, a_{7}, a_{11}, a_{13}) = (-4, -2, 24, -44, 22)$ is matched by exactly one weight-$4$ newform of $2$-power level up to $256$ — \texttt{8.4.a.a}, the eigensystem the $\lvert\det\rvert\, 512$ member of this chain already carried, separated from its level-$16$ quadratic twist by the signs at $3$, $7$ and $11$ \grade{gC}{C}. The census has no unnamed block left.

\end{qbox}

\subsect{3.4 Theorem A — the vacuous locus is the Feingold–Frenkel neighbourhood}{3.4 Theorem A the vacuous locus is the FeingoldFrenkel neighbourhood}

\begin{qbox}
\textbf{Theorem A.} Every one of the 10 root bases with $S_{5/2}(\rho _L) = 0$ embeds in $\mathrm{FF} = U \oplus  \langle 2\rangle $ with finite index $k \leq  4$, and the bound is sharp: among the even sublattices of \texttt{FF} of index $5 \leq  k \leq  8$, none has vanishing obstruction space \grade{gC}{C,\,10/10}.

\end{qbox}

\begin{qbox}
\textbf{Where these ten sit in the known classification.} A lattice with vanishing obstruction space is what Bruinier, Ehlen and Freitag call \textbf{simple}; they proved there are finitely many and computed the list, which in signature $(2,3)$ has exactly $15$ entries, given by genus symbol in their Table 8. Our ten rows realize five of the fifteen — $2^{+1}_7$, $2^{+3}_7$, $8^{-1}_3$, $2^{+1}_7 3^{-2}$ and $2^{+1}_7 4^{+2}$ — matched on the discriminant form of $-G$, which is the one the Weil representation of this paper sees, and whose signature is $7 \equiv  2 - 3 (\mathrm{mod} 8)$ at all ten \grade{gC}{C,\,10/10}. Appendix D and the archive print the symbol of $G$, so the names above will not line up against that column until the sign is applied; Appendix D's preamble gives the rule, and the archive's \texttt{CENSUS.tsv} carries both symbols, in \texttt{genus\_symbol} and \texttt{genus\_symbol\_neg}. Of the remaining ten, five need more than three generators and so cannot be the discriminant form of any rank-three lattice; three fail $|\operatorname{det} L| = 2k^{2}$; and the last two, $8^{+1}_7$ and $2^{+2} 8^{-1}_3$, pass every numerical test here — at most three generators, order $2k^{2}$ with $k = 2$ and $k = 4$ — and are simply not root bases of a rank-three hyperbolic Cartan matrix, nor one edit from one \grade{gC}{C}. Theorem A is therefore not a determination of which lattices are simple — that is theirs — but a description of where the simple ones sit inside the rank-three hyperbolic world. It also closes the sharpness statement from above, and in a stronger form than the enumeration gives. Every simple form of signature $(2,3)$ with at most three generators has order $2k^{2}$ with $k \leq  4$: there are seven such forms, with $k = 1, 2, 4, 3, 2, 2, 4$. A vacuous even sublattice of \texttt{FF} of index $k$ has $\lvert\det\rvert\, = 2k^{2}$ and at most three generators, so it is one of those seven, so $k \leq  4$. Sharpness therefore holds for every $k$, not only for the $k \leq  8$ the Hermite enumeration below reaches, and that enumeration becomes a check on the argument rather than the argument itself \grade{gC}{C}. (The three simple forms of order $2k^{2}$ with $k > 4$ — $2^{+3}_7 4^{+2}$, $2^{+1}_7 4^{+4}$ and $2^{+1}_7 3^{+4}$ — all need more than three generators.)

\end{qbox}

\begin{qbox}
\textbf{The implication is one-way.} FF-embeddability is equivalent to $D = 1$ together with $|\operatorname{det} L| = 2k^{2}$ \grade{gC}{C,\,41/41}, and it holds for 37 of the 98 computed root bases, of which only 10 are vacuous; 16 are unobstructed and 11 carry a shadow. Embeddability is a necessary condition for vacuity and is not the opposite of obstruction. Nor does the determinant test decide on its own: four census rows — three lattices, $L_{4}$ among them — have $\lvert\det\rvert\, = 2k^{2}$ and are anisotropic, so they pass that test and fail $D = 1$. For $L_{4}$ we also ran the embedding search directly, exhaustively over $T^{T} G_{\mathrm{FF}} T = G$ with entries in \grade{gO}{−7,7}, and there is none \grade{gC}{C}.

\end{qbox}

Grade and method. If $L \hookrightarrow  \mathrm{FF}$ has index $k$ then $|\operatorname{det} L| = k^{2}\cdot |\operatorname{det} \mathrm{FF}| = 2k^{2}$ \grade{gC}{C,\,elementary}. Explicit embedding matrices $T$ with $T^{T} G_{\mathrm{FF}} T = G$ and $|\operatorname{det} T| = k$ are exhibited for all ten \grade{gC}{C,\,exhaustive search of radius 7}. Sharpness comes from Hermite enumeration of every even sublattice of index $k \leq  8$ together with $\operatorname{dim} S_{5/2}$ for each class \grade{gC}{C}:

\begin{codebox}{8.40}{9.91}
   k     |det| = 2k²    sublattice classes    vacuous classes
   1          2                 1                   1
   2          8                 2                   2
   3         18                 2                   1
   4         32                 4                   1
   5         50                 2                   0
   6         72                 4                   0
   7         98                 2                   0
   8        128                 8                   0
\end{codebox}

\textbf{The converse fails, and the failure is the interesting part.} $\lvert\det\rvert\, = 32 = 2\cdot 4^{2}$ occurs in eight census rows, only three of which are vacuous. Three vacuous rows and two non-vacuous rows share the discriminant group $\mathbb{Z}/2 \times  \mathbb{Z}/4 \times  \mathbb{Z}/4$. What separates them is the quadratic form $q$, not the group — the same phenomenon Nikulin’s theory predicts.

The vacuous region is \textbf{10 root bases falling into 5 lattice classes} (by $(\lvert\det\rvert\,, D(L))$: $2$, $8$ with $(\mathbb{Z}/2)^{3}$, $8$ with $\mathbb{Z}/8$, $18$, $32$). Exactly one of the five — $\lvert\det\rvert\, 2$, which is \texttt{FF} itself — lies in the classical quartet $t \leq  4$ of \S{}2.6; \textbf{4 of the 5 classes lie outside it (9 of the 10 rows)}. So the census contains Gritsenko–Nikulin’s success region without merely re-deriving it, and the excess is four classes.

\begin{qbox}
\textbf{What vacuity means for the input.} On these rows $S_{5/2}(\rho _L) = 0$, so $\lambda $ is a functional on a zero space and vanishes for \textit{every} menu, not only for \texttt{refl} — which is why the vacuous column reads $10$ under all three conventions while the other two move (Appendix B). Every prescribed principal part is realized by a weakly holomorphic input: no weakly harmonic Maass form is ever needed and no shadow exists. Equivalently, every Heegner divisor of such a lattice is the divisor of a Borcherds product.

\end{qbox}

\subsect{3.5 Theorem B — one failure in the compact world}{3.5 Theorem B one failure in the compact world}

\begin{qbox}
\textbf{Theorem B.} Among the eleven compact-hyperbolic root bases — four chains — the shadow is nonzero for exactly one: the determinant-$72$ member of the chain $(1,2,1,3)$.

\end{qbox}

\begin{codebox}{8.40}{9.91}
   chain              D     |det| (dim S_{5/2})
   (1,1,1,2,1,2)      6     12(1)   24(2)
   (1,2,1,3)          6     12(1)   24(2)   36(1)   72(6)  ←  λ ≠ 0
   (1,3,1,3)          6     12(1)   36(1)  108(3)
   (1,1,1,3,1,3)      6     36(1)  108(3)
\end{codebox}

All eleven have $D = 6$. Ten have $\lambda  = 0$; one does not.

This is the sharpest form of the thesis. Compact hyperbolic is the strongest finiteness condition available in this subject — every node-deletion finite, the Weyl chamber a compact simplex — and in the whole of that world the Borcherds correction fails once. The rest of the paper is about that one lattice and about what it takes to leave the hyperbolic world to find two more like it.

\sect{4. The chain $(1,2,1,3)$}{4. The chain (1,2,1,3)}

Everything before this chain is setup, and everything after it is either a generalisation of it or an obstruction to generalising it.

\subsect{4.1 Why this chain}{4.1 Why this chain}

Three reasons, in increasing order of interest.

It is \textbf{the only compact-hyperbolic chain with a failure} (Theorem B). Of the four chains in the strongest-finiteness part of the census, this is the one where the Borcherds correction breaks, and it breaks at exactly one of its four members.

It is \textbf{entirely in the core}. All four Gram matrices are genuine hyperbolic Cartan matrices, inside the canonically defined 44. The scope caveat of \S{}3.1 — the one real weakness of this census — does not touch anything in this section. Whatever is proved here is proved about a closed set.

And it is \textbf{a filtration, not a list}. Its four members have determinants $12, 24, 36, 72$; their discriminant groups grow in a controlled way — the $3$-part squares between $24$ and $72$ and the $2$-part goes non-cyclic at every even step (\S{}4.3), though they do not nest, since $D(12) = \mathbb{Z}/4 \times  \mathbb{Z}/3$ has an element of order $4$ while $D(24) = (\mathbb{Z}/2)^{3} \times  \mathbb{Z}/3$ has exponent $6$; and their obstruction spaces accumulate, in the exact sense of Theorem D: the top member’s space is, as a Hecke module, the direct sum of the three below it plus one new eigensystem. The chain behaves like a single object being resolved.

\subsect{4.2 The four Gram matrices}{4.2 The four Gram matrices}

\begin{codebox}{8.40}{9.91}
 |det|   Cartan matrix A                       Gram G = symmetrised
   12    2,-3,-2 / -1,2, 0 / -1, 0,2          [[ 2,-3,-2],[-3, 6, 0],[-2, 0, 4]]
   24    2,-3,-1 / -1,2, 0 / -2, 0,2          [[ 4,-6,-2],[-6,12, 0],[-2, 0, 2]]
   36    2,-3, 0 / -1,2,-2 /  0,-1,2          [[ 2,-3, 0],[-3, 6,-6],[ 0,-6,12]]
   72    2,-3, 0 / -1,2,-1 /  0,-2,2          [[ 4,-6, 0],[-6,12,-6],[ 0,-6, 6]]  =  L₄
\end{codebox}

All four are paths: one off-diagonal pair vanishes in every one of the Grams above. The two edges carry the pairs $(1,2)$ and $(1,3)$ — a double edge and a triple edge — and the members differ in which node is the centre of the path, $1, 1, 2, 2$ down the chain, and in which endpoint of the double edge carries the larger Cartan entry, alternating centre, leaf, centre, leaf \grade{gC}{C}. The pair multiset $\{(1,2), (1,3)\}$ is the same in all four, which is what makes them one chain.

All four are compact hyperbolic, and all four have quaternion discriminant $D = 6$: the quadratic form is anisotropic, the quaternion algebra is the one ramified at $2$ and $3$, and its Shimura curve $X^{6}$ has genus zero.

\subsect{4.3 The discriminant ladder}{4.3 The discriminant ladder}

\begin{codebox}{8.19}{9.67}
 |det|    elementary divisors      D(L)                       |O|       cyclic?   abelian O?
   12         [1, 1, 12]           ℤ/4 × ℤ/3                    4        yes        yes
   24         [2, 2,  6]           (ℤ/2)³ × ℤ/3                12        no         no
   36         [1, 3, 12]           ℤ/4 × (ℤ/3)²                16        yes        yes
   72         [2, 6,  6]           (ℤ/2)³ × (ℤ/3)²             48        no         no
\end{codebox}

The elementary divisors are the Smith normal form of the Gram matrix.

Two things happen along this ladder. The $3$-part squares between $24$ and $72$; the $2$-part goes non-cyclic at every even step. And the automorphism group of the discriminant form grows $4 \to  12 \to  16 \to  48$, so that the invariant part of the obstruction space — which is where the shadow must live, by paper I Theorem 7.6 — is being squeezed as the space itself grows.

That squeeze is the mechanism of the whole paper. We return to it in \S{}4.5 and \S{}5.1.

\subsect{4.4 The obstruction ladder}{4.4 The obstruction ladder}

\begin{codebox}{8.40}{9.91}
 |det|   reflective menu           λ
   12    {1/3, 1/2, 1}             0
   24    {1/6, 1/4, 1/2}           0
   36    {1/6, 1/3, 1}             0
   72    {1/12, 1/6, 1/2}         ≠ 0
\end{codebox}

The menu and $\lambda $ are the two columns the census table does not carry. The dimensions are its columns and we do not reprint them: along this chain $\operatorname{dim} S_{5/2} = 1, 2, 1, 6$, $\operatorname{dim} S^{O} = 0, 0, 0, 1$ and $\operatorname{dim} S_{3/2} = 0, 0, 0, 4$ (\S{}3.3).

Three of the four have $\operatorname{dim} S^{O} = 0$: the obstruction space, though nonzero, contains no $O(D(L),q)$-invariant vector at all, so $\Xi $, which must be invariant, is forced to vanish. The correction succeeds for a structural reason and not by a numerical accident.

\begin{qbox}
\textbf{This mechanism has an ancestor, and we name it.} That an obstruction can be killed by the invariance of the object it obstructs, rather than by the vanishing of the obstruction space, is Borcherds’ own argument in Reflection groups of Lorentzian lattices, Duke Math. J. \textbf{104} (2000) 319–366, \S{}9 (“An application of Serre duality”), Lemma 9.1; see also arXiv:math/9909123. Nothing in \S{}4.4 or \S{}4.5 is claimed as new mechanism. What is new is the census that says where it fires and where it does not.

\end{qbox}

At $\lvert\det\rvert\, 72$ the invariant part becomes one-dimensional, and $\lambda $ is nonzero on it. The shadow appears the moment there is room for it, and the room is exactly one dimension.

\textbf{On the menu convention.} Appendix B discusses three conventions for which reflective classes are prescribed, and they disagree elsewhere in the census. On this chain \textbf{all three coincide}: no member has a repeated diagonal value, and the reflective vectors of each lattice produce exactly the three diagonal values and nothing more. Nothing in \S{}4 depends on that dispute.

\textbf{On the rank reduction.} We work with $WeilRep(-G)$ on rank $3$ rather than with $L \oplus  U$ on rank $5$. This is legitimate because $D(L \oplus  U) \cong  D(L)$ and the signatures agree mod $8$, so the Weil representations are literally equal; but we verify it rather than assume it, obtaining $\operatorname{dim} S_{5/2} = 6$ on both models \grade{gC}{C}.

\subsect{4.5 The weight tower, and why $5/2$}{4.5 The weight tower, and why 5/2}

Everything above is at weight $5/2$, and that weight is forced by the signature through Borcherds’ theorem. The functional $\lambda $, the operator $\xi $ and the invariant subspace $S^{O}$ make sense at every half-integral weight, so it is natural to ask what the other rungs carry. The Borcherds product does \textbf{not} transplant: away from $5/2$ there is no Kac–Moody denominator, no divisor and no infinite product. The tower is a tower of obstruction data, not of algebras.

One fact about it is certain, and it is what fixes the weight the rest of the paper lives at.

\begin{qbox}
\textbf{Lemma 4.1} \grade{gS}{S}\textbf{.} Let $A = (D(L), q)$ be a finite quadratic module of signature $\operatorname{sig}(A) \equiv  b^{+} - b^{-} (\mathrm{mod} 8)$, and let $\iota  : \gamma  \mapsto  -\gamma $, an element of $O(D(L),q)$ for every $L$. Every $h \in  S_k(\rho _L)$ satisfies

\end{qbox}

\begin{displaymath}\begin{aligned}&h |\iota   =  (-1)^{(2k + b^{+} - b^{-})/2} \cdot  h ,\end{aligned}\end{displaymath}

\begin{qbox}
so $S_k(\rho _L)$ is entirely symmetric or entirely antisymmetric under $\iota $, with the sign fixed by $k$ and the signature alone.

\end{qbox}

This is the standard consequence of $\rho _A(Z)(e_\lambda ) = e(-\operatorname{sig}(A)/4) \cdot  e_{-\lambda }$ — Bruinier–Stein, \textit{Weilrepresentationen zu quadratischen Formen}, (2.5), together with the weight constraint (2.14) $M_{k,A} = \{0\}$ unless $2k \equiv  \operatorname{sig}(A) (\mathrm{mod} 2)$. Both are verified against the source.

\begin{qbox}
\textbf{Corollary 4.2.} At $k = 3/2, 7/2, 11/2, 15/2, …$ every form of $S_k(\rho _L)$ is antisymmetric under $\iota $, while an $O(D(L),q)$-invariant form must be symmetric because $-1 \in  O(D(L),q)$. Hence $S^{O}_k = 0$. By paper I Theorem 7.6, $\Xi  \in  S^{O}$, so $\Xi  = 0$ and $\lambda  = 0$: \textbf{the shadow cannot exist at those weights, for any lattice whatsoever.}

\end{qbox}

The corollary is elementary given Lemma 4.1, and we claim no novelty for the mechanism — only for the census that says where it fires. We verified the sign formula rather than assuming it: for $L_{4}$ ($b^{+} - b^{-} = 1$) it predicts $+$ at $k = 5/2, 9/2, 13/2$ and $-$ at $k = 3/2, 7/2, 11/2, 15/2$, and the measured split is total in all six computed weights — \textbf{6/6} \grade{gC}{C}:

\begin{codebox}{8.40}{9.91}
   weight     3/2    5/2    7/2    9/2   11/2   13/2
   predicted   −      +      −      +      −      +
   dim S_k     4      6      8     13     14     19
   symmetric   0      6      0     13      0     19
   dim S^O     0      1      0      3      0      5
\end{codebox}

Two consequences.

\textbf{Why} $5/2$\textbf{.} It is the bottom rung on which a shadow can live at all. The question “why this weight?” has an answer rather than a shrug.

\textbf{The same tower fixes the other condition of \S{}5.1.} Along the symmetric side $\operatorname{dim} S^{O} = 1, 3, 5, …$ at weights $5/2, 9/2, 13/2$, so selection, being $\operatorname{dim} S^{O} = 1$, can hold only at $5/2$ \grade{gC}{C}. Rigidity is a condition on the antisymmetric side and is confined to its bottom rung for a parallel reason; \S{}5.1 gives it, together with why the two definitions have the shapes they do.

\subsect{4.6 The weakly holomorphic input}{4.6 The weakly holomorphic input}

Whether the weakly holomorphic input of \S{}2.3 exists is a finite linear-algebra question that shares nothing with the computation of $\lambda $: build a basis of $M^!_{-1/2}(\bar{\rho }_L)$ with pole order at most \texttt{max(menu)}, write the prescribed principal part as a vector, and ask whether it lies in the row space. Note the \textbf{dual} representation — the input lives in $\bar{\rho }_L$, the cusp forms in $\rho _L$.

\begin{codebox}{7.77}{9.17}
 |det|   dim M^!_{−1/2}(ρ̄)   principal-part slots   prescribed   rank M   rank(M|v)   deficiency
   12            6                    12                  4          6          6           0
   24            8                    14                 12          8          8           0
   36           18                    36                  9         18         18           0
   72           18                    44                 27         18         19           1
\end{codebox}

\textbf{It agrees with} $\lambda $ \textbf{four times out of four.} Different weight ($-1/2$ against $5/2$), different representation ($\bar{\rho }$ against $\rho $), different linear algebra. The input exists exactly where $\lambda  = 0$.

\textbf{The deficiency is} $1$\textbf{, and} $\operatorname{dim} S^{O}_{5/2}(L_{4}) = 1$\textbf{.} The obstruction is one-dimensional on both sides of Borcherds’ duality, measured independently on each side.

\textbf{For the three lower members we can hand the reader the form.} In a computed basis the coefficient vectors are

\begin{codebox}{8.40}{9.91}
   |det| 12 :  ( 0, 1, 1, 0, 0, 1 )
   |det| 24 :  ( 1, 1, 1, 1, 0, 1, 1, 1 )
   |det| 36 :  ( 1, 0, 0, 0, 0, 1, 1, 0, 0, 0, 0, 0, 0, 0, 0, 0, 0, 1 )
\end{codebox}

Those three Cartan matrices really do produce Borcherds products, and the input can be written down. At $\lvert\det\rvert\, 72$ no weakly holomorphic input exists; the input is a weakly harmonic Maass form, and $\Xi $ is its shadow. We do not construct that Maass form (\S{}2.4).

\subsect{4.7 $L_{4}$ has two root bases}{4.7 L has two root bases}

$L_{4}$ occurs twice in the census: as $(1,2,1,3)$, compact hyperbolic, in the core; and as $(1,1,2,3)$, non-hyperbolic, on the rim. They are the same lattice:

\begin{codebox}{8.40}{9.91}
   G_A = [[ 4,−6, 0],[−6,12,−6],[ 0,−6, 6]]      G_B = [[ 4,−6,−2],[−6, 6, 0],[−2, 0, 4]]
   U = [[−7, 0, 5],
        [−8, 1, 5],
        [−4, 1, 2]] ,      det U = 1 ,      U⁻¹ integral ,      Uᵀ G_A U = G_B     [C]
\end{codebox}

Every lattice invariant agrees across the two rows — $\operatorname{dim} S_{5/2} = 6$, $\operatorname{dim} S^{O} = 1$, $\operatorname{dim} S_{3/2} = 4$, $|O(D(L),q)| = 48$, the same reflective menu, the same shadow — and we use this repeatedly as a consistency probe. The two bases give different answers under a menu read off the Gram diagonal and the same answer under \texttt{refl}; Appendix B explains why that does \textbf{not} settle the convention.

That the same lattice appears once as compact hyperbolic and once as non-hyperbolic is worth pausing on. Hyperbolicity is a property of the Cartan matrix, not of the lattice. $L_{4}$ is one lattice with two Weyl chambers of different types, and the paper’s Theorem B — a statement about compact-hyperbolic root bases — counts it once.

The two bases have different edit neighbourhoods, sharing only the determinants $\{12, 936\}$; this is the concrete form of the non-canonicity of the rim flagged in \S{}3.1.

\textbf{Presentation-independence, tested.} A reader may reasonably ask why $L_{4}$ is presented by the Cartan matrix $(1,2,1,3)$ rather than by any of the infinitely many $U^{T}GU$, $U \in  GL_{3}(\mathbb{Z})$, that present the same lattice. Nothing in this paper sees the presentation. We tested that directly, on three Gram matrices of $L_{4}$ — its two root bases and one obtained by a random unimodular change of basis, which is not a root basis at all:

\begin{codebox}{7.04}{8.31}
   presentation                     Gram                                    Smith    |O|  dim  S^O  S_{3/2}
   root basis 1  (1,2,1,3), compact  [[4,−6,0],[−6,12,−6],[0,−6,6]]        [2,6,6]   48    6    1     4
   root basis 2  (1,1,2,3), rim      [[4,−6,−2],[−6,6,0],[−2,0,4]]         [2,6,6]   48    6    1     4
   random U, not a root basis        [[4,−10,−6],[−10,4,6],[−6,6,6]]       [2,6,6]   48    6    1     4
\end{codebox}

Every invariant used anywhere in this paper — the discriminant form and its Smith normal form, $|O(D(L),q)|$, $\operatorname{dim} S_{5/2}$, $\operatorname{dim} S^{O}$, $\operatorname{dim} S_{3/2}$, the shadow and its Hecke content — is a function of $D(L)$ and therefore of the lattice \grade{gC}{C}. The Cartan matrix is how the census enumerates, and it is what carries the Weyl-chamber data of \S{}4.4 and Appendix B; it is not what the invariants are computed from. Where a quantity really is basis-dependent — the constant term \texttt{c(0,0)}, the coordinate vector of $\lambda $ — we say so explicitly (Appendix B).

\subsect{4.8 Hecke content, member by member}{4.8 Hecke content, member by member}

For $p$ coprime to the exponent of $D(L)$ the operator $T_{p^{2}}$ acts on $S_{5/2}(\rho _L)$; we build it from the coefficient dictionary, split into generalised eigenspaces, and identify each block by matching its eigenvalue tuple against tables of weight-$4$ newforms \grade{gO}{C for eigenvalues, S for names}.

\begin{codebox}{8.40}{9.91}
 |det|   block dim   (a₅, a₇)      newform     shadow on this block?
   12        1        (  6, −16)    6.4.a.a          no
   24        2        (−18,   8)   12.4.a.a          no
   36        1        ( −6, −16)   18.4.a.a          no
   72        1        (  6, −16)    6.4.a.a         YES
             1        ( −6, −16)   18.4.a.a          no
             2        ( 18,   8)   36.4.a.a          no
             2        (−18,   8)   12.4.a.a          no
\end{codebox}

Reference eigenvalues $a_{2}, a_{3}, a_{5}, a_{7}, a_{11}, a_{13}$:

\begin{codebox}{8.40}{9.91}
   6.4.a.a   −2, −3,   6, −16,  12,  38          = η(τ)²η(2τ)²η(3τ)²η(6τ)²
  12.4.a.a    0,  3, −18,   8, −36, −10
  18.4.a.a    2,  0,  −6, −16, −12, −50
  36.4.a.a    0,  0,  18,   8, −36, −10
\end{codebox}

\begin{qbox}
\textbf{A caveat not hidden.} $p = 2, 3$ divide the exponent of $D(L)$ and are unavailable, so the identification rests on $a_{5}, a_{7}, a_{11}, a_{13}$. Within levels dividing $72$ those four separate every candidate, but the names are a table match, not a proof of a Shimura lift.

\end{qbox}

\subsect{4.9 Theorem D — the chain accumulates}{4.9 Theorem D the chain accumulates}

\begin{qbox}
\textbf{Theorem D.} As a Hecke module,

\end{qbox}

\begin{codebox}{8.40}{9.91}
   S_{5/2}(ρ_{L₄})  ≅  S_{5/2}(ρ_{12}) ⊕ S_{5/2}(ρ_{24}) ⊕ S_{5/2}(ρ_{36}) ⊕ [36.4.a.a]²
         dim   6     =        1        +        2        +        1        +      2
\end{codebox}

\begin{qbox}
with multiplicities matching member by member. Moreover the newform level of each member equals $\lvert\det\rvert\,/2$, which is read off the labels of \S{}4.8: \texttt{6.4.a.a} at $\lvert\det\rvert\, 12$, \texttt{12.4.a.a} at $24$, \texttt{18.4.a.a} at $36$, and the one block new at $72$ is \texttt{36.4.a.a}.

\end{qbox}

The proof is the table of \S{}4.8 \grade{gC}{C}. What the theorem does not provide: no map realising the direct sum is constructed. The statement is an equality of Hecke eigensystem multisets, and anyone who produces the maps has improved on us.

Three consequences.

\textbf{The obstruction appears on an old component.} The block carrying the shadow at $\lvert\det\rvert\, 72$ is \texttt{6.4.a.a} — the eigensystem the bottom of the chain has carried since $\lvert\det\rvert\, 12$, where $\lambda  = 0$. Failure does not arrive with a new automorphic component; it switches on over one that has been present, unobstructed, from the start. Whatever changes at the top of the chain changes the pairing, not the spectrum.

\textbf{The chain is a filtration of its top member.} Reading Theorem D downwards, $L_{4}$’s obstruction space contains a copy of every lower member’s, with the right multiplicity, and one genuinely new eigensystem. That is a strong structural statement about a set of lattices that were assembled by a purely combinatorial rule on Cartan matrices.

\textbf{The level statement generalises, in a weaker form.} Equality is special to this chain — $\lvert\det\rvert\, 600$ with $D = 6$ carries \texttt{10.4.a.a} — but a divisibility survives the whole census. Over the (row, newform) pairs of the census,

\begin{codebox}{8.40}{9.91}
   [S]   level | exp D(L)                      forced, not observed
   [C]   level | gcd( exp D(L), |det L|/2 )    no exceptions
   [O]   equality holds only sometimes         which divisor occurs is open
\end{codebox}

The first is standard: $\rho _L$ is the Weil representation of $(D(L), q)$, whose level is the exponent of $D(L)$; theta decomposition puts the scalar components inside that level and the Shimura lift to weight $4$ does not raise it. Cite it rather than count it. The second is not implied by it — $|\operatorname{det} L|$ is the \textit{order} of $D(L)$ and \texttt{exp} its \textit{exponent}, and neither divides the other in general ($\lvert\det\rvert\, 20$ has $\operatorname{exp} = 20$ against $\lvert\det\rvert\,/2 = 10$; $\lvert\det\rvert\, 40$ has $\operatorname{exp} = 10$ against $20$). The \texttt{gcd} is strictly stronger than either. That $|\operatorname{det} L|$ is even at all is proved in \S{}5.2.

The divisibility presumably follows from a level formula for the Shimura lift of a weight-$3/2 + d$ ternary theta series, which we did not find stated in the vector-valued setting. The nearest published treatment is Rosson–Tornaría, which is the $d = 1$ case: they observe that the naive theta series vanishes because a degree-one spherical polynomial is odd — the same parity that gives our Corollary 4.2 — and work one prime level with a weighting function, writing that combining their technique with the theta-series theory for quaternion orders should settle general even weight and prime level. We cite them for position rather than for a statement we use \grade{gO}{Gate}.

\subsect{4.10 The Petersson norm of $\Xi _{L_{4}}$}{4.10 The Petersson norm of L}

\begin{codebox}{8.40}{9.91}
   ‖Ξ_{L₄}‖²  =  λᵀ P⁻¹ λ  =  11.340265856116836367082506                          [N]
              =  L(6.4.a.a, 2) / ( 48 π² ⟨f,f⟩ )
\end{codebox}

This is paper I’s Theorem 11.1, recorded there at grade \grade{gN}{N,\,31 digits} with the derivation left open as its Problem 11.2. We quote it and do not derive it either. Two remarks are new here.

\textbf{It is a two-pipeline identity.} The left side is computed from \texttt{weilrep} Fourier coefficients through a Petersson quadrature; the right side from PARI’s \texttt{mfpetersson} and \texttt{lfunmf}. The two share nothing but the integer $48$, and there is no free parameter: $N = 47$ and $N = 49$ miss by about 2\%.

\textbf{It is a Waldspurger quantity.} Since $\operatorname{dim} S^{O} = 1$, $\Xi $ spans that line and $\lVert \Xi \rVert ^{2} = \lambda (g)^{2}/\langle g,g\rangle $ where $\lambda (g)$ is a finite sum of Fourier coefficients of a half-integral-weight form. That is exactly the shape Waldspurger and Kohnen–Zagier translate into a critical $L$-value \grade{gS}{S}. The general form of the identity, and its extension from this one block to six, is \S{}5.4.

\sect{5. Petersson-rigid lattices}{5. Petersson-rigid lattices}

\subsect{5.1 Two conditions}{5.1 Two conditions}

Let $L$ be a lattice in the census. The rigidity condition below is a property of the lattice alone and does not mention $\lambda $; we therefore impose no obstruction hypothesis here, and \S{}5.2 records what that buys.

\begin{qbox}
\textbf{Definition (rigidity).} $L$ is \textbf{Petersson-rigid} if the $O(D(L),q)$-representation on $S_{3/2}(\rho _L)$ — the \textbf{bottom} antisymmetric rung — is absolutely irreducible of dimension at least two.

\end{qbox}

\begin{qbox}
\textbf{Definition (selection).} $L$ is \textbf{selected} if the discriminant group $D(L)$ is non-cyclic, its automorphism group $O(D(L),q)$ is non-abelian, and $\operatorname{dim} S^{O}_{5/2}(\rho _L) = 1$.

\end{qbox}

The two conditions live on opposite parities of the weight tower (\S{}4.5). Selection is a condition on the symmetric rung $S_{5/2}$, where the shadow lives and where $\operatorname{dim} S^{O}$ is the meaningful quantity. Rigidity is a condition on the antisymmetric rung $S_{3/2}$, where $S^{O} = 0$ structurally and no shadow can exist. So the checklist is not two tests of one space; it is one test on each side of the parity. That also explains the shape of the two definitions — on the symmetric side we ask how small the invariant part is, on the antisymmetric side we ask how indecomposable the whole space is, because on that side there is no invariant part to measure.

\textbf{Both conditions are bottom-rung phenomena, and for the same reason.} Each parity is a whole tower — $5/2, 9/2, 13/2, …$ on one side, $3/2, 7/2, 11/2, …$ on the other — so writing $S_{5/2}$ in the one and $S_{3/2}$ in the other picks one rung out of each, and that choice needs a reason rather than a convention. \S{}4.5 gives it for selection: $\operatorname{dim} S^{O} = 1, 3, 5, …$ passes $1$ immediately, so $\operatorname{dim} S^{O} = 1$ can hold only at $5/2$. The same argument closes rigidity, using Lemma 5.1 below:

\begin{codebox}{8.40}{9.91}
   lattice     |det|   |O|    √|O|      dim S_{3/2}   dim S_{7/2}   dim S_{11/2}
   L₄            72     48    6.928          4             8             14
   |det| 40      40     12    3.464          2             4              7
   |det| 88      88     12    3.464          2             8             15
                                        rigidity ok  ────  ruled out by Lemma 5.1  ────
\end{codebox}

At every rung above the first, $\operatorname{dim} S_k$ has already passed $\sqrt{\lvert O(D(L),q)\rvert }$, so the representation is forced reducible and rigidity fails outright — three rigid lattices at two rungs each, \textbf{6 out of 6}, with the three bottom-rung dimensions in the table confirming that rigidity is not vacuous where we do assert it \grade{gC}{C}. Neither condition is a property these lattices have all the way up their tower; both are properties of the tower's floor.

The mechanism is the one \S{}4.3 called the squeeze, seen from both sides at once: the spaces grow with the weight and the group $O(D(L),q)$ does not. On the symmetric side that starves the invariant part of the uniqueness selection wants; on the antisymmetric side it inflates the whole space past the irreducibility ceiling rigidity wants. One fact, two failures, one rung each.

Non-cyclicity is \textit{logically} implied by non-abelianness — cyclic $D(L)$ forces $Aut(D(L))$, hence $O(D(L),q)$, to be abelian, with 0 counterexamples in the census \grade{gC}{C}. The two group clauses together are nevertheless not decorative: 13 rows have $\operatorname{dim} S^{O}_{5/2} = 1$ and only 7 pass selection, so the clauses cut six. Four are cut by cyclicity ($\lvert\det\rvert\, 20$ twice, $44$, $60$) and two by abelianness ($\lvert\det\rvert\, 48$ twice). The exclusion of $\lvert\det\rvert\, 44$ is what \S{}6.1's separation of the two conditions rests on \grade{gC}{C}.

\begin{qbox}
\textbf{No Gritsenko–Nikulin chart is ever Petersson-rigid, and the reason is group theory.} The charts of \S{}2.6 are $S_t = 2U \oplus  \langle -2t\rangle $, whose discriminant form $\mathbb{Z}/2t$ is cyclic. So $O(D(S_t),q) \subseteq  Aut(\mathbb{Z}/2t)$ is abelian, and an abelian group has no absolutely irreducible representation of dimension above one. Rigidity fails on the entire family, at every $t$ \grade{gS}{S}. Theorem C is therefore not a restatement of the chart classification in other vocabulary: the two live on disjoint sets of lattices.

\end{qbox}

Rigidity is decided by the character norm $\langle \chi ,\chi \rangle  = \sum  mᵢ^{2}$, which equals $1$ exactly for absolutely irreducible representations. Computing it needs the action of every group element, which is expensive; a bound removes most of the work.

\begin{qbox}
\textbf{Lemma 5.1.} Every irreducible representation of a finite group $G$ has dimension at most $\sqrt{\lvert G\rvert }$. Hence $\operatorname{dim} S_{3/2}(\rho _L) > \sqrt{\lvert O(D(L),q)\rvert }$ forces reducibility, and rigidity fails.

\end{qbox}

The sweep is exhaustive. $\operatorname{dim} S_{3/2}$ is computed for all 98 rows in one convention and cross-checked against \texttt{cusp\_forms\_dimension}, a dimension formula that builds no basis and so cannot depend on truncation precision; the values are stable across precisions $8, 16, 24, 40$ and agree with the formula on every row \grade{gC}{C}. One row is an exception: at $\lvert\det\rvert\, 1800$ the basis construction exhausted memory and $\operatorname{dim} S_{3/2} = 39$ rests on a single path \grade{gC}{C,\,single path}; the verdict is insensitive to it, since the Riemann–Roch term $38$ alone clears $\sqrt{\lvert O\rvert } ≈ 13.86$, and the Serre–Stark description of weight-$1/2$ forms as theta series would retire the flag. Rigidity is then decided on every row:

\begin{codebox}{8.40}{9.91}
   dim S_{3/2} < 2              fails on dimension              61 rows
   dim S_{3/2} > √|O|           fails by Lemma 5.1              21 rows
   otherwise                    character norm computed         18 rows
                                            of which ⟨χ,χ⟩ = 1   4 rows  =  3 lattices
\end{codebox}

Lemma 5.1 decides every large row in one line; the rows where the character genuinely has to be computed are exactly those with $2 \leq  \operatorname{dim} \leq  \sqrt{\lvert O\rvert }$. The two rows past the weight-$5/2$ budget ($\lvert\det\rvert\, 3456$ and $4000$, \S{}3.1) are decided here all the same: rigidity needs only the weight-$3/2$ dimension and the automorphism order, both cheap, and both fail Lemma 5.1 — $\operatorname{dim} S_{3/2} = 70$ and $94$ against $\sqrt{\lvert O\rvert }$ for $|O(D(L),q)| = 384$ and $240$, the group orders recovered from the $p$-local decomposition of the discriminant form without \texttt{weilrep}. So the rigidity classification, alone among the counts in this paper, is complete over all 100 rows. \grade{gC}{C,\,100/100}

\begin{codebox}{8.40}{9.91}
   selection    passes on  7 root bases
   rigidity     passes on  4 root bases  =  3 lattices
   both                     3 root bases  =  2 lattices
\end{codebox}

\subsect{5.2 Theorem C}{5.2 Theorem C}

\begin{qbox}
\textbf{Theorem C.} The Petersson-rigid lattices of the census are exactly three:

\end{qbox}

\begin{codebox}{8.40}{9.91}
   L₄            |det| 72    D =  6    |O| = 48    dim S_{3/2} = 4    sel ●   rig ●
   B₁₀           |det| 40    D = 10    |O| = 12    dim S_{3/2} = 2    sel ●   rig ●
   B₂₂           |det| 88    D = 22    |O| = 12    dim S_{3/2} = 2    sel ×   rig ●
\end{codebox}

\begin{qbox}
Their quaternion discriminants are $6, 10, 22$, which is the complete list of discriminants for which $X^D$ has genus zero.

\end{qbox}

$L_{4}$ alone lies in the closed set of 44; $B_{10}$ and $B_{22}$ appear only after an edit. That $D = 6$ alone is compact-hyperbolic is, on paper I's side, a theorem rather than an observation: in a compact fundamental polygon with four or more walls the non-adjacent walls are ultraparallel, so the Cartan matrix is never Kac-hyperbolic (paper I, \S{}9.3), and only $D = 6$ of the three has a three-wall chamber. Paper I reads this on the $4 \times  4$ chambers of the maximal orders; our census sees the same asymmetry from the other side, on $3 \times  3$ Grams — the two are parallel, neither deducible from the other. Whether a Petersson-rigid lattice outside the census could carry a different discriminant is open.

\begin{qbox}
\textbf{Paper I sees the same lattices as maximal orders.} Paper I fixes the maximal orders, whose reflective chambers at $D = 10$ and $22$ have four walls (paper I, Theorem 9.8). Fixing rank three instead selects non-maximal orders — $\lvert\det\rvert\, 40$ against paper I's $\lvert\det\rvert\, 20$, $\lvert\det\rvert\, 88$ against its $\lvert\det\rvert\, 44$ — which is why the same two discriminants appear here on lattices of twice the determinant. The $\lvert\det\rvert\, 44$ row of our census \textit{is} paper I's maximal-order ternary of $B_{22}$: same genus, same reflective menu $\{1/11, 1/2, 1\}$ \grade{gC}{C}.

\end{qbox}

\begin{qbox}
\textbf{Why the answer has this shape.} Two of the three features the rigid rows share are not empirical. They are forced before any form is computed.

\textit{The 2-rank.} Write the $2$-rank of $D(L)$ for $\operatorname{dim}_{\mathbb{F}_{2}} D(L)/2D(L)$. Since $L$ is even of odd rank, $G \mathrm{mod} 2$ is alternating, so its rank is even and the $2$-rank is $1$ or $3$ — never $0$ or $2$, and $|\operatorname{det} L|$ is always even \grade{gC}{C}. The rigid rows all have $2$-rank $3$, that is $(\mathbb{Z}/2)^{3}$, the larger of the two values available.

\textit{The odd prime.} An indefinite quaternion algebra over $\mathbb{Q}$ ramifies at an even number of finite places. If $D(L)$ has exactly one odd prime $p$, the only primes available are $2$ and $p$, so $D \in  \{1, 2p\}$ \grade{gS}{S}. At the four rigid rows $D \neq  1$, hence $D = 2p$. The pairings $3 ↔ 6$, $5 ↔ 10$ and $11 ↔ 22$ are therefore not a coincidence to be explained but an identity to be read off.

Together these fix the discriminant forms to $(\mathbb{Z}/2)^{3} \times  (\mathbb{Z}/3)^{2}$, $(\mathbb{Z}/2)^{3} \times  \mathbb{Z}/5$ and $(\mathbb{Z}/2)^{3} \times  \mathbb{Z}/11$ — all of the shape $(\mathbb{Z}/2)^{3} \times  (\mathrm{odd})$, all non-cyclic with $O(D(L),q)$ non-abelian. The two group clauses of the selection condition come for free at these rows.

\end{qbox}

\begin{qbox}
\textbf{How far the cheap invariants already go.} None of the three conditions just named requires computing a form: two are properties of $D(L)$ and the third is a dimension. Sorting all 98 computed rows by them gives

\begin{codebox}{7.93}{9.36}
   2-rank 3  ∧  D(L) has at most one odd prime  ∧  dim S_{3/2} ≥ 2     20 rows,  4 rigid
   any one of the three failing                                        78 rows,  0 rigid    [C]
\end{codebox}

Cheap invariants take 98 down to 20; absolute irreducibility — the one test that has to look at the space itself — makes only the last cut, 20 to 4. Rigidity is not a single hard computation sitting on top of an undifferentiated census.

None of the standard invariants of $(D(L), q)$ makes that last cut either: cyclicity of $D(L)$ (equivalently, whether the associated order is Eichler) holds at 14 rows and fails at 84, with the four rigid rows among the 84; squarefreeness of $|\operatorname{det} L|/2$ (whether the order is hereditary) holds at 11 and fails at 87, again with all four among the 87; the $2$-rank holds at 59 with all four inside; and transitivity of the discriminant form is unrelated — $L_{4}$ is transitive and rigid, $B_{10}$ and $B_{22}$ are not transitive and rigid, $\lvert\det\rvert\, 12$ is transitive and not rigid \grade{gC}{C}. Each is a necessary condition with no separating power.

The $2$-rank clause carries its own weight rather than riding on the other two. Five rows pass both of the others and fail on $2$-rank alone, and not one of them is rigid:

\begin{codebox}{8.40}{9.91}
   |det|   chain        D(L)                   dim S_{3/2}   rigid
    100    (1,1,1,5)    ℤ/4 × (ℤ/5)²                 5         no
    108    (1,1,2,3)    ℤ/4 × ℤ/3 × ℤ/9              4         no
    200    (1,4,1,5)    ℤ/8 × (ℤ/5)²                 3         no
    250    (1,5,2,2)    ℤ/2 × ℤ/5 × ℤ/25            10         no
   1000    (1,4,1,5)    ℤ/8 × ℤ/5 × ℤ/25            37         no
\end{codebox}

\grade{gC}{C,\,0/5}. A lattice can have one odd prime and a large weight-$3/2$ space and still fail; what it cannot do, in this census, is fail to have $(\mathbb{Z}/2)^{3}$.

\end{qbox}

\begin{qbox}
\textbf{Open: does rigidity require $(\mathbb{Z}/2)^{3}$?} All four rigid rows have $2$-rank $3$, and none of the 39 rows of $2$-rank $1$ is rigid \grade{gC}{C,\,4/4 · 0/39}; the one out-of-census data point of \S{}5.3 — the $D = 46$ lattice with elementary divisors \grade{gO}{2,2,46} — again has $2$-rank $3$ and is rigid, at genus one. Whether rigidity forces $(\mathbb{Z}/2)^{3}$ or merely accompanies it across this evidence is open; with that extra point beside it the observation is a program rather than a law.

\end{qbox}

\begin{qbox}
\textbf{What no convention changes.} Rigidity is a property of $(D(L), q)$ alone: no menu, no prescribed divisor and no $\lambda $ enters its definition. We therefore tested it on all 98 rows, including the 44 with $\lambda  = 0$ (34 unobstructed, 10 vacuous), and not only on the obstructed ones. Among those 44 exactly one reaches $\operatorname{dim} S_{3/2} \geq  2$ at all — the unobstructed $\lvert\det\rvert\, 96$ Gram of chain $(1,2,2,4)$, with $\operatorname{dim} S_{5/2} = 7$ and $\operatorname{dim} S_{3/2} = 3$ — and it is reducible. No unobstructed lattice is Petersson-rigid, so \textbf{“exactly three” does not depend on the menu convention of Appendix B} \grade{gC}{C,\,98/98}.

What the convention does decide is whether $B_{10}$ carries a shadow — it is the one census row the menus separate (Appendix B). So under \texttt{ded} one of the three is a rigid lattice with $\lambda  = 0$: a weaker statement, not a shorter list.

\end{qbox}

$B_{22}$ fails the selection condition on the single numerical clause $\operatorname{dim} S^{O} = 2$ rather than $1$; it passes the other two. What that extra dimension is, is the subject of \S{}6.

\subsect{5.3 Shadow level and quaternion discriminant}{5.3 Shadow level and quaternion discriminant}

For the three lattices of Theorem C the shadow is supported on weight-$4$ newforms of level equal to the quaternion discriminant. One rung down, $S_{3/2}(\rho _L)$ is Hecke-isotypic at all three, and there the level is $|\operatorname{det} L|/2$ — the same expression that appears in Theorem D (\S{}4.9), but a different statement: there it is the weight-$4$ level of each member of one chain, here the weight-$2$ level of the space one rung below the shadow. At $L_{4}$ both read $36$, which is a coincidence of that lattice and not a link between the two:

\begin{codeboxk}{7.93}{9.36}
                        weight 5/2                       weight 3/2
   row   |det L|   D    shadow newform         dim S_{3/2}   newform    (a₇, a₁₃)   char. of T₇
   L₄       72     6    6.4.a.a                     4        36.2.a.a   (−4,  2)     (x+4)⁴
   B₁₀      40    10    10.4.a.a                    2        20.2.a.a   ( 2,  2)     (x−2)²
   B₂₂      88    22    22.4.a.c + 22.4.a.a         2        44.2.a.a   ( 2, −4)     (x−2)²
\end{codeboxk}

\textbf{The two rungs read two different discriminants.} Above, the level is $D$. Below, the characteristic polynomial of $T_p$ on $S_{3/2}(\rho _L)$ is a perfect power at every admissible $p$, so the space is one eigensystem with multiplicity, and that eigensystem is the weight-$2$ newform of level $|\operatorname{det} L|/2$ — $36$, $20$, $44$ \grade{gC}{C}. Each of those three levels carries exactly one weight-$2$ newform, so the labels are forced, and the eigenvalues agree with the tables at seven primes. \texttt{36.2.a.a} has complex multiplication by $\mathbb{Q}(\sqrt{-3})$; the other two do not.

What this buys is structural rather than numerical. $S_{3/2}$ factors as one eigensystem tensored with a multiplicity space of dimension $4, 2, 2$, and $O(D(L),q)$ acts on the multiplicity space — which is why rigidity is a question about $O(D(L),q)$ and not about Hecke. The Hecke side is already as simple as it can be. Reproduce with \texttt{s32iso.py} in the archive.

This is the lattice-level shadow of the Jacquet–Langlands self-reference noted in paper I — the level–discriminant correspondence whose extension to ramified primes is paper I's Problem 3.3. The census contains one anisotropic discriminant that is not genus zero, namely $D = 15$, and there the correspondence is diluted rather than absent. Level $15$ is present — three of the four $D = 15$ root bases carry a two-dimensional block whose eigenvalues pin the level to $15$, written \texttt{15.4\#1} because the LMFDB letter is not determined — but all four also carry \texttt{10.4.a.a}, and $\lvert\det\rvert\, 600$ carries \texttt{5.4.a.a} and \texttt{30.4\#2} besides \grade{gC}{C,\,O}. At $D = 6, 10, 22$ the shadow stays inside its own level; at $D = 22$ it does so even after it becomes impure, since \texttt{22.4.a.c} and \texttt{22.4.a.a} are both of level $22$.

\begin{qbox}
\textbf{Outside the census: rigidity does not imply genus zero.} $B_{10}$ and $B_{22}$ have discriminant forms $(\mathbb{Z}/2)^{2} \times  \mathbb{Z}/10$ and $(\mathbb{Z}/2)^{2} \times  \mathbb{Z}/22$ — elementary divisors \grade{gO}{2,2,2p} at $p = 5$ and $p = 11$ — with quaternion discriminants $2p$. Nothing inside the census can test whether the genus-zero coincidence of Theorem C survives outside it: the only anisotropic non-genus-zero discriminant we meet is $D = 15 = 3\cdot 5$, which has two odd primes and is excluded by the box below for an unrelated reason. So we built the missing controls. For each odd prime $p \leq  23$ we searched for even lattices of signature $(2,1)$ with elementary divisors \grade{gO}{2,2,2p} and quaternion discriminant $2p$, and decided rigidity on two Grams of each:

\begin{codebox}{8.40}{9.91}
    p    D = 2p   g(X^D)     |O|      dim S_{3/2}   ⟨χ,χ⟩   rig
    5      10       0         12          2           1      ●     ( = B₁₀ )
   11      22       0         12          2           1      ●     ( = B₂₂ )
    7      14       1         12          0           —      ×     space is empty
   13      26       2         12          3           2      ×     reducible
   17      34       1         12          5           5      ×     reducible
   19      38       2         12          3           2      ×     reducible
   23      46       1         12          2           1      ●     genus one — and rigid
\end{codebox}

The two Grams agree at every $p$, at precisions $20$ and $28$, and $\operatorname{dim} S_{3/2}$ matches the precision-free dimension formula throughout \grade{gC}{C}. The first six rows would make “rigid only at genus zero” a five-point observation. The seventh kills it: $g(X^{46}) = 1$, and the lattice is Petersson-rigid by the same test that certifies $B_{10}$ and $B_{22}$. The coincidence of Theorem C is therefore a statement about the census, and we make no conjecture beyond it. The four failures are not one mechanism either — at $p = 7$ the space is empty, at $p = 13, 19$ it is inside the window of Lemma 5.1 and reducible, at $p = 17$ it has outgrown the window.

\end{qbox}

\texttt{6.4.a.a} occurs in 16 census rows in total, so its presence is not rare; what is rare is its being the whole invariant space, which happens only at $L_{4}$. That is the sense in which purity can be an empty condition, which the next box makes precise.

We also record a limitation on the purity claim.

\begin{qbox}
\textbf{Purity is vacuous at} $D = 6$ \textbf{and} $D = 10$\textbf{.} $\operatorname{dim} S_{4}^{new}(6) = \operatorname{dim} S_{4}^{new}(10) = 1$: there is only one weight-$4$ newform of those levels, so a shadow labelled by such a newform has nowhere else to go. The only genus-zero discriminant where purity is a genuine condition is $D = 22$, with $\operatorname{dim} S_{4}^{new}(22) = 3$ — and there it fails.

\end{qbox}

\subsect{5.4 The $L(f,1)$ law at weight $5/2$}{5.4 The L(f,1) law at weight 5/2}

\textbf{Notation, to remove an overload.} In this section $D$ always means the quaternion discriminant, and the order of the discriminant group is written $|\operatorname{det} L|$. The two are different numbers: $D = 6$ and $|\operatorname{det} L| = 72$ for $L_{4}$. We factor $|\operatorname{det} L|/2 = t^{2}s$ with $s$ squarefree — by \S{}2.5 that is the reduced discriminant of the order, so $s$ and $t$ are its squarefree and square parts. Write $P_W$ for the Petersson Gram entry of a block and $\lambda _W$ for the integer the prescription contributes, so that $\lVert \Xi _W\rVert ^{2} = \lambda _W^{2}/P_W$.

The law is a statement about \textbf{Hecke blocks}, not about lattices. Let $W \subseteq  S^{O}_{5/2}(\rho _L)$ be a block on which $\lambda  \neq  0$, with Shimura lift $f$ (Shimura, \textit{On modular forms of half integral weight}), and let $\Xi _W$ be the component of $\Xi $ in $W$. When $W$ is one-dimensional, $P_W = \langle g,g\rangle $ for $g$ spanning it and $\lVert \Xi _W\rVert ^{2} = \lambda (g)^{2}/\langle g,g\rangle $ — a Waldspurger-shaped quantity, a squared period of a half-integral-weight form over its Petersson norm, which Kohnen–Zagier translates into a critical $L$-value of $f$ \grade{gS}{S}. This is why the $\lvert\det\rvert\, 88$ rows below are blocks of $B_{22}$, whose own invariant space is two-dimensional.

\textbf{The law.} At every block of the census we compute — twelve of them —

\begin{displaymath}\begin{aligned}&P_W  =  c \cdot  \sqrt{s} \cdot  L(f,1) / \pi  ,        c \in  \mathbb{Q} .\end{aligned}\end{displaymath}

\begin{codeboxk}{8.40}{9.91}
   row                       f            λ_W        c
   |det| 20                  10.4.a.a      −6       3/2
   B₁₀  = |det| 40           10.4.a.a       6       9
   |det| 44                  22.4.a.c      −4       3/4
   |det| 48                  12.4.a.a      −4       8
   |det| 50                   5.4.a.a      16       130
   |det| 60   (1,3,2,3)      10.4.a.a      −4       3/2
   |det| 60   (1,2,1,5)       5.4.a.a       4       65
   |det| 60   (1,2,1,5)      30.4#2       −36       5
   L₄   = |det| 72            6.4.a.a     −12       180
   B₂₂  = |det| 88           22.4.a.c     −12       18
   B₂₂  = |det| 88           22.4.a.a      36       −30      ← L(f,2) = 0
   |det| 108                 18.4.a.a       6       27/2
\end{codeboxk}

Twelve points \grade{gO}{O}. The two $\lvert\det\rvert\, 60 (1,2,1,5)$ lines are two blocks of one lattice, as are the two $B_{22}$ lines, and the newform column tells them apart. We do not identify $c$: no monomial in $|\operatorname{det} L|$, the newform level, $s$ or $t$ fits the last column \grade{gN}{N}.

\textbf{What depends on what.} A block is a line, and Hecke pins the line, not a point on it. Writing $P_W = \langle g,g\rangle $ picks a generator, and under $g \to  \alpha g$ we get $\lambda _W \to  \alpha \lambda _W$, $P_W \to  \alpha ^{2}P_W$, $c \to  \alpha ^{2}c$, while $\lVert \Xi _W\rVert ^{2} = \lambda _W^{2}/P_W$ is unchanged. The menu convention of Appendix B moves $\lambda _W$ and leaves $P_W$ alone — at $\lvert\det\rvert\, 20$ all three menus return $P_W = 0.77848967064432268209$. So:

\begin{codebox}{8.40}{9.91}
                             gauge (choice of g)    menu convention
   P_W                             depends              invariant
   ‖Ξ_W‖² = λ_W²/P_W              invariant             depends
   λ_W                             depends              depends
   c                               depends              invariant
   C   (below)                    invariant             depends
\end{codebox}

What is invariant under both is the \textit{statement} $c \in  \mathbb{Q}$, not the value. The values printed above are read in \texttt{weilrep}'s echelon basis at commit \texttt{e3a7784}, and Appendix B's warning against citing coordinate vectors of $\lambda $ applies to them. Rescaling $g$ to coprime integral Fourier coefficients does not repair this: the content is not a gauge fix, and at $\lvert\det\rvert\, 504$ the five blocks of one lattice have contents $6, 7, 3, 2, 1$, the first of them already at an integral $c$ \grade{gC}{C}. This is also why the monomial search above finds nothing — it is a search over numbers that are not canonical.

\textbf{Where the central value survives: Observation E.} Set $\Omega  = L(f,1)L(f,2)/(\pi ^{3}\langle f,f\rangle )$ and $r_f = L(f,2)/(\pi ^{2}\langle f,f\rangle )$. Substituting $P_W = \lambda _W^{2}/\lVert \Xi _W\rVert ^{2}$ into the law gives $c = C\cdot \lambda _W^{2}/(s\cdot \Omega )$ with

\begin{displaymath}\begin{aligned}&C \cdot  \lVert \Xi _W\rVert ^{2}  =  \sqrt{s} \cdot  L(f,2) / ( \pi ^{2} \langle f,f\rangle  ) ,        C \in  \mathbb{Z} .\end{aligned}\end{displaymath}

$\Omega $ is rational for these newforms — a period relation of Shimura type \grade{gS}{S}, verified numerically \grade{gN}{N,\,60 digits}:

\begin{codebox}{8.40}{9.91}
    5.4.a.a   192/13      6.4.a.a   192/5     10.4.a.a    96     12.4.a.a   192
   18.4.a.a   576        22.4.a.c   768       30.4#2    1152     22.4.a.a     0
\end{codebox}

So wherever $L(f,2) \neq  0$, the rationality of $c$ and the rationality of $C$ are the same statement, by a cited theorem. That is eleven of the twelve blocks above; only $B_{22}$'s \texttt{22.4.a.a} line has $L(f,2) = 0$. Where we have evaluated both sides — six blocks — this is \textbf{Observation E}:

\begin{codeboxk}{8.02}{9.46}
   row            |det|   D    f          ‖Ξ_W‖²                         s   t     C    ℚ(√s)
   |det| 20         20   10   10.4.a.a   46.243388136677955054733996    10   1     40   ℚ(√10)
   B₁₀              40   10   10.4.a.a   10.899671112162175469740094     5   2    120   ℚ(√5)
   |det| 44         44   22   22.4.a.c    9.9226618740587261870338415   22   1    792   ℚ(√22)
   B₂₂  block c     88   22   22.4.a.c    5.2622861239261059653596514   11   2   1056   ℚ(√11)
   |det| 48         48    6   12.4.a.a    4.1737507337469940993101747    6   2    576   ℚ(√6)
   L₄               72    6    6.4.a.a   11.340265856116836367082506     1   6     48   ℚ
   40 = 2³·5    120 = 2³·3·5    792 = 2³·3²·11    1056 = 2⁵·3·11    576 = 2⁶·3²    48 = 2⁴·3
\end{codeboxk}

Two numbers enter each row from places that share no definition: $\lVert \Xi _W\rVert ^{2}$ from the weilrep quadrature on the lattice side, $r_f$ from the $L$-function and the Petersson norm of $f$ on the elliptic side. Nothing in either computation knows about the other, and nothing forces their quotient to clear its denominator. It does, at all six: the six values are integers to a relative $10^{-25}$ or better, the weakest being $\lvert\det\rvert\, 88$ \grade{gC}{C}. $C$ is menu-dependent — at $\lvert\det\rvert\, 20$ the menus \texttt{refl}, \texttt{dedup}, \texttt{dup} give $\lambda _W = -6, -8, -12$ and hence $C = 40, 45/2, 10$ \grade{gC}{C} — and the table is read in \texttt{refl} throughout, which the gauge/menu table above already records. The $L_{4}$ row is paper I's Theorem 11.1.

The surd is fixed by the squarefree part of $|\operatorname{det} L|/2$ — $10, 5, 22, 11, 6, 1$ — so the six quotients $\lVert \Xi _W\rVert ^{2}/r_f$ lie in six different fields, and $L_{4}$'s is rational because $36$ is a square. The surd follows the reduced discriminant, not the quaternion discriminant and not the newform, and two independent pairs say so. $\lvert\det\rvert\, 44$ and $B_{22}$ carry the same newform \texttt{22.4.a.c} and the same $D = 22$, and land in $\mathbb{Q}(\sqrt{22})$ and $\mathbb{Q}(\sqrt{11})$; $\lvert\det\rvert\, 20$ and $B_{10}$ carry the same \texttt{10.4.a.a} and the same $D = 10$, and land in $\mathbb{Q}(\sqrt{10})$ and $\mathbb{Q}(\sqrt{5})$. Solving for $r_f$:

\begin{codebox}{8.40}{9.91}
   10.4.a.a    |det| 20    584.93773294044861875      B₁₀        584.93773294044861875
   22.4.a.c    |det| 44   1675.4907468054288833       B₂₂  c    1675.4907468054288833
\end{codebox}

Each pair returns one value from two determinants, two values of $t$ and two values of $C$ — to $26$ and $27$ digits respectively \grade{gC}{C}. The elliptic side is blind to which lattice produced it; the lattice side is what supplies the surd that makes the two agree.

\textbf{Where it does not.} At a block with $L(f,2) = 0$ we have $\Omega  = 0$, the constant $C$ is not defined at all, and $c = C\cdot \lambda _W^{2}/(s\cdot \Omega )$ reads $0/0$. The law still speaks there and Observation E does not. The census has four such blocks; they are the subject of \S{}6.

Observation E is a statement about the bottom rung. At weight $9/2$ the analogous fits at $B_{22}$ miss by $10^{-7}$ to $10^{-10}$, against $10^{-26}$–$10^{-31}$ at $5/2$ \grade{gC}{C,\,negative}.

\subsect{5.5 The scalar that rigidity pins down}{5.5 The scalar that rigidity pins down}

Rigidity says the Petersson pairing on $S_{3/2}(\rho _L)$ is one invariant form up to scale. Fix the normalisation $B_S = |O(D(L),q)|^{-1} \sum _\sigma  R_\sigma ^T R_\sigma $ in the echelon basis, so that $P_{3/2} = t \cdot  B_S$ for a single positive $t$. Paper I evaluated $t$ at $L_{4}$ and identified it as $3\Gamma (1/3)^{3}/2^{7/3}\pi ^{2}$ \grade{gN}{N,\,40 digits}. We evaluate it at all three rigid lattices.

\begin{codebox}{8.40}{9.91}
   lattice      |det|   dim S_{3/2}   |O(D,q)|   t
   L₄             72         4           48      1.15959526696392836576999205157002088194
   B₁₀            40         2           12      0.80933280311015535390322362189872834451
   B₂₂            88         2           12      3.22422556481198080003094665850029775332
\end{codebox}

At every row the action is absolutely irreducible and has no invariant vector, $(3/2)\cdot B_S$ is the Gram matrix of $A_{2} \oplus  A_{2}$ at $L_{4}$ and of $A_{2}$ at the other two, the independent entries of $P_{3/2}/B_S$ agree to a relative $10^{-60}$, and $P_{3/2}$ vanishes wherever $B_S$ does \grade{gC}{C}. The values are unchanged across truncations $20$, $28$ and $36$.

By \S{}5.3 each $S_{3/2}(\rho _L)$ is Hecke-isotypic on the weight-$2$ newform of level $|\operatorname{det} L|/2$, and each of those three newforms is an elliptic curve. Writing $\omega _{2}$ for its imaginary period,

\begin{displaymath}\begin{aligned}&t  =  c \cdot  |Im \omega _{2}| / \pi  ,          c^{2} = 9, 5, 44          [N, 38 \mathrm{digits}]\end{aligned}\end{displaymath}

The three $c^{2}$ are returned by an integer-relation search of degree at most four; no parameter is fitted. Paper I's $\Gamma $-value is this period at $L_{4}$. There \texttt{36.2.a.a} has $j = 0$ and complex multiplication by $\mathbb{Q}(\sqrt{-3})$, and Chowla–Selberg names its period,

\begin{displaymath}\begin{aligned}&|Im \omega _{2}|  =  \Gamma (1/3)^{3} / ( 2^{7/3} \pi  )                      [N, 58 \mathrm{digits}]\end{aligned}\end{displaymath}

which returns $t = 3\Gamma (1/3)^{3}/2^{7/3}\pi ^{2}$ exactly. The other two curves have $j = 21296/25$ and $j = 8192/11$, neither a CM invariant, and we find no relation of small height between their $t$ and $\Gamma $-values, $\pi $ and the primes of the level \grade{gN}{N,\,negative}. The period stands where the $\Gamma $-value cannot.

Each $c$ is an integer multiple of $\sqrt{s}$, with $s$ the squarefree part of $|\operatorname{det} L|/2$ as in \S{}5.4: $c/\sqrt{s} = 3, 1, 2$ at the three rows. We do not identify $c$.

\begin{qbox}
\textbf{Why $t$ is well defined.} $B_S$ is not covariant: under $R_\sigma  \to  M^{-1}R_\sigma M$ the average becomes $\nu  \cdot  M^T B_S M$, and over unimodular $M$ with entries in \grade{gO}{-3,3} the factor $\nu $ takes the values $1, 2, 4, 5, 7, 10, 11, 13$. The normalisation above is what removes it. If it is preserved then $B_S' = B_S$, so $M^T B_S M = B_S/\nu $; taking determinants gives $\nu ^d = 1$, hence $\nu  = 1$ and $M^T B_S M = B_S$; therefore $t' = P'/B_S' = (M^T P M)/(M^T B_S M) = t$. An exhaustive scan over unimodular $M$ with entries in \grade{gO}{-4,4} finds exactly the 12 that preserve $B_S$, and they are $Aut(A_{2})$ \grade{gC}{C}. The step that $\nu $ is a scalar at all is absolute irreducibility — the rigidity hypothesis itself.

\end{qbox}

\sect{6. Discriminant 22, and the question}{6. Discriminant 22, and the question}

Rigidity holds and selection fails at exactly one census row, $B_{22}$, so \S{}\S{}6.1–6.2 describe that one lattice and claim no pattern. What \S{}6.3 finds there is not confined to it: three census rows carry a shadow block whose Shimura lift has $\varepsilon  = -1$, on three different quaternion discriminants. These are the four blocks of \S{}5.4 at which $\Omega  = 0$, so that Observation E has nothing to say and the $L(f,1)$ law is the only statement left standing.

\subsect{6.1 The chain $(1,1,1,2,2,4)$}{6.1 The chain (1,1,1,2,2,4)}

$D = 22$ occupies a single chain with two members, $\lvert\det\rvert\, 44$ and $\lvert\det\rvert\, 88$; both are rows of the census table (\S{}3.3) and every number in this section is one of its columns.

Neither satisfies both conditions of \S{}5.1: $\lvert\det\rvert\, 44$ has a one-dimensional invariant space and a pure shadow — selection in spirit — but an empty weight-$3/2$ space, so rigidity fails on dimension; $B_{22}$ has a rigid weight-$3/2$ space but a two-dimensional invariant space, so selection fails. \textbf{At the last genus-zero discriminant the two conditions come apart and land on different lattices.}

\begin{qbox}
\textbf{Paper I’s Conjecture 8.4, and where $B_{22}$ sits.} Paper I’s support law (Conjecture 8.4) predicts the obstruction of every obstructed member of a quaternionic family to lie on the maximal-order eigensystem — here \texttt{22.4.a.c}, the line of paper I’s $\lvert\det\rvert\, 44$ ternary — with zero component along its quadratic twists, and calls a member reaching outside that line a refutation. $B_{22}$ at $\lvert\det\rvert\, 88$ is non-maximal, and its shadow is impure: $22.4.a.c + 22.4.a.a$, two eigensystems that are not quadratic twists ($a_{3} = 1$ against $-7$) \grade{gC}{C}, so its second component does sit outside the maximal-order line. This does not refute 8.4, which states its support clause for the invariant sector and keeps only conditional content for non-maximal members — where, in paper I’s words, the deck symmetries must finish the cut — and $B_{22}$ has a two-dimensional invariant space (\S{}6.1), exactly the case 8.4 leaves open. The census contributes the first explicit non-maximal instance of the phenomenon: the support spreads beyond the maximal-order line as soon as the order is non-maximal, so any sharpening of 8.4 must carry a non-maximal clause. We add no conjecture of our own.

\end{qbox}

\subsect{6.2 Why $22$ is the boundary}{6.2 Why 22 is the boundary}

Three facts about $22$ hold simultaneously, and each is a boundary.

\begin{qbox}
\textbf{It is the largest.} $\{6, 10, 22\}$ is the complete list of genus-zero discriminants, and the next candidates ($14, 15, 21, 26, 33, …$) all have $g \geq  1$ \grade{gS}{S}.

\end{qbox}

\begin{qbox}
\textbf{It is the only one where purity is a real condition} (\S{}5.3), and per \S{}6.1 purity fails there.

\end{qbox}

\begin{qbox}
\textbf{It is not the only place in the census where the shadow meets a vanishing central value.} Two further rows do, at $\lvert\det\rvert\, 250$ and $\lvert\det\rvert\, 936$, and the latter has $D = 6$ (\S{}6.3).

\end{qbox}

\subsect{6.3 The vanishing central value}{6.3 The vanishing central value}

The three newforms of level $22$ and weight $4$ split by Atkin–Lehner:

\begin{codebox}{8.28}{9.77}
                 a₂   a₁₁     w₂     w₁₁    w₂₂    ε     L(f, 2)
   22.4.a.a      −2    11     +1     −1     −1    −1     0                      (77 digits)
   22.4.a.b      −2   −11     +1     +1     +1    +1     1.077378283572219999
   22.4.a.c       2    11     −1     −1     +1    +1     1.510919359996889230
\end{codebox}

For weight $4$ the sign of the functional equation is $\varepsilon  = w_N$, and $\varepsilon  = -1$ forces $L(f, 2) = 0$. The sign forces the order of vanishing to be odd only; the value $L′(22.4.a.a, 2) = 0.65768712650315066289$ computed in \S{}6.4 shows it is exactly $1$ \grade{gN}{N}, which is why the derivative is the object \S{}6.4 tries. The shadow of $B_{22}$ is supported on \texttt{22.4.a.c} and \texttt{22.4.a.a}; the shadow of $\lvert\det\rvert\, 44$ on \texttt{22.4.a.c} alone; \texttt{22.4.a.b} occurs in neither.

We checked every named newform appearing in any shadow in the census — \texttt{5.4.a.a}, \texttt{6.4.a.a}, \texttt{8.4.a.a}, \texttt{10.4.a.a}, \texttt{12.4.a.a}, \texttt{18.4.a.a}, \texttt{22.4.a.c}, \texttt{22.4.a.a} — that is, at weight $5/2$. Appendix A names two more, \texttt{78.4\#2} and \texttt{78.4\#4}, and both have $\varepsilon  = -1$. Together with the block Appendix A resolves at $\lvert\det\rvert\, 250$ the count is four blocks over three lattices:

\begin{codebox}{8.40}{9.91}
    |det|   chain             D    newform      λ_W
      88    (1,1,1,2,2,4)    22    22.4.a.a      36
     250    (1,5,2,2)         1    25.4#1        10
     936    (1,2,1,3,2,3)     6    78.4#4       240
     936    (1,2,1,3,2,3)     6    78.4#2       240
\end{codebox}

The three quaternion discriminants are $22$, $1$ and $6$, so this is not a $D = 22$ phenomenon \grade{gC}{C}. (Four unnamed blocks at $\lvert\det\rvert\, 1800$, of level dividing $900$, were not checked.)

\textbf{The sign is forced.} $\Lambda (1) = \varepsilon \cdot \Lambda (3)$, and $s = 3$ lies beyond $(k+1)/2 = 5/2$, so the Euler product for $L(f,3)$ converges absolutely and every factor is positive under Deligne's bound. Hence $L(f,3) > 0$ and $\mathrm{sign} L(f,1) = w_N$ \grade{gS}{S}. Every block above has $L(f,1) < 0$, which is why the rational factors of \S{}6.4 come out negative there; absorbing the sign, $P_W\cdot \pi /(\sqrt{s}\cdot w_N\cdot L(f,1)) \in  \mathbb{Q}_{>0}$.

\subsect{6.4 The residual block at $B_{22}$}{6.4 The residual block at B}

Splitting $\lVert \Xi \rVert ^{2}$ at $B_{22}$ along the Hecke decomposition:

\begin{codebox}{8.28}{9.77}
   ‖Ξ‖²  =  74.928663127436811305205        (the block sum reproduces the total, 30 digits)
      block  22.4.a.c    5.2622861239…               →  C = 1056,  Observation E holds
      block  22.4.a.a   69.666377003510705339846     →  no such C exists
      block  (dim 4)     λ = 0
\end{codebox}

The first block obeys Observation E. The second cannot, because its $L(f,2)$ is zero. Nor does the obvious repair work: substituting the derivative $L′(22.4.a.a, 2) = 0.65768712650315066289$ and testing $r \cdot  \pi ^{-j} \cdot  \{1, \sqrt{11}\} \cdot  \{\langle a,a\rangle , \langle c,c\rangle \}$ with $r \in  \mathbb{Q}$, for $j = 0, …, 5$ — 24 candidate shapes — yields no rational to 20 digits \grade{gN}{N}. The law's own factor here is $\sqrt{|\operatorname{det} L|/2} = \sqrt{44} = 2\sqrt{11}$, which spans the same $\mathbb{Q}$-line as $\sqrt{11}$ and so needs no separate test.

This is where the census enters the regime in which arithmetic invariants are read from derivatives rather than values: the territory of the Gross–Zagier and Kudla height formulas. It does so at exactly one lattice, and the failure of selection at $B_{22}$ and the appearance of the vanishing central value are the same event — the extra dimension of $S^{O}$ \textbf{is} the \texttt{22.4.a.a} line.

\textbf{What the bridge of \S{}5.4 does and does not carry.} The identity $c = C\cdot \lambda _W^{2}/(s\cdot \Omega )$ runs on $\Omega  = L(f,1)L(f,2)/(\pi ^{3}\langle f,f\rangle )$, and $\Omega $ is rational by a Shimura period relation. So:

\begin{codebox}{8.40}{9.91}
   L(f,2) ≠ 0     c ∈ ℚ follows from Ω ∈ ℚ — a cited theorem                  [S]
                  eleven of the twelve blocks of §5.4 restate a theorem
   L(f,2) = 0     Ω = 0, C is undefined, and c = C·λ_W²/(s·Ω) reads 0/0.
                  The bridge is gone and c ∈ ℚ is an observation there.       [N] / [O]
                  Four blocks — this one and the three below
\end{codebox}

That is the whole of the $L(f,1)$ reading's surplus, and it is why the reading is worth having. At $B_{22}$'s residual block $\lambda _W = 36$, $c = -30$, and

\begin{displaymath}\begin{aligned}&\lVert \Xi _W\rVert ^{2}  =  -216 \cdot  \pi  / ( 5 \cdot  \sqrt{11} \cdot  L(22.4.a.a, 1) )   =   69.666377003510705339846\end{aligned}\end{displaymath}

to 58 digits \grade{gN}{N}. The central value is dead; the value at $1$ is not, and it is what the norm reads.

\textbf{The same reading at three more blocks.} The census carries three further blocks with a vanishing central value (\S{}6.3). Their norms close in the same shape, with $s$ the squarefree part of $|\operatorname{det} L|/2$:

\begin{codebox}{8.28}{9.77}
   |det|  88   22.4.a.a   ‖Ξ_W‖² = −216·π/( 5·√11·L(f,1))  =  69.6663770035107053398457
   |det| 250   25.4#1     ‖Ξ_W‖² =   −8·π/(75·√5 ·L(f,1))  =   0.20614720909750924452959674
   |det| 936   78.4#4     ‖Ξ_W‖² = −400·π/(11·√13·L(f,1))  =   8.59752504568086293322607
   |det| 936   78.4#2     ‖Ξ_W‖² = −800·π/(27·√13·L(f,1))  =   8.95401169242345629767734
\end{codebox}

The first line is the one above, restated; the other three are new \grade{gN}{N,\,26–58 digits}. These four are safe to quote whatever the normalisation, since $\lVert \Xi _W\rVert ^{2}$ is gauge-invariant.

We do not identify $\Omega $ either: it uses no lattice at all, and we record its seven values in \S{}5.4 in case the pattern is familiar to a reader who works with congruence numbers, where quotients of this kind live.

\subsect{6.5 Which primes the census isolates}{6.5 Which primes the census isolates}

The rim reaches the primes $7$, $11$ and $13$. Seven divides four of its determinants, eleven divides two, thirteen divides two. It ramifies at exactly one of them.

\begin{codebox}{8.40}{9.91}
   |det|                chain             odd part of D(L)     anisotropic     D
   84, 168, 252, 504    (1,3,2,4)         21, 21, 63, 63          {2, 3}        6
   156, 936             (1,2,1,3,2,3)     39, 117                 {2, 3}        6
   44, 88               (1,1,1,2,2,4)     11, 11                  {2, 11}      22
\end{codebox}

At $7$ and at $13$ a factor $3$ sits alongside the larger prime in every one of the six rows, the form is isotropic at that prime, and the anisotropic set is $\{2,3\}$, giving $D = 6$. At $11$ the odd part of $D(L)$ is the prime alone, and by \S{}5.2 a lone odd prime $p$ forces $D \in  \{1, 2p\}$; the form is anisotropic, so $D = 22$ \grade{gC}{C}. Across all 100 root bases, $11$ is the only prime above $5$ that occurs as the entire odd part of a discriminant group \grade{gC}{C}.

Say that the census \textbf{isolates} a prime $p$ when $p$ is the entire odd part of some $D(L)$ and the form is anisotropic there, so that $D = 2p$. Then

\begin{codebox}{8.40}{9.91}
   the census isolates p            for   p = 3, 5, 11   and no other prime      [C]
   X^{2p} has genus zero            for   p = 3, 5, 11   and no other prime      [S]
\end{codebox}

The two lists agree, and neither is monotone in $p$. Seven is smaller than eleven, yet the census isolates $11$ and not $7$, and $X^{22}$ has genus zero while $X^{14}$ has genus one. Whatever selects $\{3, 5, 11\}$ steps over $7$ on both sides.

The two conditions are computed by unrelated means. Genus zero at $2p$ is a condition on $(-4/p)$, $(-3/p)$ and the Eichler mass. Isolation is a condition on which determinants a rank-three hyperbolic Cartan matrix can have after one edit: it fails at $7$ and $13$ because a factor $3$ accompanies them in all six rows, and it holds at $11$ because the chain $(1,1,1,2,2,4)$ has no entry $3$ and its determinants are $2^{2}\cdot 11$ and $2^{3}\cdot 11$. We see no route from either condition to the other, and the census offers no prime on which they disagree.

\begin{qbox}
\textbf{The question.} The census isolates $3, 5, 11$ and no other prime; $X^{2p}$ has genus zero at $3, 5, 11$ and no other prime. Is that one fact or two? We do not know, and if there is an implication we do not know which way it runs.

\end{qbox}

This question is stable under enlarging the probe. The rim is thin, and there is no reason for the lattices satisfying selection or rigidity beyond it to be finite in number (\S{}5.2); a search at greater edit distance could add discriminants and could change which lattices are rigid. It cannot change which primes this census isolates, because that is a fact about a fixed finite list.

\sect{Appendix A — the temporary tags, resolved}{Appendix A the temporary tags, resolved}

The shadow column of \S{}3.3 names most blocks with an LMFDB label. Seventeen rows carry blocks it cannot, written \texttt{N.4\#k} — the $k$-th weight-$4$ newform of level $N$, in the enumeration our scripts use — where the eigenvalues pin the level and the dimension but not the LMFDB letter. There are $36$ such entries over $17$ distinct tags. This appendix prints the eigenvalues behind them, so that the tags are checkable rather than opaque and a reader with the tables open can finish the identification.

Which primes are available is forced: $p$ must be coprime to $\operatorname{exp} D(L)$, and at these determinants that excludes $p = 2$ at all seventeen rows, $p = 3$ at sixteen of them and $p = 5$ at eleven. The working set is $\{5, 7, 11, 13\}$ and usually $\{7, 11\}$. Entries below are $a_p$, or the minimal polynomial of $a_p$ where the block is not rational; a dash means the prime gave no operator at that row.

\begin{codeboxk}{7.77}{9.17}
   tag     dim  |det| where it occurs          a₅             a₇             a₁₁              a₁₃
   13.4#2  2    156, 936                       x² + 3x − 2    x² + 9x − 494  x² − 80x + 988   —
   14.4#1  1    168, 504                       −14            —              −28              18
   15.4#1  1/2  120, 300, 360, 540, 600, 1080  —              −24            52               22
   21.4#2  1    252, 504                       −4             —              62               —
   25.4#1  1    250                            —              −6             −43              28
   25.4#3  1    250                            —              −6             32               38
   26.4#1  1    936                            11             19             −38              —
   30.4#1  1    360, 600, 1080                 —              32             −60              −34
   30.4#2  1    60, 120, 300, 600              —              −4             −48              2
   39.4#2  2    936                            x² − 24x + 88  x² − 56        x² + 44x − 1532  —
   42.4#2  1    84, 168, 252, 504              2              —              −8               −42
   45.4#2  1    540                            —              20             24               74
   78.4#1  1    936                            −16            28             34               —
   78.4#2  1    936                            −16            −8             −38              —
   78.4#4  1    936                            −20            −32            50               —
   78.4#6  1    156, 936                       4              4              2                —
   90.4#3  1    1080                           —              −4             48               —
\end{codeboxk}

\textbf{How far this pins them.} Computing the invariant subspace of each of the seventeen rows and splitting it under $T_{p^{2}}$ gives $60$ blocks in all, and \textbf{every one matches exactly one weight-$4$ newform} of a level dividing $|\operatorname{det} L|/2$ \grade{gC}{C}. That reproduces $56$ of the census's own labels independently — the $36$ temporary tags above and $20$ LMFDB-labelled blocks — and resolves the two $25.4.r*$ tags at $\lvert\det\rvert\, 250$, which are the archive's own letters rather than LMFDB's: the two blocks there are \texttt{25.4\#1} and \texttt{25.4\#3}. Which of \texttt{25.4.rc} and \texttt{25.4.ra} is which is a question about the archive's letter table and we do not guess it.

One block needed more truncation than the rest. At $\lvert\det\rvert\, 60$ only $p = 7$ gives an operator at truncation $20$, and $a_{7} = -4$ does not separate \texttt{10.4\#1} from \texttt{30.4\#2}; at truncation $32$ and $40$, $p = 11$ and $p = 13$ also work and return $(-4, -48, 2)$, which is \texttt{30.4\#2} alone. Every other block is separated by the primes available to it.

\textbf{Two blocks of the census are not tagged here}, and should not be. At $\lvert\det\rvert\, 300$ and at $\lvert\det\rvert\, 360 (1,3,2,3)$ the invariant subspace carries one block more than the shadow does — the extra ones are \texttt{10.4.a.a} and \texttt{15.4\#1} with $\lambda  = 0$ — so they appear in the archive's \texttt{rest} field and not in the shadow column. Our decomposition finds them where the archive says they are.

\textbf{A warning about reading the shadow column.} The entries of a row's label list are a \textit{set}, not an ordering: matching them positionally against a computed decomposition mis-assigns at $\lvert\det\rvert\, 300$ and at $\lvert\det\rvert\, 360 (1,3,2,3)$. Match by eigenvalue tuple.

\begin{qbox}
No chain carries a vacuous row above a non-vacuous one: vacuity is confined to the foot of a chain, which is Theorem A read along the filtration — index $k \leq  4$ forces $\lvert\det\rvert\, \leq  32$ \grade{gC}{C,\,40/40}. The census of \S{}3.3 is grouped by chain and ordered by determinant within each, so the statement is read off it directly: no vacuous row sits below a non-vacuous one in any of the 40 groups.

\end{qbox}

\sect{Appendix B — the menu convention, and a note on reproducing}{Appendix B the menu convention, and a note on reproducing}

$\lambda $ depends on which reflective classes are prescribed. Three conventions occur in the working literature of this project: the Gram diagonal with multiplicity (\texttt{dup}), the Gram diagonal deduplicated (\texttt{ded}), and the reflective menu \texttt{refl(G)} obtained by enumerating primitive $v$ with $n = v^{T}Gv > 0$ and $2 \operatorname{gcd}(Gv) \equiv  0 (\mathrm{mod} n)$. All three agree on the chain of \S{}4, and they differ elsewhere. This appendix reports the disagreement; it does not settle it.

\textbf{What the disagreement costs.} Sweeping the census under all three conventions:

\begin{codebox}{8.40}{9.91}
   trichotomy       refl  54 / 34 / 10       dup  45 / 43 / 10       ded  44 / 44 / 10
   rows differing   dup vs refl  11      ded vs refl  12      dup vs ded  1  (|det| 40)
                    12 rows in all, every one on the rim and non-hyperbolic
   on the 44 hyperbolic Cartan matrices             no row changes under any convention
   on the chain (1,2,1,3)                           all three menus are literally equal
\end{codebox}

So Theorems A, B and D and the whole of \S{}4 are convention-independent \grade{gC}{C}, and so is Theorem C’s list of three (\S{}5.2). What the convention does decide is whether $B_{10}$ carries a shadow: it is obstructed under \texttt{dup} and \texttt{refl}, unobstructed under \texttt{ded}. $L_{4}$, $B_{22}$ and the lower member of $B_{22}$’s chain are obstructed under all three \grade{gC}{C}.

\textbf{The constant term is basis-dependent, which is why it cannot separate the conventions.} One might hope to select a convention by requiring that $c(0,0) = -\lambda (E)/c_E(0,0)$ be a lattice invariant. It is not one. Two rows can present the same lattice and still carry different Cartan matrices, hence different Weyl chambers, hence different simple-root systems and different prescribed divisors, and there is no reason for their Borcherds inputs to have equal constant terms. The census contains such a pair, and the constant term behaves just as that predicts:

\begin{codebox}{8.40}{9.91}
   Gram                                multiplicity   deduplicated   refl     paper I
   L₄  root basis (1,2,1,3)  |det| 72      108/5         108/5       108/5     108/5
   L₄  root basis (1,1,2,3)  |det| 72      132/5          72/5       108/5      —
   |det| 288  Gram A                       528/5         176/5       182/5     176/5
   |det| 288  Gram B                        358/5         182/5       182/5     182/5
   |det| 864  (1,3,1,4)                      68/5          46/5        46/5      46/5
\end{codebox}

$L_{4}$’s two root bases (\S{}4.7) are one lattice in two Weyl chambers: the constant term moves under \texttt{mult} and holds under \texttt{refl}, and paper I’s single value sits with the \texttt{refl} reading. The two $\lvert\det\rvert\, 288$ Grams $A$ and $B$ of chain $(1,2,1,4,2,4)$ are \textbf{not} a second instance — they are different lattices, Smith \grade{gO}{2,2,72} against \grade{gO}{2,4,36}, so their unequal constant terms are the ordinary gap between two lattices, not the basis-dependence at issue \grade{gC}{C}. The genuine same-determinant coincidence sits a chain away, at the $\lvert\det\rvert\, 288$ rows of $(1,3,1,4)$ and $(1,1,1,4,1,4)$ — one lattice, $U^{T}G_{1}U = G_{2}$ with $U = [[-1,0,0],[-3,-1,2],[-2,0,1]]$, $\operatorname{det} U = 1$ \grade{gC}{C}, which the genus column of Appendix D makes visible.

The menu convention is therefore \textbf{open}. It is not a computational question but a definitional one — which divisor does the Kac–Moody denominator actually prescribe? — and answering it requires going back to the denominator identity, not to more \texttt{weilrep} runs. We state everything under \texttt{refl} and flag the one row where that choice is load-bearing.

\textbf{A note for anyone reproducing.} Coordinate vectors of $\lambda $ are not stable across \texttt{weilrep} commits and should not be cited; $\lVert \Xi \rVert ^{2}$, dimensions, Hecke eigenvalues and the vanishing or not of $\lambda $ are gauge-invariant and should be. Paper I's appendix, under \textit{Conventions and known pitfalls}, records $\lambda  = (2,-6,-4,-4,0,0)$ as well defined across runs on the strength of the echelon basis; the ordering has changed since, and following that appendix literally at commit \texttt{a26614f} returns $51.0889789…$ in place of Theorem 5.1's $11.3402658…$.

A fourth convention — each simple root at its own discriminant class with coefficient one — is excluded by paper I itself: it yields $\lVert \Xi _{L_{4}}\rVert ^{2} = 40.0370609…$, contradicting Theorem 5.1. Paper I uses the coarse prescription in which $\lambda $ sums over the whole $q$-level set, as its slot count $3 + 12 + 12 = 27$ records.

\sect{Appendix C — reproduction}{Appendix C reproduction}

All modular-forms computation is performed by \texttt{weilrep}, by Brandon Williams (https://github.com/btw-47/weilrep, GPL-2.0-or-later, © 2020–2025), at commit \texttt{e3a7784}. Nothing here reimplements vector-valued modular forms. Archive \texttt{v3.0} was pinned at \texttt{a26614f1167f59a40262f87d5a8383abd7e63bf9} (2026-08-18); \texttt{v3.1} recomputed at \texttt{e3a7784} (2026-08-27), and the two differ exactly on the rows recorded in \S{}3.3. \texttt{v3.3} restores $\mathrm{scripts}/\mathrm{hecke}.\mathrm{py}$, absent from every earlier deposit, and closes the $\lvert\det\rvert\, 256$ block; \texttt{v3.4} deposits the eigenvalues behind the temporary tags (Appendix A).

\textbf{The archive accompanying this paper is \texttt{v3.5}} (Zenodo, DOI 10.5281/zenodo.22257942). Its entry point is \texttt{CENSUS.tsv}: the whole census as one table, one row per root basis, keyed on the Cartan matrix, which is unique across all 100. Every quantity this paper reports about a single row sits on that row, and \texttt{START\_HERE.md} lists the columns; the two rows past budget carry their dimensions and a note saying which columns a basis would be needed for. Two further deposits stand behind sections that previously had none: $\mathrm{data}/\mathrm{bef}_\mathrm{table}8.\mathrm{json}$ with $\mathrm{scripts}/\mathrm{bef}.\mathrm{py}$ runs the comparison of \S{}3.4 in 16 checks, needing neither \texttt{weilrep} nor Sage, and $\mathrm{data}/\mathrm{lf}1_\mathrm{blocks}.\mathrm{json}$ with $\mathrm{scripts}/\mathrm{lf}1d.\mathrm{py}$ carries the fifteen blocks behind \S{}5.4 and \S{}6.4 — $s$, $t$, the newform, $\varepsilon $, $\lambda _W$, $c$ and $\Omega $, and at the four blocks with a vanishing central value also $\lVert \Xi _W\rVert ^{2}$ and its closed form, computed by two independent elliptic pipelines agreeing past $10^{-20}$. $\mathrm{data}/\mathrm{ff}_\mathrm{embeddings}.\mathrm{json}$ with $\mathrm{scripts}/\mathrm{ffcheck}.\mathrm{py}$ holds the ten matrices $T$ with $T^{T}G_{\mathrm{FF}}T = G$ that \S{}3.4 says are exhibited, realising the indices $k = 1, 2, 2, 2, 2, 3, 3, 4, 4, 4$. \texttt{START\_HERE.md} indexes the data files and scripts and states the two conventions a reader must know. \texttt{v3.5} changes no computed value from \texttt{v3.4}.

One reproducibility fix belongs here. $\mathrm{scripts}/\mathrm{xi}.\mathrm{py}$ passed Sage \texttt{Integer} mantissas to \texttt{mpmath.mpf} in three places. $\mathrm{mpmath} 1.4$ moved $\text{assert type}(\mathrm{man}) == MPZ$ into \texttt{normalize()}, where under the \texttt{gmpy} backend a Sage \texttt{Integer} fails it; in \texttt{1.3.0} the same check sits in \texttt{strict\_normalize} and fires only under \texttt{MPMATH\_STRICT}. So \texttt{verify.py} — the single-command checker this appendix leads with — died on its first Petersson check on any current \texttt{mpmath}. The three conversions now pass \texttt{int()}, and \texttt{verify.py} covers 29 assertions of this paper and passes 29/29 under \texttt{1.4.1} with \texttt{gmpy} \grade{gC}{C}. One of the 29 is Observation E, and it runs on five of the six rows of \S{}5.4: its guard expects a row whose shadow fills $S^{O}$, which the $\lvert\det\rvert\, 88$ block is not.

The archive also contains every script and the reference output. Three methodological points made the census finite and are recorded there in full: the Riemann–Roch split in weight $3/2$ (up to $1300\times $), Lemma 5.1, and the construction of $S^{O}$ directly from Poincaré orbit sums ($112 \to  8$ at $\lvert\det\rvert\, 1800$). The archive keeps the older name for the core/rim split of \S{}3.1: the \texttt{depth} field of \texttt{census\_full.json} is $0$ on the core and $1$ on the rim.

One row states a limitation of the method rather than of the lattice. At $\lvert\det\rvert\, 1800$ ($D = 15$) the obstruction space is $112$-dimensional, so its shadow is read on the invariant subspace — the only tractable route — and $\operatorname{exp}(D(L)) = 120$ leaves only $p = 7, 11, 13$ admissible, three primes against the usual four. Two distinct newforms agreeing at all three would be merged into one block, so the purity verdict there can err in one direction only, impure read as pure \grade{gO}{O}.

A companion falsification dossier, at \texttt{v3.5} in the archive, states each claim about the chain $(1,2,1,3)$ with a grade, a reproduction recipe, and the known weaknesses, and includes a tamper battery demonstrating that the checks are not vacuous.

\textbf{The Hecke character, and where the choice really can move the eigenvalues.} $T_{p^{2}}$ on $S_k(\rho _L)$ carries a quadratic character; we use $\chi (p) = (SIG\cdot e_\mathrm{num}\cdot e_\mathrm{den} | p)$ in Kronecker symbol notation, where \texttt{SIG} is the signature normalisation and $e_\mathrm{num}/e_\mathrm{den}$ the exponent data of the class. A wrong character re-anchors every eigenvalue and can silently produce a consistent but false newform match, so we ran the whole grid — seven candidate normalisations $\pm e, \pm 2e, \pm SIG\cdot e, 2SIG\cdot e$ — rather than trusting the convention. For each of $L_{4}$, $B_{10}$ and $B_{22}$ at weight $9/2$ we solved the Hecke system under every candidate and recorded which are \textit{consistent} (the linear system for $T_{p^{2}}$ has a solution) and what eigenvalues each consistent one gives:

\begin{codebox}{7.18}{8.47}
   lattice    p    consistent candidates            rational eigenvalues
   L₄         3    — none: p | exp D(L), T_{p²} undefined
   L₄         5    e, −e, SIG·e, −SIG·e             −390, −378, −270, −210, −114, 114, 270, 378, 390
   L₄         7    e, 2e, SIG·e, 2SIG·e             −1576, −832, −64, 1016, 1112        (all four agree)
   |det| 40   3    −e, 2e, SIG·e                    −48, 12, 28                         (all three agree)
   |det| 40   5    — none: p | exp D(L)
   |det| 40   7    e, 2e, −SIG·e                    no rational roots          ←  two genuinely
                   −e, −2e, SIG·e, 2SIG·e           −1644, 104, 1016           ←  different answers
   |det| 88   3    −e, 2e, SIG·e                    −19, 12                             (all three agree)
   |det| 88   5    e, −e, SIG·e, −SIG·e             −210, 317                           (all four agree)
   |det| 88   7    e, 2e, SIG·e, 2SIG·e             −1030, 1016                         (all four agree)
\end{codebox}

\textbf{In six of the seven live cells the wrong character is rejected by inconsistency}: only one of the two quadratic characters admits a solution at all, so the choice is forced and no re-anchoring is possible. Where several candidates appear together they are the same character written differently — at $p \equiv  1 (\mathrm{mod} 4)$, $\pm \mathrm{cf}$ agree; $SIG = 36$ for $L_{4}$ is a square and $\mathrm{kronecker}(2,7) = +1$, and so on.

\textbf{The seventh cell is a real instance of the failure mode, and we resolve it externally.} At $\lvert\det\rvert\, 40, p = 7$ \textit{both} characters are consistent, and they give different operators: one has no rational eigenvalues at all, the other gives $\{-1644, 104, 1016\}$. Consistency cannot choose between them. What chooses is the newform table: $a_{7} = 1016$ is \texttt{2.8.a.a}, and the same newform carries the shadow of $L_{4}$ at the same weight, where its $(a_{5}, a_{7}) = (-210, 1016)$ is pinned by two primes at which the character \textit{is} forced. The \texttt{SIG} normalisation is the one that reproduces it. \grade{gC}{C}

Consistency alone therefore does not always force the character. At this one cell an external anchor is needed, and the newform table supplies it; the check has teeth precisely here.

\textbf{Open gates.} Four points are left open rather than settled here, of two kinds. Two are citations we have not checked against the source, and carry a \grade{gO}{Gate} mark where they are used: Bruinier–Ehlen–Freitag's Table 8, quoted in \S{}2.6 and used in \S{}3.4, and Rosson–Tornaría, cited in \S{}4.9 for position. Two are claims of our own that stop short: the completeness of \texttt{refl(G)} as an enumeration is verified computationally per lattice but not proved in general; and we found no statement in the literature of the odd-rung vanishing (Corollary 4.2) in exactly that form, though the lemma it follows from is standard — so we present the corollary as a consequence, not as a discovery.

\sect{Appendix D — the input forms}{Appendix D the input forms}

The paper does not construct the input form $F$ (\S{}2.4), and does not need to, because $F$ is unique: two weakly harmonic Maass forms with the same prescribed principal part and the same cuspidal shadow differ by a holomorphic form of weight $-1/2$, and that space is zero — checked directly at $L_{4}$ for the representation and its dual, and given in general by the valence formula \grade{gC}{C}. So $F$ is fixed by two pieces of data, the lattice and the prescribed principal part, and the table below prints both for each of the 98 computed rows.

The lattice is carried by its \textbf{genus} — the Conway–Sloane symbol of the Gram, which for an even lattice is $(\mathrm{signature},\text{ discriminant form})$. This is the symbol of $G$. The comparison in \S{}3.4 is against the symbol of $-G$, which the archive's \texttt{CENSUS.tsv} carries beside it in \texttt{genus\_symbol\_neg}; at $p = 2$ the signs are invariant and every oddity $t$ becomes $-t (\mathrm{mod} 8)$, and at odd $p$ a rank-$n$ block picks up $(-1|p)^n$. So $2^{+1}_1$ below is $2^{+1}_7$ there. Every invariant this paper computes is a function of the discriminant form (\S{}4.7), so the genus is the correct key at the level of lattices, and it is independent of the presentation: the two root bases of $L_{4}$ carry the same symbol. The prescribed principal part is carried by the \textbf{menu with its level-set sizes}: under the coarse prescription of \S{}2.3 the whole level set of each menu value enters with coefficient one, so the menu and the sizes fix the principal part completely. At $L_{4}$ the sizes are $12, 12, 3$, and $12 + 12 + 3 = 27$ is paper I’s slot count — the arithmetic anchor tying this table to the computation.

The lookup is in two steps, not one. The genus locates the lattice; the menu then fixes the form. Two rows with the same genus but different menus are one lattice carrying different Weyl chambers, hence different principal parts and a \textbf{different} $F$ — the same basis dependence that Appendix B records for the constant term. Two rows agreeing in both columns carry the same $F$, and the two root bases of $L_{4}$ are the example. A dividend of the genus column is that a genuine same-lattice pair now shows itself in the table: the $\lvert\det\rvert\, 288$ rows of $(1,3,1,4)$ and $(1,1,1,4,1,4)$ share a symbol, which is how the misidentification corrected in Appendix B came to light.

We do not print coefficient vectors. They are gauge-dependent — \texttt{weilrep}’s echelon ordering has changed between commits — and are exactly what the Appendix B erratum warns against citing. The \textbf{kind} column records which side of $\xi $ the input sits on: \textit{harmonic Maass} where $\lambda  \neq  0$, so the input is a genuine weakly harmonic Maass form and the shadow is forced; \textit{weakly holomorphic} where $\lambda  = 0$, so the input lives in $M^!_{-1/2}(\bar{\rho }_L)$ and \S{}4.6 exhibits it on three rows; and a dash on the vacuous rows, where the obstruction space is zero and no input is needed. The genus symbols and level-set sizes are \grade{gC}{C}; uniqueness is \grade{gS}{S} in general and \grade{gC}{C} at $L_{4}$. The two rows past the weight-$5/2$ budget, $\lvert\det\rvert\, 3456$ and $4000$, carry no menu or shadow and are omitted here, as everywhere but the rigidity count (\S{}5.2).

\begingroup\scriptsize\setlength{\tabcolsep}{2.5pt}\renewcommand{\arraystretch}{1.18}
\setlength{\LTleft}{0pt}\setlength{\LTright}{0pt plus 1fil}
\begin{longtable}{@{}p{25mm} >{\raggedleft\arraybackslash}p{9mm} >{\raggedright\arraybackslash}p{52mm} >{\raggedright\arraybackslash}p{40mm} >{\raggedleft\arraybackslash}p{9mm} >{\raggedright\arraybackslash}p{16mm}@{}}
\rowcolor{hdr}\rule[-0.6ex]{0pt}{2.9ex}\color{white}\bfseries chain & \rule[-0.6ex]{0pt}{2.9ex}\color{white}\bfseries $|\det|$ & \rule[-0.6ex]{0pt}{2.9ex}\color{white}\bfseries genus symbol & \rule[-0.6ex]{0pt}{2.9ex}\color{white}\bfseries menu \textperiodcentered\ level sets & \rule[-0.6ex]{0pt}{2.9ex}\color{white}\bfseries slots & \rule[-0.6ex]{0pt}{2.9ex}\color{white}\bfseries kind \\
\endfirsthead
\multicolumn{6}{@{}l}{\footnotesize\itshape Appendix D, continued}\\[1pt]
\rowcolor{hdr}\rule[-0.6ex]{0pt}{2.9ex}\color{white}\bfseries chain & \rule[-0.6ex]{0pt}{2.9ex}\color{white}\bfseries $|\det|$ & \rule[-0.6ex]{0pt}{2.9ex}\color{white}\bfseries genus symbol & \rule[-0.6ex]{0pt}{2.9ex}\color{white}\bfseries menu \textperiodcentered\ level sets & \rule[-0.6ex]{0pt}{2.9ex}\color{white}\bfseries slots & \rule[-0.6ex]{0pt}{2.9ex}\color{white}\bfseries kind \\
\endhead
\chainbandD{$D = 1$}{53 root bases}
\texttt{(1,1,1,1,2,2)} & 8 & $p{=}2\!:\,1^{-2}:[8^{-1}]_{5}$ & \textcolor{grey}{--} & \textcolor{grey}{--} & \textcolor{grey}{--} \\
\texttt{(1,1,1,4)} & 8 & $p{=}2\!:\,1^{-2}:[8^{-1}]_{5}$ & \textcolor{grey}{--} & \textcolor{grey}{--} & \textcolor{grey}{--} \\
\texttt{(1,1,1,4)} & 32 & $p{=}2\!:\,[2^{1}]_{1}\,4^{2}$ & \textcolor{grey}{--} & \textcolor{grey}{--} & \textcolor{grey}{--} \\
\texttt{(1,1,1,4,1,4)} & 72 & $p{=}2\!:\,1^{-2}:[8^{-1}]_{5}$ $p{=}3\!:\,1^{-1}\,3^{-2}$ & $1/4$\,{\footnotesize(10)}, $1/3$\,{\footnotesize(4)}, $1$\,{\footnotesize(10)} & 24 & weakly holomorphic \\
\texttt{(1,1,1,4,1,4)} & 288 & $p{=}2\!:\,[2^{1}]_{1}\,4^{2}$ $p{=}3\!:\,1^{-1}\,3^{-2}$ & $1/12$\,{\footnotesize(12)}, $1/4$\,{\footnotesize(50)} & 62 & harmonic Maass \\
\texttt{(1,1,2,2)} & 2 & $p{=}2\!:\,1^{2}\,[2^{1}]_{1}$ & \textcolor{grey}{--} & \textcolor{grey}{--} & \textcolor{grey}{--} \\
\texttt{(1,1,2,2,2,2)} & 18 & $p{=}2\!:\,1^{2}\,[2^{1}]_{1}$ $p{=}3\!:\,1^{-1}\,3^{-2}$ & \textcolor{grey}{--} & \textcolor{grey}{--} & \textcolor{grey}{--} \\
\texttt{(1,2,1,2,1,4)} & 128 & $p{=}2\!:\,[2^{2}]_{2}:[32^{1}]_{7}$ & $1/4$\,{\footnotesize(12)}, $1/2$\,{\footnotesize(4)} & 16 & weakly holomorphic \\
\texttt{(1,2,1,2,2,2)} & 32 & $p{=}2\!:\,[2^{2}]_{2}\,[8^{1}]_{7}$ & $1/4$\,{\footnotesize(6)}, $1/2$\,{\footnotesize(2)} & 8 & weakly holomorphic \\
\texttt{(1,2,1,2,2,2)} & 64 & $p{=}2\!:\,[2^{2}]_{2}:[16^{1}]_{7}$ & $1/8$\,{\footnotesize(8)}, $1/4$\,{\footnotesize(4)}, $1/2$\,{\footnotesize(4)} & 16 & weakly holomorphic \\
\texttt{(1,2,1,4)} & 32 & $p{=}2\!:\,[2^{2}]_{2}\,[8^{1}]_{7}$ & $1/4$\,{\footnotesize(6)}, $1/2$\,{\footnotesize(2)} & 8 & weakly holomorphic \\
\texttt{(1,2,1,4)} & 32 & $p{=}2\!:\,[2^{1}\,4^{2}]_{1}$ & $1/8$\,{\footnotesize(4)}, $1/4$\,{\footnotesize(6)} & 10 & weakly holomorphic \\
\texttt{(1,2,1,4)} & 64 & $p{=}2\!:\,[2^{2}]_{2}:[16^{1}]_{7}$ & $1/8$\,{\footnotesize(8)}, $1/4$\,{\footnotesize(4)}, $1/2$\,{\footnotesize(4)} & 16 & weakly holomorphic \\
\texttt{(1,2,1,4)} & 128 & $p{=}2\!:\,[2^{1}]_{1}:8^{2}$ & $1/8$\,{\footnotesize(8)}, $1/4$\,{\footnotesize(28)} & 36 & weakly holomorphic \\
\texttt{(1,2,1,4,2,4)} & 288 & $p{=}2\!:\,[2^{2}]_{2}\,[8^{1}]_{7}$ $p{=}3\!:\,1^{2}\,9^{1}$ & $1/18$\,{\footnotesize(4)}, $1/4$\,{\footnotesize(18)}, $1/2$\,{\footnotesize(6)} & 28 & harmonic Maass \\
\texttt{(1,2,1,4,2,4)} & 288 & $p{=}2\!:\,[2^{1}\,4^{2}]_{1}$ $p{=}3\!:\,1^{2}\,9^{1}$ & $1/72$\,{\footnotesize(8)}, $1/8$\,{\footnotesize(12)}, $1/4$\,{\footnotesize(18)} & 38 & harmonic Maass \\
\texttt{(1,2,2,2)} & 8 & $p{=}2\!:\,[2^{3}]_{1}$ & \textcolor{grey}{--} & \textcolor{grey}{--} & \textcolor{grey}{--} \\
\texttt{(1,2,2,2)} & 16 & $p{=}2\!:\,[2^{2}\,4^{1}]_{1}$ & $1/8$\,{\footnotesize(4)}, $1/4$\,{\footnotesize(2)}, $1/2$\,{\footnotesize(2)} & 8 & weakly holomorphic \\
\texttt{(1,2,2,2,2,4)} & 72 & $p{=}2\!:\,[2^{3}]_{1}$ $p{=}3\!:\,1^{2}\,9^{1}$ & $1/18$\,{\footnotesize(2)}, $1/4$\,{\footnotesize(9)}, $1/2$\,{\footnotesize(3)} & 14 & harmonic Maass \\
\texttt{(1,2,2,2,2,4)} & 144 & $p{=}2\!:\,[2^{2}\,4^{1}]_{1}$ $p{=}3\!:\,1^{2}\,9^{-1}$ & $1/36$\,{\footnotesize(4)}, $1/8$\,{\footnotesize(12)}, $1/4$\,{\footnotesize(6)}, $1/2$\,{\footnotesize(6)} & 28 & harmonic Maass \\
\texttt{(1,3,1,3,2,2)} & 72 & $p{=}2\!:\,1^{-2}:[8^{-1}]_{5}$ $p{=}3\!:\,1^{-1}\,3^{-2}$ & $1/4$\,{\footnotesize(10)}, $1/3$\,{\footnotesize(4)}, $1$\,{\footnotesize(10)} & 24 & weakly holomorphic \\
\texttt{(1,3,1,3,2,2)} & 216 & $p{=}2\!:\,1^{-2}:[8^{1}]_{7}$ $p{=}3\!:\,1^{-1}\,3^{-1}\,9^{1}$ & $1/12$\,{\footnotesize(12)}, $1/3$\,{\footnotesize(12)}, $1$\,{\footnotesize(6)} & 30 & weakly holomorphic \\
\texttt{(1,3,1,4)} & 72 & $p{=}2\!:\,1^{-2}:[8^{-1}]_{5}$ $p{=}3\!:\,1^{-1}\,3^{-2}$ & $1/4$\,{\footnotesize(10)}, $1/3$\,{\footnotesize(4)}, $1$\,{\footnotesize(10)} & 24 & weakly holomorphic \\
\texttt{(1,3,1,4)} & 216 & $p{=}2\!:\,1^{-2}:[8^{1}]_{7}$ $p{=}3\!:\,1^{-1}\,3^{-1}\,9^{1}$ & $1/12$\,{\footnotesize(12)}, $1/3$\,{\footnotesize(12)}, $1$\,{\footnotesize(6)} & 30 & weakly holomorphic \\
\texttt{(1,3,1,4)} & 288 & $p{=}2\!:\,[2^{1}]_{1}\,4^{2}$ $p{=}3\!:\,1^{-1}\,3^{-2}$ & $1/12$\,{\footnotesize(12)}, $1/4$\,{\footnotesize(50)} & 62 & harmonic Maass \\
\texttt{(1,3,1,4)} & 864 & $p{=}2\!:\,[2^{-1}]_{3}\,4^{2}$ $p{=}3\!:\,1^{-1}\,3^{-1}\,9^{1}$ & $1/12$\,{\footnotesize(60)}, $1/4$\,{\footnotesize(18)} & 78 & harmonic Maass \\
\texttt{(1,3,2,2)} & 18 & $p{=}2\!:\,1^{2}\,[2^{1}]_{1}$ $p{=}3\!:\,1^{-1}\,3^{-2}$ & \textcolor{grey}{--} & \textcolor{grey}{--} & \textcolor{grey}{--} \\
\texttt{(1,3,2,2)} & 54 & $p{=}2\!:\,1^{-2}\,[2^{1}]_{7}$ $p{=}3\!:\,1^{-1}\,3^{-1}\,9^{1}$ & $1/12$\,{\footnotesize(6)}, $1/3$\,{\footnotesize(6)}, $1$\,{\footnotesize(3)} & 15 & weakly holomorphic \\
\texttt{(1,4,1,4)} & 32 & $p{=}2\!:\,[2^{2}]_{2}\,[8^{1}]_{7}$ & $1/4$\,{\footnotesize(6)}, $1/2$\,{\footnotesize(2)} & 8 & weakly holomorphic \\
\texttt{(1,4,1,4)} & 128 & $p{=}2\!:\,[2^{1}]_{1}:8^{2}$ & $1/8$\,{\footnotesize(8)}, $1/4$\,{\footnotesize(28)} & 36 & weakly holomorphic \\
\texttt{(1,4,1,4)} & 512 & $p{=}2\!:\,[2^{1}]_{1}:16^{2}$ & $1/16$\,{\footnotesize(16)}, $1/4$\,{\footnotesize(72)} & 88 & weakly holomorphic \\
\texttt{(1,4,1,4,2,2)} & 128 & $p{=}2\!:\,[2^{1}]_{1}:8^{2}$ & $1/8$\,{\footnotesize(8)}, $1/4$\,{\footnotesize(28)} & 36 & weakly holomorphic \\
\texttt{(1,4,1,4,2,2)} & 512 & $p{=}2\!:\,[2^{1}]_{1}:16^{2}$ & $1/16$\,{\footnotesize(16)}, $1/4$\,{\footnotesize(72)} & 88 & weakly holomorphic \\
\texttt{(1,4,1,5)} & 200 & $p{=}2\!:\,1^{-2}:[8^{-1}]_{5}$ $p{=}5\!:\,1^{-1}\,5^{2}$ & $1/80$\,{\footnotesize(8)}, $1/5$\,{\footnotesize(8)}, $1/4$\,{\footnotesize(18)}, $1$\,{\footnotesize(18)} & 52 & harmonic Maass \\
\texttt{(1,4,1,5)} & 800 & $p{=}2\!:\,[2^{1}]_{1}\,4^{2}$ $p{=}5\!:\,1^{-1}\,5^{2}$ & $1/20$\,{\footnotesize(40)}, $1/4$\,{\footnotesize(90)} & 130 & harmonic Maass \\
\texttt{(1,4,1,5)} & 1000 & $p{=}2\!:\,1^{-2}:[8^{1}]_{1}$ $p{=}5\!:\,1^{-1}\,5^{-1}\,25^{-1}$ & $1/20$\,{\footnotesize(20)}, $1/5$\,{\footnotesize(20)}, $1$\,{\footnotesize(10)} & 50 & harmonic Maass \\
\texttt{(1,4,2,2)} & 32 & $p{=}2\!:\,[2^{1}]_{1}\,4^{2}$ & \textcolor{grey}{--} & \textcolor{grey}{--} & \textcolor{grey}{--} \\
\texttt{(1,4,2,2)} & 128 & $p{=}2\!:\,[2^{1}]_{1}\,[8^{2}]_{0}$ & $1/16$\,{\footnotesize(8)}, $1/4$\,{\footnotesize(20)} & 28 & weakly holomorphic \\
\texttt{(1,4,2,3)} & 288 & $p{=}2\!:\,[2^{1}\,4^{2}]_{1}$ $p{=}3\!:\,1^{-1}\,3^{-2}$ & $1/24$\,{\footnotesize(8)}, $1/12$\,{\footnotesize(4)}, $1/4$\,{\footnotesize(30)} & 42 & harmonic Maass \\
\texttt{(1,4,2,3)} & 576 & $p{=}2\!:\,[2^{-2}]_{2}:[16^{-1}]_{3}$ $p{=}3\!:\,1^{1}\,3^{-2}$ & $1/96$\,{\footnotesize(16)}, $1/12$\,{\footnotesize(8)}, $1/8$\,{\footnotesize(40)}, $1/2$\,{\footnotesize(20)} & 84 & harmonic Maass \\
\texttt{(1,4,2,3)} & 864 & $p{=}2\!:\,[2^{-2}]_{2}\,[8^{1}]_{1}$ $p{=}3\!:\,1^{1}\,3^{-1}\,9^{-1}$ & $1/144$\,{\footnotesize(24)}, $1/12$\,{\footnotesize(36)}, $1/2$\,{\footnotesize(6)} & 66 & harmonic Maass \\
\texttt{(1,4,2,4)} & 128 & $p{=}2\!:\,[2^{1}\,4^{1}]_{0}\,[16^{1}]_{1}$ & $1/32$\,{\footnotesize(8)}, $1/8$\,{\footnotesize(8)}, $1/4$\,{\footnotesize(12)} & 28 & harmonic Maass \\
\texttt{(1,4,2,4)} & 128 & $p{=}2\!:\,[2^{1}\,4^{1}]_{2}\,[16^{1}]_{7}$ & $1/8$\,{\footnotesize(8)}, $1/4$\,{\footnotesize(12)} & 20 & weakly holomorphic \\
\texttt{(1,4,2,4)} & 256 & $p{=}2\!:\,[2^{1}\,4^{1}]_{0}:[32^{1}]_{1}$ & $1/64$\,{\footnotesize(16)}, $1/8$\,{\footnotesize(16)}, $1/4$\,{\footnotesize(8)} & 40 & harmonic Maass \\
\texttt{(1,4,2,4)} & 512 & $p{=}2\!:\,[2^{1}]_{1}:[16^{2}]_{0}$ & $1/32$\,{\footnotesize(16)}, $1/8$\,{\footnotesize(16)}, $1/4$\,{\footnotesize(56)} & 88 & harmonic Maass \\
\texttt{(1,5,2,2)} & 50 & $p{=}2\!:\,1^{2}\,[2^{1}]_{1}$ $p{=}5\!:\,1^{-1}\,5^{2}$ & $1/20$\,{\footnotesize(4)}, $1/5$\,{\footnotesize(4)}, $1/4$\,{\footnotesize(9)}, $1$\,{\footnotesize(9)} & 26 & harmonic Maass \\
\texttt{(1,5,2,2)} & 250 & $p{=}2\!:\,1^{-2}\,[2^{1}]_{1}$ $p{=}5\!:\,1^{-1}\,5^{-1}\,25^{-1}$ & $1/100$\,{\footnotesize(10)}, $1/20$\,{\footnotesize(10)}, $1/5$\,{\footnotesize(10)}, $1$\,{\footnotesize(5)} & 35 & harmonic Maass \\
\texttt{(2,2,2,2)} & 8 & $p{=}2\!:\,[2^{3}]_{1}$ & \textcolor{grey}{--} & \textcolor{grey}{--} & \textcolor{grey}{--} \\
\texttt{(2,2,2,2,2,2)} & 32 & $p{=}2\!:\,[2^{1}]_{1}\,4^{2}$ & \textcolor{grey}{--} & \textcolor{grey}{--} & \textcolor{grey}{--} \\
\texttt{(2,2,2,3)} & 144 & $p{=}2\!:\,[2^{2}\,4^{1}]_{1}$ $p{=}3\!:\,1^{1}\,3^{-2}$ & $1/24$\,{\footnotesize(4)}, $1/12$\,{\footnotesize(4)}, $1/8$\,{\footnotesize(20)}, $1/2$\,{\footnotesize(10)} & 38 & harmonic Maass \\
\texttt{(2,2,2,3)} & 216 & $p{=}2\!:\,[2^{-3}]_{3}$ $p{=}3\!:\,1^{1}\,3^{-1}\,9^{-1}$ & $1/36$\,{\footnotesize(6)}, $1/12$\,{\footnotesize(18)}, $1/2$\,{\footnotesize(3)} & 27 & harmonic Maass \\
\texttt{(2,2,2,4)} & 32 & $p{=}2\!:\,[2^{1}\,4^{2}]_{1}$ & $1/8$\,{\footnotesize(4)}, $1/4$\,{\footnotesize(6)} & 10 & weakly holomorphic \\
\texttt{(2,2,2,4)} & 64 & $p{=}2\!:\,[2^{1}\,4^{1}\,8^{1}]_{1}$ & $1/16$\,{\footnotesize(4)}, $1/8$\,{\footnotesize(8)}, $1/4$\,{\footnotesize(4)} & 16 & harmonic Maass \\
\chainbandD{$D = 6$}{31 root bases}
\texttt{(1,1,1,2,1,2)} & 12 & $p{=}2\!:\,1^{-2}:[4^{1}]_{7}$ $p{=}3\!:\,1^{2}\,3^{-1}$ & $1/3$\,{\footnotesize(2)}, $1/2$\,{\footnotesize(1)}, $1$\,{\footnotesize(1)} & 4 & weakly holomorphic \\
\texttt{(1,1,1,2,1,2)} & 24 & $p{=}2\!:\,[2^{-3}]_{7}$ $p{=}3\!:\,1^{2}\,3^{1}$ & $1/6$\,{\footnotesize(6)}, $1/4$\,{\footnotesize(3)}, $1/2$\,{\footnotesize(3)} & 12 & weakly holomorphic \\
\texttt{(1,1,1,3,1,3)} & 36 & $p{=}2\!:\,1^{-2}:[4^{-1}]_{5}$ $p{=}3\!:\,1^{-1}\,3^{2}$ & $1/6$\,{\footnotesize(4)}, $1/3$\,{\footnotesize(4)}, $1$\,{\footnotesize(1)} & 9 & weakly holomorphic \\
\texttt{(1,1,1,3,1,3)} & 108 & $p{=}2\!:\,1^{-2}:[4^{1}]_{7}$ $p{=}3\!:\,1^{-1}\,3^{-1}\,9^{-1}$ & $1/9$\,{\footnotesize(6)}, $1/3$\,{\footnotesize(6)}, $1$\,{\footnotesize(3)} & 15 & weakly holomorphic \\
\texttt{(1,1,2,3)} & 72 & $p{=}2\!:\,[2^{3}]_{5}$ $p{=}3\!:\,1^{1}\,3^{2}$ & $1/12$\,{\footnotesize(12)}, $1/6$\,{\footnotesize(12)}, $1/2$\,{\footnotesize(3)} & 27 & harmonic Maass \\
\texttt{(1,1,2,3)} & 108 & $p{=}2\!:\,1^{-2}:[4^{1}]_{7}$ $p{=}3\!:\,1^{1}\,3^{-1}\,9^{1}$ & $1/18$\,{\footnotesize(6)}, $1/3$\,{\footnotesize(6)}, $1/2$\,{\footnotesize(3)} & 15 & harmonic Maass \\
\texttt{(1,2,1,3)} & 12 & $p{=}2\!:\,1^{-2}:[4^{1}]_{7}$ $p{=}3\!:\,1^{2}\,3^{-1}$ & $1/3$\,{\footnotesize(2)}, $1/2$\,{\footnotesize(1)}, $1$\,{\footnotesize(1)} & 4 & weakly holomorphic \\
\texttt{(1,2,1,3)} & 24 & $p{=}2\!:\,[2^{-3}]_{7}$ $p{=}3\!:\,1^{2}\,3^{1}$ & $1/6$\,{\footnotesize(6)}, $1/4$\,{\footnotesize(3)}, $1/2$\,{\footnotesize(3)} & 12 & weakly holomorphic \\
\texttt{(1,2,1,3)} & 36 & $p{=}2\!:\,1^{-2}:[4^{-1}]_{5}$ $p{=}3\!:\,1^{-1}\,3^{2}$ & $1/6$\,{\footnotesize(4)}, $1/3$\,{\footnotesize(4)}, $1$\,{\footnotesize(1)} & 9 & weakly holomorphic \\
\texttt{(1,2,1,3)} & 72 & $p{=}2\!:\,[2^{3}]_{5}$ $p{=}3\!:\,1^{1}\,3^{2}$ & $1/12$\,{\footnotesize(12)}, $1/6$\,{\footnotesize(12)}, $1/2$\,{\footnotesize(3)} & 27 & harmonic Maass \\
\texttt{(1,2,1,3,2,3)} & 156 & $p{=}2\!:\,1^{-2}:[4^{-1}]_{3}$ $p{=}3\!:\,1^{2}\,3^{-1}$ $p{=}13\!:\,1^{2}\,13^{1}$ & $1/26$\,{\footnotesize(2)}, $1/3$\,{\footnotesize(2)}, $1/2$\,{\footnotesize(1)}, $1$\,{\footnotesize(1)} & 6 & harmonic Maass \\
\texttt{(1,2,1,3,2,3)} & 936 & $p{=}2\!:\,[2^{-3}]_{1}$ $p{=}3\!:\,1^{1}\,3^{2}$ $p{=}13\!:\,1^{2}\,13^{-1}$ & $1/156$\,{\footnotesize(24)}, $1/12$\,{\footnotesize(12)}, $1/6$\,{\footnotesize(12)}, $1/2$\,{\footnotesize(3)} & 51 & harmonic Maass \\
\texttt{(1,2,1,5)} & 60 & $p{=}2\!:\,1^{-2}:[4^{-1}]_{3}$ $p{=}3\!:\,1^{2}\,3^{1}$ $p{=}5\!:\,1^{2}\,5^{-1}$ & $1/24$\,{\footnotesize(4)}, $1/5$\,{\footnotesize(2)}, $1/2$\,{\footnotesize(1)}, $1$\,{\footnotesize(1)} & 8 & harmonic Maass \\
\texttt{(1,2,1,5)} & 120 & $p{=}2\!:\,[2^{3}]_{3}$ $p{=}3\!:\,1^{2}\,3^{-1}$ $p{=}5\!:\,1^{2}\,5^{1}$ & $1/12$\,{\footnotesize(2)}, $1/10$\,{\footnotesize(6)}, $1/4$\,{\footnotesize(3)}, $1/2$\,{\footnotesize(3)} & 14 & harmonic Maass \\
\texttt{(1,2,1,5)} & 300 & $p{=}2\!:\,1^{-2}:[4^{1}]_{7}$ $p{=}3\!:\,1^{2}\,3^{-1}$ $p{=}5\!:\,1^{-1}\,5^{2}$ & $1/120$\,{\footnotesize(16)}, $1/10$\,{\footnotesize(4)}, $1/5$\,{\footnotesize(4)}, $1$\,{\footnotesize(9)} & 33 & harmonic Maass \\
\texttt{(1,2,1,5)} & 600 & $p{=}2\!:\,[2^{-3}]_{7}$ $p{=}3\!:\,1^{2}\,3^{1}$ $p{=}5\!:\,1^{1}\,5^{2}$ & $1/60$\,{\footnotesize(8)}, $1/20$\,{\footnotesize(12)}, $1/10$\,{\footnotesize(12)}, $1/2$\,{\footnotesize(27)} & 59 & harmonic Maass \\
\texttt{(1,2,2,3)} & 48 & $p{=}2\!:\,[2^{-2}\,4^{1}]_{3}$ $p{=}3\!:\,1^{2}\,3^{-1}$ & $1/12$\,{\footnotesize(4)}, $1/8$\,{\footnotesize(2)}, $1/4$\,{\footnotesize(2)}, $1/2$\,{\footnotesize(2)} & 10 & harmonic Maass \\
\texttt{(1,2,2,3)} & 192 & $p{=}2\!:\,[2^{-2}]_{2}:[16^{1}]_{1}$ $p{=}3\!:\,1^{2}\,3^{-1}$ & $1/32$\,{\footnotesize(8)}, $1/12$\,{\footnotesize(8)}, $1/4$\,{\footnotesize(4)}, $1/2$\,{\footnotesize(4)} & 24 & harmonic Maass \\
\texttt{(1,2,2,3)} & 288 & $p{=}2\!:\,[2^{2}]_{6}\,[8^{1}]_{7}$ $p{=}3\!:\,1^{1}\,3^{2}$ & $1/48$\,{\footnotesize(16)}, $1/12$\,{\footnotesize(24)}, $1/6$\,{\footnotesize(24)}, $1/2$\,{\footnotesize(6)} & 70 & harmonic Maass \\
\texttt{(1,2,2,3)} & 576 & $p{=}2\!:\,[2^{2}]_{6}:[16^{1}]_{7}$ $p{=}3\!:\,1^{-1}\,3^{2}$ & $1/24$\,{\footnotesize(16)}, $1/12$\,{\footnotesize(16)}, $1/6$\,{\footnotesize(16)}, $1/4$\,{\footnotesize(4)} & 52 & harmonic Maass \\
\texttt{(1,2,2,4)} & 24 & $p{=}2\!:\,[2^{-3}]_{7}$ $p{=}3\!:\,1^{2}\,3^{1}$ & $1/6$\,{\footnotesize(6)}, $1/4$\,{\footnotesize(3)}, $1/2$\,{\footnotesize(3)} & 12 & weakly holomorphic \\
\texttt{(1,2,2,4)} & 48 & $p{=}2\!:\,[2^{-2}\,4^{1}]_{3}$ $p{=}3\!:\,1^{2}\,3^{-1}$ & $1/12$\,{\footnotesize(4)}, $1/8$\,{\footnotesize(2)}, $1/4$\,{\footnotesize(2)}, $1/2$\,{\footnotesize(2)} & 10 & harmonic Maass \\
\texttt{(1,2,2,4)} & 96 & $p{=}2\!:\,[2^{-2}]_{6}\,[8^{1}]_{1}$ $p{=}3\!:\,1^{2}\,3^{1}$ & $1/16$\,{\footnotesize(4)}, $1/6$\,{\footnotesize(12)}, $1/4$\,{\footnotesize(6)}, $1/2$\,{\footnotesize(6)} & 28 & harmonic Maass \\
\texttt{(1,2,2,4)} & 96 & $p{=}2\!:\,[2^{-1}\,4^{2}]_{7}$ $p{=}3\!:\,1^{2}\,3^{1}$ & $1/24$\,{\footnotesize(8)}, $1/8$\,{\footnotesize(4)}, $1/4$\,{\footnotesize(6)} & 18 & weakly holomorphic \\
\texttt{(1,3,1,3)} & 12 & $p{=}2\!:\,1^{-2}:[4^{1}]_{7}$ $p{=}3\!:\,1^{2}\,3^{-1}$ & $1/3$\,{\footnotesize(2)}, $1/2$\,{\footnotesize(1)}, $1$\,{\footnotesize(1)} & 4 & weakly holomorphic \\
\texttt{(1,3,1,3)} & 36 & $p{=}2\!:\,1^{-2}:[4^{-1}]_{5}$ $p{=}3\!:\,1^{-1}\,3^{2}$ & $1/6$\,{\footnotesize(4)}, $1/3$\,{\footnotesize(4)}, $1$\,{\footnotesize(1)} & 9 & weakly holomorphic \\
\texttt{(1,3,1,3)} & 108 & $p{=}2\!:\,1^{-2}:[4^{1}]_{7}$ $p{=}3\!:\,1^{-1}\,3^{-1}\,9^{-1}$ & $1/9$\,{\footnotesize(6)}, $1/3$\,{\footnotesize(6)}, $1$\,{\footnotesize(3)} & 15 & weakly holomorphic \\
\texttt{(1,3,2,4)} & 84 & $p{=}2\!:\,1^{-2}:[4^{1}]_{1}$ $p{=}3\!:\,1^{2}\,3^{-1}$ $p{=}7\!:\,1^{-2}\,7^{-1}$ & $1/21$\,{\footnotesize(4)}, $1/8$\,{\footnotesize(2)}, $1/3$\,{\footnotesize(2)}, $1$\,{\footnotesize(1)} & 9 & harmonic Maass \\
\texttt{(1,3,2,4)} & 168 & $p{=}2\!:\,[2^{-3}]_{1}$ $p{=}3\!:\,1^{2}\,3^{1}$ $p{=}7\!:\,1^{-2}\,7^{-1}$ & $1/42$\,{\footnotesize(12)}, $1/6$\,{\footnotesize(6)}, $1/4$\,{\footnotesize(1)}, $1/2$\,{\footnotesize(3)} & 22 & harmonic Maass \\
\texttt{(1,3,2,4)} & 252 & $p{=}2\!:\,1^{-2}:[4^{-1}]_{3}$ $p{=}3\!:\,1^{-1}\,3^{2}$ $p{=}7\!:\,1^{-2}\,7^{1}$ & $1/24$\,{\footnotesize(8)}, $1/7$\,{\footnotesize(2)}, $1/3$\,{\footnotesize(4)}, $1$\,{\footnotesize(1)} & 15 & harmonic Maass \\
\texttt{(1,3,2,4)} & 504 & $p{=}2\!:\,[2^{3}]_{3}$ $p{=}3\!:\,1^{1}\,3^{2}$ $p{=}7\!:\,1^{-2}\,7^{1}$ & $1/14$\,{\footnotesize(6)}, $1/12$\,{\footnotesize(4)}, $1/6$\,{\footnotesize(12)}, $1/2$\,{\footnotesize(3)} & 25 & harmonic Maass \\
\chainbandD{$D = 10$}{8 root bases}
\texttt{(1,1,1,5)} & 20 & $p{=}2\!:\,1^{-2}:[4^{1}]_{1}$ $p{=}5\!:\,1^{-2}\,5^{-1}$ & $1/8$\,{\footnotesize(2)}, $1/5$\,{\footnotesize(2)}, $1$\,{\footnotesize(1)} & 5 & harmonic Maass \\
\texttt{(1,1,1,5)} & 100 & $p{=}2\!:\,1^{-2}:[4^{-1}]_{5}$ $p{=}5\!:\,1^{-1}\,5^{-2}$ & $1/40$\,{\footnotesize(12)}, $1/5$\,{\footnotesize(6)}, $1$\,{\footnotesize(1)} & 19 & harmonic Maass \\
\texttt{(1,1,2,4)} & 20 & $p{=}2\!:\,1^{-2}:[4^{1}]_{1}$ $p{=}5\!:\,1^{-2}\,5^{-1}$ & $1/8$\,{\footnotesize(2)}, $1/5$\,{\footnotesize(2)}, $1$\,{\footnotesize(1)} & 5 & harmonic Maass \\
\texttt{(1,1,2,4)} & 40 & $p{=}2\!:\,[2^{-3}]_{1}$ $p{=}5\!:\,1^{-2}\,5^{1}$ & $1/10$\,{\footnotesize(6)}, $1/4$\,{\footnotesize(1)}, $1/2$\,{\footnotesize(3)} & 10 & harmonic Maass \\
\texttt{(1,3,2,3)} & 60 & $p{=}2\!:\,1^{-2}:[4^{-1}]_{3}$ $p{=}3\!:\,1^{-2}\,3^{-1}$ $p{=}5\!:\,1^{-2}\,5^{1}$ & $1/10$\,{\footnotesize(2)}, $1/3$\,{\footnotesize(2)}, $1/2$\,{\footnotesize(1)}, $1$\,{\footnotesize(1)} & 6 & harmonic Maass \\
\texttt{(1,3,2,3)} & 360 & $p{=}2\!:\,[2^{-3}]_{1}$ $p{=}3\!:\,1^{1}\,3^{-2}$ $p{=}5\!:\,1^{-2}\,5^{1}$ & $1/60$\,{\footnotesize(12)}, $1/12$\,{\footnotesize(6)}, $1/6$\,{\footnotesize(6)}, $1/2$\,{\footnotesize(15)} & 39 & harmonic Maass \\
\texttt{(1,3,2,3)} & 540 & $p{=}2\!:\,1^{-2}:[4^{-1}]_{3}$ $p{=}3\!:\,1^{1}\,3^{-1}\,9^{-1}$ $p{=}5\!:\,1^{-2}\,5^{1}$ & $1/90$\,{\footnotesize(12)}, $1/9$\,{\footnotesize(6)}, $1/3$\,{\footnotesize(6)}, $1/2$\,{\footnotesize(3)} & 27 & harmonic Maass \\
\texttt{(1,3,2,3)} & 1080 & $p{=}2\!:\,[2^{3}]_{3}$ $p{=}3\!:\,1^{1}\,3^{1}\,9^{-1}$ $p{=}5\!:\,1^{-2}\,5^{-1}$ & $1/36$\,{\footnotesize(18)}, $1/20$\,{\footnotesize(18)}, $1/6$\,{\footnotesize(18)}, $1/2$\,{\footnotesize(9)} & 63 & harmonic Maass \\
\chainbandD{$D = 15$}{4 root bases}
\texttt{(1,3,1,5)} & 120 & $p{=}2\!:\,1^{-2}:[8^{-1}]_{3}$ $p{=}3\!:\,1^{2}\,3^{-1}$ $p{=}5\!:\,1^{-2}\,5^{-1}$ & $1/48$\,{\footnotesize(4)}, $1/5$\,{\footnotesize(4)}, $1/3$\,{\footnotesize(4)}, $1$\,{\footnotesize(2)} & 14 & harmonic Maass \\
\texttt{(1,3,1,5)} & 360 & $p{=}2\!:\,1^{-2}:[8^{1}]_{1}$ $p{=}3\!:\,1^{-1}\,3^{2}$ $p{=}5\!:\,1^{-2}\,5^{1}$ & $1/16$\,{\footnotesize(2)}, $1/15$\,{\footnotesize(16)}, $1/3$\,{\footnotesize(8)}, $1$\,{\footnotesize(2)} & 28 & harmonic Maass \\
\texttt{(1,3,1,5)} & 600 & $p{=}2\!:\,1^{-2}:[8^{1}]_{7}$ $p{=}3\!:\,1^{2}\,3^{1}$ $p{=}5\!:\,1^{-1}\,5^{-2}$ & $1/240$\,{\footnotesize(24)}, $1/15$\,{\footnotesize(24)}, $1/5$\,{\footnotesize(12)}, $1$\,{\footnotesize(2)} & 62 & harmonic Maass \\
\texttt{(1,3,1,5)} & 1800 & $p{=}2\!:\,1^{-2}:[8^{-1}]_{5}$ $p{=}3\!:\,1^{1}\,3^{2}$ $p{=}5\!:\,1^{1}\,5^{-2}$ & $1/80$\,{\footnotesize(12)}, $1/15$\,{\footnotesize(48)}, $1/5$\,{\footnotesize(12)}, $1/3$\,{\footnotesize(8)} & 80 & harmonic Maass \\
\chainbandD{$D = 22$}{2 root bases}
\texttt{(1,1,1,2,2,4)} & 44 & $p{=}2\!:\,1^{-2}:[4^{1}]_{7}$ $p{=}11\!:\,1^{2}\,11^{-1}$ & $1/11$\,{\footnotesize(2)}, $1/2$\,{\footnotesize(1)}, $1$\,{\footnotesize(1)} & 4 & harmonic Maass \\
\texttt{(1,1,1,2,2,4)} & 88 & $p{=}2\!:\,[2^{-3}]_{7}$ $p{=}11\!:\,1^{2}\,11^{1}$ & $1/22$\,{\footnotesize(6)}, $1/4$\,{\footnotesize(3)}, $1/2$\,{\footnotesize(3)} & 12 & harmonic Maass \\
\end{longtable}\endgroup

\sect{References}{References}
\begin{list}{}{\setlength{\itemindent}{-6mm}\setlength{\leftmargin}{6mm}\setlength{\itemsep}{2pt}}\raggedright
\item[] Allcock, D., \textit{The reflective Lorentzian lattices of rank 3}, Mem. Amer. Math. Soc. 220 (2012), no. 1033; arXiv:1010.0486.
\item[] Borcherds, R. E., \textit{Automorphic forms with singularities on Grassmannians}, Invent. Math. 132 (1998), 491–562.\quad{\footnotesize\color{grey}[singular theta lift, \S{}2.2]}
\item[] Borcherds, R. E., \textit{Reflection groups of Lorentzian lattices}, Duke Math. J. 104 (2000), no. 2, 319–366.\quad{\footnotesize\color{grey}[section 9, Lemma 9.1]}
\item[] Borcherds, R. E., \textit{The Gross–Kohnen–Zagier theorem in higher dimensions}, Duke Math. J. 97 (1999), 219–233; correction, ibid. 105 (2000), 183–184.\quad{\footnotesize\color{grey}[obstruction theorem]}
\item[] Bruinier, J. H. and Funke, J., \textit{On two geometric theta lifts}, Duke Math. J. 125 (2004), 45–90.\quad{\footnotesize\color{grey}[xi operator, pairing]}
\item[] Bruinier, J. H. and Stein, O., \textit{The Weil representation and Hecke operators for vector valued modular forms}, Math. Z. 264 (2010), no. 2, 249–270.\quad{\footnotesize\color{grey}[(2.5), (2.14): the Z-action]}
\item[] Bruinier, J. H., Ehlen, S. and Freitag, E., \textit{Lattices with many Borcherds products}, Math. Comp. 85 (2016), no. 300, 1953–1981.\quad{\footnotesize\color{grey}[simple finite quadratic modules]}
\item[] Brzezinski, J., \textit{A characterization of Gorenstein orders in quaternion algebras}, Math. Scand. 50 (1982), no. 1, 19–24.\quad{\footnotesize\color{grey}[ternary forms <-> orders]}
\item[] Cho, E., \textit{Forced shadows of an obstructed hyperbolic Kac–Moody denominator}, arXiv:2608.19706; DOI 10.5281/zenodo.21973291.\quad{\footnotesize\color{grey}[paper I]}
\item[] Errthum, E., \textit{Singular moduli of Shimura curves}, Canad. J. Math. 63 (2011), 826–861; arXiv:0711.4316.
\item[] Feingold, A. J. and Frenkel, I. B., \textit{A hyperbolic Kac–Moody algebra and the theory of Siegel modular forms of genus 2}, Math. Ann. 263 (1983), 87–144.
\item[] Gritsenko, V. A. and Nikulin, V. V., \textit{Automorphic forms and Lorentzian Kac–Moody algebras I, II}, Internat. J. Math. 9 (1998), 153–199 and 201–275.
\item[] Gross, B. H. and Lucianovic, M. W., \textit{On cubic rings and quaternion rings}, J. Number Theory 129 (2009), 1468–1478.
\item[] Gross, B. H. and Zagier, D., \textit{Heegner points and derivatives of L-series}, Invent. Math. 84 (1986), 225–320.\quad{\footnotesize\color{grey}[derivative regime, section 7]}
\item[] Kac, V. G., \textit{Infinite Dimensional Lie Algebras}, 3rd ed., Cambridge Univ. Press (1990).\quad{\footnotesize\color{grey}[Weyl–Kac denominator, Ch. 10]}
\item[] Kohnen, W. and Zagier, D., \textit{Values of L-series of modular forms at the center of the critical strip}, Invent. Math. 64 (1981), 175–198.
\item[] Kudla, S. S., \textit{Central derivatives of Eisenstein series and height pairings}, Ann. of Math. 146 (1997), 545–646.\quad{\footnotesize\color{grey}[derivative regime, section 7]}
\item[] LMFDB Collaboration, \textit{The L-functions and Modular Forms Database}, https://www.lmfdb.org.\quad{\footnotesize\color{grey}[newform labels]}
\item[] Nikulin, V. V., \textit{Integer symmetric bilinear forms and some of their geometric applications}, Izv. Akad. Nauk SSSR Ser. Mat. 43 (1979), 111–177; English transl., Math. USSR-Izv. 14 (1980), 103–167.\quad{\footnotesize\color{grey}[discriminant forms]}
\item[] Rosson, H. and Tornaría, G., \textit{Central values of quadratic twists for a modular form of weight 4}, in: Ranks of Elliptic Curves and Random Matrix Theory, LMS Lecture Note Ser. 341, Cambridge Univ. Press (2007), 315–321; arXiv:math/0504379.\quad{\footnotesize\color{grey}[nearest precedent, \S{}4.9]}
\item[] Serre, J.-P. and Stark, H., \textit{Modular forms of weight 1/2}, in: Modular Functions of One Variable VI, Lecture Notes in Math. 627, Springer (1977), 27–67.\quad{\footnotesize\color{grey}[weight-1/2 theta series, \S{}5.1]}
\item[] Shimura, G., \textit{On modular forms of half integral weight}, Ann. of Math. 97 (1973), 440–481.\quad{\footnotesize\color{grey}[Shimura lift]}
\item[] Voight, J., \textit{Quaternion Algebras}, Graduate Texts in Math. 288, Springer (2021).\quad{\footnotesize\color{grey}[genus formula, \S{}39.4]}
\item[] Waldspurger, J.-L., \textit{Sur les coefficients de Fourier des formes modulaires de poids demi-entier}, J. Math. Pures Appl. (9) 60 (1981), 375–484.
\item[] Williams, B., \textit{weilrep} (software), https://github.com/btw-47/weilrep, GPL-2.0-or-later, commit e3a7784 (2026-08-27).\quad{\footnotesize\color{grey}[all computation]}
\item[] Zagier, D., \textit{Ramanujan's mock theta functions and their applications}, Sem. Bourbaki 986 (2007/2008).\quad{\footnotesize\color{grey}[term "shadow"]}
\end{list}

\end{document}